\documentstyle[12pt,titlepage,times,diss,amssymb,latexsym,
verbatim,amsmath]{book}

\input epsf.sty

\newcommand{\al}{{\alpha}}

\newcommand{\bbr}{{\Bbb R}}
\newcommand{\bbq}{{\Bbb Q}}
\newcommand{\bbz}{{\Bbb Z}}
\newcommand{\bc}{{\Bbb C}}
\newcommand{\bd}{{\textrm{Bd}\,}}
\newcommand{\be}{\begin{enumerate}}
\newcommand{\bl}{{\rm{Bl}\,}}

\newcommand{\bo}{\partial\,}
\newcommand{\br}{{\Bbb R}}
\newcommand{\bu}{\bullet}

\newcommand{\ca}{{\cal A}}

\newcommand{\cat}{\text{\bf Cat}}
\newcommand{\catr}{{\cal R}}
\newcommand{\cb}{{\cal B}}
\newcommand{\cc}{{\cal C}}
\newcommand{\cd}{{\cal D}}

\newcommand{\cf}{{\cal F}}
\newcommand{\cg}{{\cal G}}

\newcommand{\ck}{{\cal K}}
\newcommand{\cl}{{\cal L}}
\newcommand{\cm}{{\cal M}}
\newcommand{\cn}{{\cal N}}
\newcommand{\co}{{\cal O}}

\newcommand{\cs}{{\cal S}}

\newcommand{\cv}{{\cal V}}
\newcommand{\ctop}{{\bf Top}}

\newcommand{\da}{\Delta}

\newcommand{\dar}{\downarrow}
\newcommand{\Dar}{\Downarrow}

\newcommand{\dc}{{\Bbb C}}

\newcommand{\dg}{\da(G)}
\newcommand{\dgn}{\da(G_n)}

\newcommand{\dr}{{\Bbb R}}
\newcommand{\dq}{{\Bbb Q}}
\newcommand{\dz}{{\Bbb Z}}

\newcommand{\ee}{\end{enumerate}}

\newcommand{\Ga}{\Gamma}
\newcommand{\GL}{\text{\bf GL}}
\newcommand{\grp}{\text{\bf Grp}}

\newcommand{\hra}{\hookrightarrow}

\newcommand{\id}{\text{id}}
\newcommand{\il}{\text{colim}\,}
\newcommand{\im}{\text{Im}\,}

\newcommand{\la}{\leftarrow}
\newcommand{\lam}{\lambda}
\newcommand{\Lam}{\Lambda}
\newcommand{\lar}{\leftarrow}

\newcommand{\lcm}{\text{lcm}\,}

\newcommand{\lla}{\longleftarrow}
\newcommand{\lmu}{(\lam,\mu)}
\newcommand{\Lmu}{(\Lam,\mu)}
\newcommand{\lmt}{\longmapsto}
\newcommand{\lr}{\leftrightarrow}
\newcommand{\lra}{\longrightarrow}
\newcommand{\lspan}{\text{span}\,}

\newcommand{\minus}{\text{Minus\,}}

\newcommand{\nin}{\noindent}

\newcommand{\pl}{\Pi_\lambda}
\newcommand{\plus}{\text{Plus\,}}
\newcommand{\pr}{\nin{\bf Proof. }}

\newcommand{\ra}{\rightarrow}
\newcommand{\Ra}{\Rightarrow}
\newcommand{\rh}{{\widetilde H}}
\newcommand{\rk}{{\text{rk}\,}}

\newcommand{\rsa}{\rightsquigarrow}

\newcommand{\sgn}{\text{sgn}\,}
\newcommand{\sig}{\sigma}
\newcommand{\sla}{\Sigma_\lam^{S^2}}
\newcommand{\slx}{\Sigma_\lambda^X}
\newcommand{\sm}{\setminus}
\newcommand{\Span}{\text{\rm span}}

\newcommand{\ssc}{\text{\bf SS}}
\newcommand{\st}{\text{Stab}\,}
\newcommand{\stab}{\text{Stab}\,}
\newcommand{\susp}{\text{susp}\,}

\newcommand{\tb}{\tilde\beta}

\newcommand{\tda}{\Delta}

\newcommand{\thra}{\twoheadrightarrow}
\newcommand{\ti}{\tilde}
\newcommand{\tp}{\text{\bf Top}}
\newcommand{\tspan}{{\text{span}\,}}

\newcommand{\U}{\text{\bf U}}
\newcommand{\un}{\text{un}\,}

\newcommand{\wti}{\widetilde }
\newcommand{\wsl}{{\wti\Sigma}_\lam}

\newcommand{\xlm}{X_{\lambda,\mu}}
\newcommand{\xLm}{X_{\Lambda,\mu}}

\newcommand{\qed}{$\square$}
\newcommand{\mqed}{\square}

\newtheorem{thm}{Theorem}[section]
\newtheorem{df}[thm]{Definition}
\newtheorem{lm}[thm]{Lemma}
\newtheorem{crl}[thm]{Corollary}
\newtheorem{prop}[thm]{Proposition}
\newtheorem{conj}[thm]{Conjecture}
\newtheorem{rem}[thm]{Remark}

\newtheorem{exam}[thm]{Example}
\newtheorem{exams}[thm]{Examples}
\numberwithin{equation}{section}

\begin{document}

\pagenumbering{roman}
\thispagestyle{empty}

\clearemptydoublepage
\thispagestyle{empty}

\mbox{ }

\vskip2cm

\begin{center}
 {\LARGE{\sc Trends in Topological Combinatorics}} 
\end{center}

\vskip1.2cm
\begin{center} 
{\large {\sc Dmitry N.\ Kozlov}}
\end{center}

\vskip10cm
\begin{center}
  {Habilitationsschrift zur Vorlage an der 
   philosophisch-naturwissenschaftlichen Fakult\"at der 
   Universit\"at Bern \\[1cm]
    Dezember 2001}
\end{center}


\clearemptydoublepage
\thispagestyle{empty}

\vspace{3cm}
\mbox{ }\hfill
\begin{minipage}{7cm}
  { I was the shadow of the waxwing slain}\\
{  By the false azure in the windowpane;}\\
{ I was the smudge of ashen fluff - and I}\\
{ Lived on, flew on, in the reflected sky.}\\[0.2cm]
    -Vladimir Nabokov, {\it Pale Fire}
\end{minipage}


\clearemptydoublepage
\chapter*{Preface}
\addcontentsline{toc}{chapter}{Preface}


\vskip10pt

This thesis is thought to be reflective of the transformation, mistily
rendered in the citation on the previous page, which has taken place
in me since I have handed in my Ph.D.\ thesis in the spring of 1996.
If this transformation, not unnatural taking into account that I was
23 at the time, was to the better or to the worse, I leave for the
reader to judge.

\vskip10pt

\noindent
Now, before indulging into mathematics, I would like to thank my
coauthors Eva-Maria Feichtner and Eric Babson, without whom the
second, resp.\ the fifth, chapters would not exist.

%
%
%
%

\vskip10pt

\noindent
In recent years I was supported in various forms by a number of
research and educational establishments. I therefore express my
graditude to:

\vskip2pt

\noindent
University of Bern,

\noindent
Swiss National Science Foundation, 

\noindent
Forschungsinstitut f\"ur Mathematik at ETH Z\"urich, 

\noindent
The Institute for Advanced Study at Princeton, 

\noindent
Massachusetts Institute of Technology at Cambridge,

\noindent
Mathematical Science Research Institute at Berkeley,
 
\noindent
Royal Institute of Technology at Stockholm,

\noindent
Swedish Natural Science Research Council.

\vskip10pt

\noindent
Finally, I would like to mention again Eva-Maria Feichtner
without whose unrequiring support this work would have been
a mere Gedankenspiel. 

\vskip3pt

\noindent
I dedicate this thesis to her.

\vskip80pt

\noindent
{\it Bern} \hskip255pt  {\sc Dmitry N.\ Kozlov}

\noindent
{\it December 14, 2001}


\clearemptydoublepage
\chapter*{Summary\markboth{Summary}{}}

\addcontentsline{toc}{chapter}{Summary}

This thesis opens with an~introductory discussion, where the reader is
gently led to the world of topological combinatorics, and, where the
results of this Habilitationsschrift are portrayed against the
backdrop of the broader philosophy of the subject. That introduction
is followed by 5 chapters, where the main body of research is
presented, and 4 appendices, where various standard tools and
notations, which we use throughout the text, are collected.

\vskip10pt

The main purpose of the first chapter is to introduce a~new category,
which we call resonance category, whose combinatorics reflects that of
canonical stratifications of $n$-fold symmetric smash products. The
study of the stratifications can then be abstracted to the study of
functors, which we name resonance functors, satisfying a~certain set
of axioms.

One frequently studied stratification is that of the set of all
degree~$n$ polynomials, defined by fixing the allowed multiplicities
of roots. We describe how our abstract combinatorial framework helps
to yield new information on the homology groups of the strata.

\vskip10pt

In the second chapter we introduce notions of combinatorial blowups,
building sets, and nested sets for arbitrary meet-semilattices. This
gives a~common abstract framework for the incidence combinatorics
occurring in the context of De~Concini-Procesi models of subspace
arrangements and resolutions of singularities in toric varieties. Our
main theorem states that a sequence of combinatorial blowups,
prescribed by a building set in linear extension compatible order,
gives the face poset of the corresponding simplicial complex of nested
sets.  As applications we trace the incidence combinatorics through
every step of the De~Concini-Procesi model construction, and we
introduce the notions of building sets and nested sets to the context
of toric varieties.

\vskip10pt

In the third chapter we study rational homology groups of one-point
compactifications of spaces of complex monic polynomials with multiple
roots. These spaces are indexed by number partitions. A~standard
reformulation in terms of quotients of orbit arrangements reduces the
problem to studying certain triangulated spaces~$\xlm$.

\newpage

We present a~combinatorial description of the cell structure of~$\xlm$
using the language of marked forests. This allows us to perform
a~complete combinatorial analysis of the topological properties of
$\xlm$ for several specific cases.  As applications we prove that
$\da(\Pi_n)/\cs_n$ is collapsible, we obtain a~new proof of a~theorem
of Arnol'd, and we find a~counterexample to a~conjecture of Sundaram
and Welker, along with a~few other smaller results.

\vskip10pt

To every directed graph G one can associate a~complex $\da(G)$
consisting of directed subforests. This construction, suggested to us
by R.~Stanley, is especially important in the~case of a~complete
double directed graph $G_n$, where it leads to studying certain
representations of the~symmetric group and relates (via
Stanley-Reisner correspondence) to an~interesting quotient ring. In
the fourth chapter we study complexes $\da(G)$ and associated quotient
complexes.

One of our results states that $\da(G_n)$ is shellable, in particular
Cohen-Ma\-cau\-lay, which can be further translated to say that
the~Stanley-Reisner ring of $\da(G_n)$ is Cohen-Macaulay. Besides
that, by computing the~ho\-mo\-lo\-gy groups of $\da(G)$ for the~cases
when $G$ is essentially a~tree, and when $G$ is a~double directed
cycle, we touch upon the~general question of the~interaction of
combinatorial properties of a~graph and topological properties of
the~associated complex.

We then continue with the study of the $\cs_n$-quotient of the complex
of directed forests, denoted $\da(G_n)/\cs_n$. It is a~simplicial
complex whose cell structure is defined combinatorially. We make use
of the machinery of spectral sequences to analyze these complexes.
In~particular, we prove that the integral homology groups of
$\da(G_n)/\cs_n$ may have torsion.

\vskip10pt

In the fifth chapter we study quotients of posets by group actions. In
order to define the quotient correctly we enlarge the considered class
of categories from posets to loopfree categories: categories without
nontrivial automorphisms and in\-ver\-ses. We view group actions as
certain functors and define the quotients as co\-li\-mits of these
functors. The advantage of this definition over studying the
quo\-ti\-ent poset (which in our language is the colimit in the poset
category) is that the realization of the quotient loopfree category is
more often homeomorphic to the quotient of the realization of the
original poset. We give conditions under which the quotient commutes
with the nerve functor, as well as conditions which guarantee that the
quotient is again a~poset.


\clearemptydoublepage
\chapter*{Zusammenfassung}
\sloppy

\addcontentsline{toc}{chapter}{Zusammenfassung}

Wir beginnen mit einer einleitenden Diskussion, die den Leser in die
Welt der topologischen Kombinatorik einf\"uhrt und die Ergebnisse
dieser Habilitations\-schrift vor den breiteren Hintergrund der das
Gebiet bestimmenden Gedankenwelt stellt. Dieser Einleitung folgen
f\"unf Kapitel, in denen wir die Hauptergebnisse ent\-wickeln, und
vier Appendizes, in denen wir Standardwerkzeuge und Notationen
zusammenstellen, wie sie im Laufe der Arbeit verwendet werden.

\vskip10pt

Ziel des ersten Kapitels ist es, eine neue Kategorie, die
Resonanzkategorie, einzuf\"uhren, deren Kombinatorik diejenige
kanonischer Stratifizierungen $n$-fa\-cher symmetrischer Smash-Produkte
widerspiegelt. Das Studium dieser Stra\-ti\-fi\-zierungen kann so zum
Studium sogenannter Resonanzfunktoren abstrahiert werden, die den
Axiomen eines bestimmten Axiomensystems gen\"ugen.

Eine h\"aufig studierte Stratifikation ist die der Menge von Polynomen
festen Grades~$n$ nach den unter den Wurzeln auftretenden
Vielfachheiten.  Wir zeigen, wie wir dank unseres abstrakten
kombinatorischen Systems zu neuen Informationen \"uber die Homologie
der Strata gelangen.

\vskip10pt

Im zweiten Kapitel f\"uhren wir f\"ur beliebige
Durchschnittshalbverb\"ande die Begriffe kombinatorischer
Aufblasungen, Bausatzmengen und vernetzter Mengen ein. Dies liefert
einen gemeinsamen abstrakten Rahmen f\"ur Inzidenzen, wie sie
einerseits im Zusammenhang mit De Concini-Procesi Modellen von
Arrangements, andererseits bei Aufl\"osungen von Singularit\"aten in
torischen Variet\"aten eine Rolle spielen.  Unser Hauptsatz besagt,
dass eine Sequenz kombinatorischer Aufblasungen, bestimmt durch eine
mittels linearer Erweiterung geordnete Bau\-satz\-menge, auf die
Seitenhalbordnung des zugeh\"origen Simplizialkomplexes vernetzter
Mengen f\"uhrt. Als Anwendung verfolgen wir die Inzidenzen in der
schritt\-wei\-sen Modellkonstruktion nach De Concini und Procesi und
f\"uhren die Begriffe von Bausatzmengen und vernetzten Mengen in den
Zusammenhang torischer Variet\"aten ein.

\vskip10pt

Im dritten Kapitel studieren wir rationale Homologiegruppen der
1-Punkt-Kompaktifizierungen von R\"aumen komplexer normierter Polynome
mit mehr\-fachen Nullstellen. Eine Umformulierung f\"ur Quotienten von
Bahnenarrangements reduziert das Problem auf das Studium gewisser
triangulierter R\"aume $\xlm$.

\newpage

Wir geben eine kombinatorische Beschreibung der Zellstruktur von
$\xlm$, indem wir uns der Terminologie sogenannter markierter B\"aume
bedienen. Dies erlaubt uns eine vollst\"andige Analyse topologischer
Eigenschaften von $\xlm$ in verschiedenen Spezialf\"allen.  Als
Anwendung beweisen wir einerseits die Kollabierbarkeit von
$\da(\Pi_n)/\cs_n$, andererseits erhalten wir, neben anderen
Ergebnissen, einen neuen Beweis f\"ur
einen Satz von Arnol'd und geben ein Gegenbeispiel f\"ur eine
Vermutung von Sundaram und Welker.

\vskip10pt

Jedem gerichteten Graphen~$G$ kann man einen Komplex~$\Delta(G)$
bestehend aus gerichteten Teilw\"aldern zuordnen. Diese Konstruktion,
die auf einen Vorschlag von R.~Stanley zur\"uckgeht, ist besonders
wichtig f\"ur vollst\"andige, doppelt ge\-rich\-te\-te Graphen $G_n$.
F\"ur solche f\"uhrt sie auf das Studium bestimmter Darstellungen der
symmetrischen Gruppe und entspricht via Stanley-Reisner Korrespondenz
interessanten Quotientenringen.  Im vierten Kapitel studieren wir die
Komplexe~$\Delta(G)$ und zugeh\"orige Quotientenkomplexe. Eines
unserer Resultate besagt, dass $\Delta(G)$ sch\"albar ist, damit
insbesondere Cohen-Macaulay, was wiederum \"ubersetzt werden kann in
die Tatsache, dass der Stanley-Reisner-Ring von~$\Delta(G)$ die
Cohen-Macaulay Eigenschaft besitzt. Dar\"uberhinaus r\"uhren wir durch
Homologieberechnungen von~$\Delta(G)$ in F\"allen, in denen $G$ im
wesentlichen ein Baum oder ein doppelt gerichteter Kreis ist, an die
allgemeinere Frage des Zusammenspiels von kombinatorischen
Eigenschaften des Graphen und topologischen Eigenschaften des
zugeh\"origen Komplexes.

Darauffolgend studieren wir den $\cs_n$-Quotienten des Komplexes
gerichteter W\"alder, $\da(G_n)/\cs_n$. Es handelt sich um einen
Simplizialkomplex mit kombinatorisch definierter Zellstruktur. Wir
benutzen Spektralsequenztechniken, um diese Komplexe zu studieren.
Insbesondere zeigen wir, dass Torsion auftreten kann in der
ganzzahligen Homologie von $\da(G_n)/\cs_n$.

\vskip10pt

Im f\"unften Kapitel studieren wir Quotienten von Halbordnungen nach
Gruppenoperationen. Um Quotienten sauber definieren zu k\"onnen,
erweitern wir die betrachtete Klasse von Kategorien von Halbordnungen
zu schleifenfreien Ka\-te\-gorien: Kategorien mit nicht-trivialen
Automorphismen und Inversen. Wir fassen Gruppenaktionen als Funktoren
auf und definieren Quotienten als Ko\-limiten dieser Funktoren. Der
Vorteil dieser Definition gegen\"uber Quotienten\-halb\-ord\-nungen (in
unserer Sprache die Kolimiten in der Kategorie der Halbordnungen) ist,
dass die Realisierungen von Quotienten schleifenfreier Kategorien
h\"aufiger hom\"oomorph sind zu den Quotienten von Realisierungen der
urspr\"unglichen Halb\-ordnung. Wir geben Bedingungen daf\"ur an, dass
die Quotientenbildung mit dem Nerv-Funktor vertauschbar ist, und
formulieren dar\"uberhinaus Bedingungen, die garantieren, dass der
Quotient wiederum eine Halbordnung ist.


\clearemptydoublepage
\tableofcontents

\clearemptydoublepage

\pagenumbering{arabic}

\clearemptydoublepage
\chapter*{Introduction\markboth{Introduction}{Introduction}}
\addcontentsline{toc}{chapter}{Introduction}

\vskip6pt
\begin{center}
{\bf Topological Combinatorics - a chapter of Discrete Mathematics} 
\end{center}
\vskip4pt

\nin For the sake of brevity, the theme of this thesis could be
defined as {\bf discrete mathematics} or {\bf combinatorics} and its
interactions with pure mathematics, in particular with algebraic
topology. This subject has long gone either without a~name, or
subordinated under other titles, which did not fully reflect its
philosophy. Lately, the term {\bf topological combinatorics} has been
used with rising frequency, so, for the matter of presentation, allow
me to extend the use of this name to the present exposition.

Despite some striking applications, see e.g.\ ~\cite{Lo78}, and the
subsequent explosive development, topological combinatorics is still
a~comparatively young subject. Its main problems ({\it
  Fragestellungen}), as well as its methods are still in their forming
stages, although a~lot of consolidating work has been undertaken in
recent years, see e.g.\ ~\cite{Bj95}. Let me highlight here two major
venues of research.

\mbox{ }
\begin{center}
{\bf Combinatorial structures in topology and geometry}
\end{center}
\vskip4pt

\nin Succinctly, the philosophy of the subject might be described as
follows. The first step consists of singling out a~structure, with
combinatorial flavor, occurring naturally in topology or geometry.
Then, one tries to find a suitable, purely combinatorial definition of
this structure. This usually results in combining known, well-studied
objects in a~way, which would be neither very natural, nor apparent
from the original context, in which these objects were defined.

Once this is done, there are two major possibilities to continue.
The~first option is to develop better structural understanding of the
new combinatorial theory, which has just been distiled. That
understanding may then be used to refine our knowledge of the original
object.

The second option is to seek out the occurrences of this combinatorial
object in other topological or geometric situations, where the
connections to the original situation are not clear, and at the moment
the link is only provided on the abstract combinatorial level.
Hopefully, one can then use this new combinatorial connection to
better understand the topological objects in question. A few instances
of the developments in this spirit are: matroid bundles and
combinatorial differential manifolds of MacPherson,
\cite{An99,GeM92,Mac93}, complex matroids, \cite{Zi93}, or, even such
a~general concept as that of oriented matroids, \cite{BLSWZ}.

\vskip4pt

The first part of this thesis can be classified into this category of
research. In the first chapter we study a~circle of problems arising
from topology, with some specific instances appearing in, among
others, singularity theory. The goal is to describe in a~functorial
way the combinatorial structure of canonical stratifications of
$n$-fold symmetric smash products. By canonical stratifications
we mean those implied by point coincidences, with strata being 
indexed by number partitions of~$n$.

The approach, which we choose to achieve that objective, is to
introduce a~new combinatorial object, a~new category, which we call
{\it resonance category}. The objects of this category are not number
partitions themselves, but rather certain equivalence classes
consisting of number partitions, which we call $n$-cuts. The
introduction of these equivalence classes, in a~characteristic for the
ideology of this subject way, is motivated by the topological fact
that strata, which are indexed by number partitions belonging to the
same $n$-cut, are homeomorphic.

Furthermore, substrata of a stratum can be glued so as to obtain a~new
stratum, and also various strata can be included into each other.
To~reflect these topological maps combinatorially, we define the
morphisms of the resonance category to be compositions of two types of
elementary morphisms, called ``gluings'' and ``inclusions,'' which are
defined purely in terms of $n$-cuts.

The original problem of studying stratifications can be viewed as
studying certain functors from the resonance category to the category
of pointed topological spaces, $\text{\bf Top}^*$. This way we split
the problem into two parts: the combinatorial one, encoded by the
resonance category, and the topological one, encoded by such
a~functor.

After that, we go on with the structural study of the resonance
category. Another, somewhat technical notion, which we need to
introduce is that of {\it a~relative resonance}. The main idea is then
to study direct product decomposition properties of relative
resonances. The general theory is rather rich, and we have succeeded
to understand just a~small part of it, but already that was sufficient
to perform various concrete computations.

For example, for a~certain class of resonances we can produce
combinatorial models to explicitly compute algebraic invariants of
the corresponding strata. More specifically, we prove that these
strata are homotopy equivalent to wedges of spheres. Then, we go ahead
to construct weighted directed graphs, whose longest directed paths
enumerate these spheres, and such that the total weights of the
longest paths give the dimensions of the spheres in question. This
provides an~explicit algorithm to compute the homotopy type of
strata for this class of resonances.

\vskip4pt

In the second chapter we turn to a~somewhat different topic.
In~\cite{DP95}, De~Concini \& Procesi have constructed certain models
for complements of subspace arrangements. The idea of their
construction is to turn any given subspace arrangement into a~divisor
with normal crossings, by means of a~sequence of blowups. To use
blowups to resolve singularities is a~standard idea of algebraic
geometry, yet the performance of such a~resolution in specific cases
often leads to intricate constructions.

One aspect of the De~Concini \& Procesi (and of the earlier Fulton \&
MacPherson, \cite{FM94}) construction, which became apparent early on,
was the fact that the combinatorics of the situation was far from
being trivial or standard material. To deal with the mounting
technicalities, De~Concini \& Procesi generalized the notions of
building sets and nested sets, originally coined by Fulton \&
MacPherson in the more specific context of combinatorics of
configuration spaces, to cover the case of subspace arrangements.

In~\cite{DP95} these notions were still defined in geometric terms,
relying on the subspace arrangements themselves, rather than on their
intersection semilattices. In the research, presented in the second
chapter, we generalize the notions of building sets and nested sets to
the purely combinatorial context of semilattices. Moreover, we
introduce a~new combinatorial construction, which we call
a~``combinatorial blowup.''

These three notions can then be combined to yield a lattice-theoretic
theorem (Theorem~\ref{thm_main}), which parallels in statement one of
the main theorems of~\cite{DP95}. As a~result, we obtain a~consistent,
nontrivial combinatorial theory, which canonically depends on the
underlying combinatorial data, and not on the geometric structures
which it encodes, with Theorem~\ref{thm_main} as the main structural
result.

Now, turning back to geometry, we are able to step-by-step trace the
construction of De~Concini \& Procesi, by using the concepts which we
introduced. Furthermore, following the ideas outlined in the beginning
of the discussion, we go on to describe how these combinatorial
notions also appear in the context of toric varieties.

\vskip6pt
\begin{center}{\bf Combinatorially defined algebraic invariants of topological 
    spaces}
\end{center}
\vskip4pt

\nin A~somewhat different direction of research is to work on
determining explicitly algebraic invariants of topological spaces
given by some combinatorial construction. This area has many aspects.
One is to derive general formulae for a~class of objects. An example
of that is provided by the Goresky \& MacPherson formula, see
Appendix~C. Among further instances one could mention the
Orlik-Solomon algebra, which describes the cohomology algebra of the
complement of a~complex hyperplane arrangement, \cite{OS80}, or the
Yuzvinsky basis for the cohomology algebra of De~Concini-Procesi
compactifications, \cite{Yu97}.

\vskip4pt

Another example is described in the third chapter, where we give
a~combinatorial formula to compute the rational Betti numbers of
spaces of complex monic polynomials with roots of fixed
multiplicities. These spaces can be also viewed as strata in the
natural stratification of the $n$-fold symmetric smash product
of~$S^2$, which were previously described in this exposition. To~give
this combinatorial description we need to introduce cell complexes,
which we call $X_{\lambda,\mu}$, which are parameterized by pairs of
number partitions $(\lambda,\mu)$, such that $\lambda\vdash\mu$.

The study of these complexes yields a~variety of applications.
To~start with, we prove collapsibility of the regular CW complex
$\da(\Pi_n)/\cs_n$. The partition lattice $\Pi_n$, and its order
complex $\da(\Pi_n)$, are of fundamental importance for encoding the
combinatorics of configuration spaces, and for computing cohomology of
the colored braid group, which motivates the appearance of the complex
$\da(\Pi_n)/\cs_n$ as a~basic object of study.

We also obtain a~new, conceptual proof of the Arnol'd Finiteness
Theorem, which boils down to a~combinatorial fact about marked
forests. Furthermore, we are able to furnish a~counterexample to
a~conjecture of Sundaram \& Welker concerning the multiplicity of the
trivial character in certain $\cs_n$-representations, as well as to
compute this multiplicity for $\cs_n$-representations previously
studied by R.\ Stanley, \cite{St82}, and P.\ Hanlon, \cite{Ha83}.

\vskip4pt

In the fourth chapter we turn to a~related aspect of this area,
namely, the combinatorial analysis of simplicial complexes with
explicit combinatorial description of the simplices. Such study
facilitates the richness of the theory, and provides a~reasonable
testing field for the computational methods of topological
combinatorics. Some instances of such investigations are various forms
of lexicographic shellability, developed by Bj\"orner et al.\ 
\cite{Bj80,Bj95,BWa83,BWa96,BWa97,Ko97}, and the study of topological
properties of complexes of not $k$-connected graphs, with implications
in knot theory, \cite{BBLSW}.

In our case, the focus is on studying complexes of directed subforests
of a~fixed graph, a~problem suggested to us by R.\ Stanley,
\cite{St97}. For the complete graph $G_n$ we are able to prove that
$\da(G_n)$ is shellable, implying that it is homotopy equivalent to
a~wedge of $(n-1)^{n-1}$ spheres of dimension $n-2$. Furthermore, we
develop machinery to analyze topological properties of $\da(G)$ for
a~class of directed graphs, called essentially trees, and perform
complete computations in the special cases of cycles and double
directed strings.

In the last section of the fourth chapter, we study the quotient
complexes $\da(G_n)/\cs_n$. These turn out to be rather complex and,
after giving the combinatorial description of the simplicial structure
of $\da(G_n)/\cs_n$ in terms of marked forests, we perform their
analysis by means of the apparatus of spectral sequences. In
particular, we prove that the homology groups of $\da(G_n)/\cs_n$ with
integer coefficients are {\it not} torsion-free in general.

\vskip4pt

Finally, in the fifth chapter, we take a~look at various quotient
constructions which appear in the context of group actions on
posets. The quotients appear in the third and the fourth chapters, so
a~short and abstract study of these constructions is commanded by the
natural yearning for completeness.

Since an~action of a~group $G$ is simply a~functor from the
one-object-category associated to $G$, it is natural to ask the
quotient to be the colimit of this functor.  The three main questions
of study, which we concentrate on, are: when are the morphisms of
$P/G$ given by the $G$-orbits of the morphisms of $P$ (a~property
called ``regularity''), when does taking the quotient commute with
Quillen's nerve functor $\da$, and, finally, which classes of
categories we generally may get as quotients.

We are able to give simple combinatorial conditions on the group
actions, which are equivalent to the desired properties. These
conditions are, in some sense, local, which implies that they are
frequently verifiable in specific situations, boding well for them
being useful in the future.


\clearemptydoublepage
\chapter*{PART I \\[2cm] Combinatorial Structures in \\[0.1cm] 
Topology and Geometry}
\addcontentsline{toc}{chapter}{PART I. Combinatorial Structures in Topology 
and Geometry}
\vspace{9cm}
\mbox{ }\hfill
\begin{minipage}{7.3cm}
  { Now entertain conjecture of a time}\\
{  When creeping murmur and the poring dark}\\
{ Fills the wide vessel of the universe.}\\[0.2cm]
    -William Shakespeare, {\it Henry V}
\end{minipage}

\clearemptydoublepage
\chapter{The Resonance Category}

\section{Canonical stratifications of symmetric smash products}

\subsection{Combinatorial stratifications of topological spaces}
\label{s1.1}
              
 Complicated combinatorial problems often arise when one studies the
homological properties of strata in some topological space with
a~given natural stratification. The examples of such stratifications
are numerous. A~very simple one is provided by taking the~$n$-fold
direct product of a topological space (possibly also taking the
quotient with respect to the $\cs_n$-action), stratified by point 
coincidences. The strata are indexed by set partitions (or number
partitions), and the biggest open stratum is a~configuration
space (ordered or unordered), whose topological properties have been 
widely studied, see e.g., \cite{FH01,FZ00}.

Another example is the stratification of a~vector space induced by
a~subspace arrangement. The strata are all possible intersections of
subspaces, they are indexed by the intersection semilattice of the
arrangement. The biggest open stratum is the complement of the
subspace arrangement, whose topological properties have also been of
quite some interest, see
e.g., \cite{Bj94,Bj95,GoM88,OS80,Vas94,Zi92,ZZ93}.

Of course, in both examples above, the main objective is to study the
open stratum, which is the complement of the largest closed stratum.
However, it was suggested by Arnol'd in a~much more general context,
see for example~\cite{Ar70a}, that in situations of this kind one should
study the problem for all closed strata. The main argument in support
of this point of view is that there is usually no immediate natural
structure on the largest open stratum, while there is one on its
complement, also known as discriminant. The structure is simply given
by stratification. To put it in philosophical terms: ``There is only
one way for the point in the stratified space to be good, but there
are many different ways for it to be bad''. Once some information has
been obtained about the closed strata, one can try to find out
something about the open stratum by means of some kind of duality.

After this general introduction we would like to describe the specific
example which will be of particular importance for this chapter. Let
$X$ be a~pointed topological space (we refer to the chosen point as
a~point at $\infty$), and denote
$$X^{(n)}=\overbrace{X\wedge X\wedge\dots\wedge X}^n/\cs_n,$$
where $\wedge$ is the smash product of pointed spaces.
In other words, $X^{(n)}$ is the set of all unordered collections of
$n$ points on $X$ with the infinity point attached in the appropriate
way. $X^{(n)}$ is naturally stratified by point coincidences, and
the strata are indexed by the number partitions of $n$. Note that 
we consider the closed strata, so, for example, the stratum indexed 
$(\underbrace{1,1,\dots,1}_n)$ is the whole space $X^{(n)}$.

If one specifies $X=S^1$, resp.~$X=S^2$, one obtains as strata the
spaces of all monic real hyperbolic, resp.~monic complex, polynomials
of degree~n with specified root multiplicities. These spaces naturally
appear in singularity theory, \cite{AGV85}. Homological invariants of
several of these strata were in particular computed by Arnol'd,
Shapiro, Sundaram, Welker, Vassiliev, and the author,
see~\cite{Ar70a,Ko99a,Ko00b,ShW98,SuW97,Vas98}.
  
\subsection{The idea of the resonance category and resonance functors}

The purpose of the research presented in this chapter is to take a
different, more abstract look at this set of problems. More
specifically, the idea is to introduce a~new canonical combinatorial
object, independent of the topology of the particular space $X$, where
the combinatorial aspects of these stratifications would be fully
reflected. This object is a~certain category, which we name {\bf
  resonance category}. The word resonance stands for certain linear
identities valid among parts of the indexing number partition for the
particular stratum. The usage of this word was suggested to the author
by B.~Shapiro,~\cite{Sh00}.

Having this canonically combinatorially defined category at hand, one
then can, for each specific topological space $X$, view the natural
stratification of $X^{(n)}$ as a~certain functor from the resonance
category to $\ctop^*$. These functors satisfy a~system of axioms,
which we take as a~definition of {\bf resonance functors}.  The
combinatorial structures in the resonance category will then project
to the corresponding structures in each specific $X^{(n)}$. This opens
the door to develop the general combinatorial theory of the resonance
category, and then prove facts valid for all resonance functors
satisfying some further conditions, such as for example acyclicity of
certain spaces.

As the main technical tools to unearth the combinatorial structures in
the resonance category, we put forward the notions of relative
resonances and direct products (most importantly of a~resonance and
a~relative resonance). Intuitively one can think of the relative
resonance as a~stratum with a~substratum shrunk to the infinity
point.

As mentioned above, to illustrate a~possible appearance of this
abstract framework we choose to use a~class of topological spaces
which come in particular from singularity theory, and whose
topological properties have been studied: spaces of polynomials (real
or complex) with prescribed root multiplicities. In particular, in
case of strata $(k^m,1^t)$, which were studied in
\cite{Ar70a,Ko99a,SuW97} for the complex case, and in
\cite{Ko00b,ShW98} for the real case, we demonstrate how the inherent
combinatorial structure of the resonance category makes this
particular resonance especially ``reducible.''

\vskip4pt
\nin
Here is the brief summary of the contents of this chapter. 

\vskip4pt
\nin
{\bf Section \ref{s1.2}.} We introduce the notion of resonance category, 
and describe the structure of its set of morphisms.

\vskip4pt
\nin
{\bf Section \ref{s1.3}.} We introduce the notions of relative resonances,
direct products of relative resonances, and resonance functors.

\vskip4pt
\nin {\bf Section \ref{s1.4}.} We formulate the problem of Arnol'd and Shapiro
which motivated this research as that concerning a~specific resonance
functor. Then, we analyze the combinatorial structure of resonances
$(a^k,b^l)$, which leads to the complete determination of homotopy
types of the corresponding strata for $X=S^1$.

\vskip4pt \nin {\bf Section \ref{s1.5}.} We analyze the combinatorial
structure of sequential resonances. For $X=S^1$, this leads to
a~complete computation of homotopy types of the strata corresponding to
resonances $(a^k,b^l,1^m)$, such that $a-bl\leq m$, as well as
resonances consisting of powers of some number. In the case of the
latter, the strata always have the homotopy type of a bouquet of
spheres. We describe a~combinatorial model to enumerate these spheres
as paths in a certain weighted directed graph, with dimensions of the
spheres being given by the total weights of the paths.

\vskip4pt
\nin
{\bf Section \ref{s1.6}.} We introduce the notion of a complexity of a~resonance
and give a~series of examples of resonances having arbitrarily high
complexity.

         \section{The resonance category} \label{s1.2}

\subsection{Resonances and their symbolic notation} 


    For every positive integer $n$, let $\{-1,0,1\}^n$ denote the set of 
all points in $\br^n$ with coordinates in the set $\{-1,0,1\}$. 
We say that a subset $S\subseteq\{-1,0,1\}^n$ is {\bf span-closed} if 
$\tspan(S)\cap\{-1,0,1\}^n=S$, where $\tspan(S)$ is the linear subspace 
spanned by the origin and points in $S$. Of course the origin lies 
in every span-closed set. For $x=(x_1,\dots,x_n)\in\{-1,0,1\}^n$,
we use the notations $\plus(x)=\{i\in[n]\,|\,x_i=1\}$ and 
$\minus(x)=\{i\in[n]\,|\,x_i=-1\}$. 


\begin{df} $\,$
  
  \nin (1) A subset $S\subseteq\{-1,0,1\}^n$ is called an~{\bf
    $n$-cut} if it is span-closed and for every $x\in
  S\setminus\{$origin$\}$ we have $\plus(x)\neq\emptyset$ and
  $\minus(x)\neq\emptyset$.  We denote the set of all
  $n$-cuts by $\catr_n$.

\nin
(2) $\cs_n$ acts on $\{-1,0,1\}^n$ by permuting coordinates, which in turn
induces $\cs_n$-action on $\catr_n$. The {\bf $n$-resonances} are defined
to be the orbits of the latter $\cs_n$-action. We let $[S]$ denote the
$n$-resonance represented by the $n$-cut $S$.
\end{df}

The cut or the resonance consisting of origin only is called {\it trivial}.

\begin{exams}\label{exsmres}
{\rm $n$-resonances for small values of $n$.}

(1) There are no nontrivial 1-resonances.

(2) There is one nontrivial 2-resonance: $[\{(0,0),(1,-1),(-1,1)\}]$.

(3) There are four nontrivial 3-resonances:
  $$[\{(0,0,0),(1,-1,0),(-1,1,0)\}],$$
  $$[\{(0,0,0),(1,-1,0),(-1,1,0),(1,0,-1),(-1,0,1),(0,1,-1),(0,-1,1)\}],$$
  $$[\{(0,0,0),(1,-1,-1),(-1,1,1)\}],$$
  $$[\{(0,0,0),(1,-1,-1),(-1,1,1),(0,1,-1),(0,-1,1)\}].$$ 

(4) Here is an example of a nontrivial 6-resonance:
$$[\{\text{origin},\pm(1,1,0,-1,-1,0),\pm(0,1,1,0,-1,-1),\pm(1,0,-1,-1,0,1)\}].$$
\end{exams}

\nin
{\bf Symbolic notation.} To describe an $n$-resonance, rather than to list
all of the elements of one of its representatives, it is more convenient
to use the following symbolic notation: we write a sequence of $n$ linear
expressions in some number (between 1 and $n$) of parameters, the order
in which the expressions are written is inessential.

Here is how to get from such a~symbolic expression to the
$n$-resonance: choose an order on the $n$ linear expressions and
observe that now they parameterize some linear subspace of ${\Bbb R}^n$, 
which we denote by $A$. The $n$-resonance is now the orbit of
$A^{\perp}\cap\{-1,0,1\}^n$.

Reversely, to go from an $n$-resonance to a~symbolic expression: choose
a~representative $n$-cut $S$, the symbolic expression can now 
be obtained as a~linear parameterization of $\tspan(S)^\perp$.

For example the 6 nontrivial resonances listed in the
Example~\ref{exsmres} are (in the same order):
$$(a,a),\,(a,a,b),\,(a,a,a),\,(a+b,a,b),\,(2a,a,a),\,
(a+b,b+c,a+d,b+d,c+d,2d).$$

\subsection{Acting on cuts with ordered set partitions} 

From now on we assume known the terminology and notations of set
partitions and ordered set partitions, as described in the Appendix~A.

\begin{df} 
  Given $\pi=(\pi_1,\dots,\pi_k)$ an ordered set partition of $[m]$
  with $k$ parts, and $\nu=(\nu_1,\dots,\nu_m)$ an ordered set
  partition of $[n]$ with $m$ parts, their {\bf composition}
  $\pi\circ\nu$ is an ordered set partition of $[n]$ with $k$ parts,
  defined by $\pi\circ\nu=(\mu_1,\dots,\mu_k)$,
  $\mu_i=\cup_{j\in\pi_i}\nu_j$, for $i=1,\dots,k$.
  
  \nin Analogously, we can define $\pi\circ\nu$ for an~ordered set
  partition $\nu$ and a~set partition $\pi$, in which case
  $\pi\circ\nu$ is a~set partition without any specified order on the
  blocks.
\end{df}

In particular, when $m=n$, and $|\pi_i|=1$, for $i=1,\dots,n$, we can
identify $\pi=(\pi_1,\dots,\pi_n)$ with the corresponding permutation
of $[n]$. The composition of two such ordered set partitions
corresponds to the multiplication of corresponding permutations, and
we denote the ordered set partition $(\{1\},\dots,\{n\})$ by $\id_n$,
or just~$\id$.

\begin{df}
  For $A\subseteq B$, let $p_{B,A}:P(B)\ra P(A)$ denote map induced by
  the restriction from $B$ to $A$. For two disjoint set $A$ and $B$,
  and $\Pi\subseteq P(A)$, $\Lambda\subseteq P(B)$, we define
  $\Pi\times\Lambda=\{\pi\in P(A\cup B)\,|\,p_{A\cup B,A}(\pi)\in\Pi,
  p_{A\cup B,B}(\pi)\in\Lambda\}$.

\end{df}

The following definition provides the combinatorial constructions
necessary to describe the morphisms of the resonance category, as well
as to define the relative resonances. 


\begin{df} \label{df2.4}
Assume $S$ is an~$n$-cut. 

\nin (1) For an~ordered set partition of $[n]$,
$\pi=(\pi_1,\dots,\pi_m)$, we define $\pi S\in\catr_m$ to be the set
of all $m$-tuples $(t_1,\dots,t_m)\in\{-1,0,1\}^m$, for which there
exists $(s_1,\dots,s_n)\in S$, such that for all $j\in[m]$, and
$i\in\pi_j$, we have $s_i=t_j$.

\nin (2) For an~unordered set partition of $[n]$,
$\pi=(\pi_1,\dots,\pi_m)$, we define $S^\pi\in\catr_n$ to be the
subset of $S$ consisting of all $(s_1,\dots,s_n)$ such that if
$i,j\in\pi_k$, $i\neq j$, for some $k=1,\dots,m$, then $s_i=s_j\neq
0$; in other words if $k=1,\dots,m$ is such that $|\pi_k|\geq 2$, then
either $\pi_k\subseteq\plus(s_1,\dots,s_n)$ or
$\pi_k\subseteq\minus(s_1,\dots,s_n)$.
\end{df}

\nin Clearly $\id S=S$, and one can see that $(\pi\circ\nu)S=\pi(\nu S)$.

\vskip4pt

\nin {\bf Verification of $(\pi\circ\nu) S=\pi(\nu S)$.}

\nin By definition we have
$$ (\pi\circ\nu)S=\{(t_1,\dots,t_k)\,|\,\exists(s_1,\dots,s_n)\in S
\text{ s.t. }\forall j\in[k],i\in\mu_j:s_i=t_j\},$$
$$\nu S=\{(x_1,\dots,x_m)\,|\,\exists(s_1,\dots,s_n)\in S
\text{ s.t. }\forall q\in[m],i\in\nu_q:s_i=x_q\},$$
\begin{multline*}
\pi(\nu S)=\{(t_1,\dots,t_k)\,|\,\exists(x_1,\dots,x_m)\in\nu S, \\
\text{ such that }\forall j\in[k],q\in\pi_j,i\in\nu_q:s_i=t_j\}.
\end{multline*}
The identity $(\pi\circ\nu) S=\pi(\nu S)$ follows now from 
the equality $\mu_j=\cup_{q\in\pi_j}\nu_q$.

\vskip4pt

There are many different ways to formulate the Definition~\ref{df2.4}.
We chose the ad hoc combinatorial language, but it is also possible to
put it in the linear-algebraic terms. An ordered set partition of
$[n]$, $\pi=(\pi_1,\dots,\pi_m)$, defines an inclusion map
$\phi:\br^m\ra\br^n$ by $\phi(e_i)= \sum_{j\in\pi_i}\tilde e_j$, where
$\{e_1,\dots,e_m\}$, resp.\ $\{\tilde e_1,\dots,\tilde e_n\}$, is the
standard orthonormal basis of $\br^m$, resp.~$\br^n$. Given
$S\in\catr_n$, $\pi S$ can then be defined as
$\phi^{-1}(\text{Im}\,\phi\cap S)$. Furthermore,
$S^\pi=\text{Im}\,\phi|_{Z}\cap S$, where $Z$ is the set of all
$(z_1,\dots,z_m)\in\br^m$, such that if $|\pi_k|\geq 2$, for some
$k=1,\dots,m$, then $z_k\neq 0$.


\subsection{The definition of the resonance category and 
the terminology for its morphisms} 


We are now ready to give the definition of the central notion of this
chapter.

\begin{df}
  The {\bf resonance category}, denoted $\catr$, is defined as follows:
  
  \nin (1) The set of objects is the set of all $n$-cuts, for
  all positive integers $n$, $\co(\catr)=\cup_{n=1}^\infty \catr_n$.
  
  \nin (2) The set of morphisms is indexed by triples $(S,T,\pi)$,
  where $S\in\catr_m$, $T\in\catr_n$, and $\pi$ is an ordered set
  partition of $[n]$ with $m$ parts, such that $S\subseteq \pi T$.
  For the reasons which will become clear later we denote the morphism
  indexed with $(S,T,\pi)$ by $S\thra\pi T\stackrel{\pi}{\hra} T$.

\nin
As the notation suggests, the initial object of the morphism 
$S\thra\pi T\stackrel{\pi}{\hra} T$ is $S$ and terminal object is $T$.
The composition rule is defined by
$$(S\thra\pi T\stackrel{\pi}{\hra}T)\circ(T\thra\nu
Q\stackrel{\nu}{\hra}Q)= S\thra\pi\nu Q\stackrel{\pi\nu}{\hra} Q,$$
where $S\in\catr_k$, $T\in\catr_m$, $Q\in\catr_n$, $\pi$ is an ordered
set partition of $[m]$ with $k$ parts, and $\nu$ is an ordered set
partition of $[n]$ with $m$ parts.
\end{df}

An~alert reader will notice that the resonances themselves did not
appear explicitly in the definition of the resonance category. In
fact, it is not difficult to notice that resonances are isomorphism
classes of objects of $\catr$. Let us now look at the set of morphisms
of $\catr$ in some more detail.

\nin
(1) For $S\in\catr_n$, the identity morphism of $S$ is
 $S\thra S\stackrel{\id}{\hra} S$. 

\nin
(2) Let us introduce short hand notations: $S\thra T$ for
$S\thra T\stackrel{\id}{\hra} T$, and $\pi T\stackrel{\pi}{\hra}T$
for $\pi T\thra \pi T\stackrel{\pi}{\hra}T$. Then we have
$$S\thra\pi T\stackrel{\pi}{\hra} T=
(S\thra\pi T)\circ(\pi T\stackrel{\pi}{\hra}T).$$
Note also that $S\thra S=S\stackrel{\id}{\hra}S$.

\nin
(3) The associativity of the composition rule can be derived from 
the commutation relation 
$$(\pi S\stackrel{\pi}{\hra}S)\circ(S\thra T)=
(\pi S\thra\pi T)\circ(\pi T\stackrel{\pi}{\hra} T)$$
as follows:
\begin{multline} \notag
(S\thra\pi T\hra T)\circ(T\thra\nu Q\hra Q)\circ(Q\thra\rho X\hra X)=\\
(S\thra\pi T)\circ(\pi T\hra T)\circ(T\thra\nu Q)\circ(\nu Q\hra Q)
\circ(Q\thra\rho X)\circ(\rho X\hra X)=\\
(S\thra\pi T)\circ(\pi T\thra\pi\nu Q)\circ(\pi\nu Q\thra\pi\nu\rho X)\circ \\
(\pi\nu\rho X\hra \nu\rho X)\circ(\nu\rho X\hra \rho X)\circ(\rho X\hra X).
\end{multline}

\nin
(4) We shall use the following names: morphisms $S\thra T$ are called
{\bf gluings} (or $n$-gluings, if it is specified that $S,T\in\catr_n$); 
morphisms $\pi T\stackrel{\pi}{\hra}T$ are called {\bf inclusions}
(or $(n,m)$-inclusions, if it is specified that $T\in\catr_n$, 
$\pi T\in\catr_m$), the inclusions are called {\bf symmetries} if $\pi$ 
is a~permutation. As observed above, the symmetries are the only 
isomorphisms in $\catr$. Here are two examples of inclusions:
$$\{(0,0),(1,-1),(-1,1)\}\stackrel{(\{1\},\{2,3\})}{\hookrightarrow}
\{(0,0,0),(1,-1,-1),(-1,1,1)\},$$
$$\{(0,0),(1,-1),(-1,1)\}\stackrel{(\{1\},\{2,3\})}{\hookrightarrow}
\{(0,0,0),\pm(1,-1,-1),\pm(0,1,-1)\}.$$

\section{Structures related to the resonance category}
\label{s1.3}

\subsection{Relative resonances}

Let $A(n)$ denote the set of all collections of non-empty multisubsets
of~$[n]$, and let $P(n)\subseteq A(n)$ be the set of all partitions
of~$[n]$. For every $S\in\catr_n$ let us define a~closure operation on
$A(n)$, resp.\ on $P(n)$.

\begin{df}\label{df3.1}
Let $\ca\in A(n)$. We define $\ca\Dar S\subseteq A(n)$
to be the minimal set satisfying the following conditions:
\begin{enumerate}
\item[(1)] $\ca\in\ca\Dar S$;
\item[(2)] if $\{B_1,B_2,\dots,B_m\}\in\ca\Dar S$, then 
$\{B_1\cup B_2,B_3,\dots,B_m\}\in\ca\Dar S$;
\item[(3)] if $\{B_1,B_2,\dots,B_m\}\in\ca\Dar S$, and there exists 
$x\in S$, such that $\plus(x)\subseteq B_1$, then
$\{(B_1\setminus\plus(x))\cup\minus(x),B_2,\dots,B_m\}\in\ca\Dar S$.
\end{enumerate}

For $\pi\in P(n)$, we define $\pi\dar S\subseteq P(n)$ as 
$\pi\dar S=(\pi\Dar S)\cap P(n)$. For a set $\Pi\subseteq P(n)$
we define $\Pi\dar S=\cup_{\pi\in\Pi}\pi\dar S$. We say that $\Pi$
is {\bf $S$-closed} if $\Pi\dar S=\Pi$.
\end{df}

The idea behind this definition comes from the context of the standard
stratification of the $n$-fold symmetric product. Given a~stratum $X$
indexed by a number partition of $n$ with $m$ parts, let us fix some
order on the parts. A~substratum $Y$ is obtained by choosing some
partition $\pi$ of $[m]$ and summing the numbers within the blocks of
$\pi$. Since the order of the parts of the number partition indexing
$X$ is fixed, $X$ gives rise to a~unique $m$-cut $S$. The set
$\pi\dar S$ describes all partitions $\nu$ of $[m]$ such that if the
numbers within the blocks of $\nu$ are summed then the obtained
stratum $Z$ satisfies $Z\subseteq Y$. In particular, if $Y$ is shrunk
to a~point, then so is $Z$. The two following examples illustrate how
the different parts of the Definition~\ref{df3.1} might be needed.

\begin{exam} {\rm The equivalences of type (2) from
  the~Definition~\ref{df3.1} are needed.} Let the stratum $X$ be
indexed by $(3,2,1,1,1)$ (fix this order of the parts), and let
$\pi=\{1\}\{23\}\{4\}\{5\}$. Then, the stratum~$Y$ is indexed by
$(3,3,1,1)$. Clearly, the stratum~$Z$, which is indexed by $(3,3,2)$,
lies inside $Y$, hence $\{1\}\{2\}\{345\}\in\pi\dar S$, where $S$ is
the cut corresponding to $(3,2,1,1,1)$. However, if one
starts from the partition $\pi$ and uses equivalences of type (3) from
the Definition~\ref{df3.1}, the only other partitions one can obtain
are $\{1\}\{24\}\{3\}\{5\}$, and $\{1\}\{25\}\{3\}\{4\}$. None of them
refines $\{1\}\{2\}\{345\}$, hence it would not be enough in the
Definition~\ref{df3.1} to just take the partitions which can be
obtained via the equivalences of type (3) and then take $\pi\dar S$ to
be the set of all the partitions which are refined by these.
\end{exam}

\begin{exam} \label{exam1.3.3}
{\rm It is necessary to view the equivalence relation
  on the larger set $A(n)$.} This time, let the stratum $X$ be indexed
by $(a+b,b+c,a+d,b+d,c+d,2d)$ (fix this order of the parts, and assume
as usual that there are no linear relations on the parts other than
those induced by the algebraic identities on the variables $a$, $b$,
$c$, and $d$). Furthermore, let $\pi=\{16\}\{23\}\{45\}$. Then the
stratum $Y$ is indexed by $(a+b+2d,a+b+c+d,b+c+2d)$. Clearly, we have
$\{34\}\{15\}\{26\}\in\pi\dar S$, where $S$ is the cut
corresponding to $(a+b,b+c,a+d,b+d,c+d,2d)$. 
\end{exam}

\noindent
A natural idea for the
Definition~\ref{df3.1} could have been to define the equivalence
relation directly on the set $P(n)$ and use ``swaps'' instead of the
equivalences of type~(3), i.e., to replace the condition (3) by:

\begin{quote} {\it
  if $\{B_1,B_2,\dots,B_m\}\in\ca\dar S$, and there exists $x\in S$,
  such that $\plus(x)\!\subseteq\! B_1$, and $\minus(x)\!\subseteq\! B_2$,
  then $\{(B_1\setminus\plus(x))\,\cup\,\minus(x),$
  $(B_2\setminus\minus(x))\cup\plus(x),B_3,\dots,B_m\}\in\ca\dar S$. }
\end{quote}

\noindent
However, this would not have been sufficient as the
Example~\ref{exam1.3.3} shows, since no swaps would be possible on
$\pi=\{16\}\{23\}\{45\}$.


\begin{df}
  Let $S$ be an~$n$-cut, $\Pi\subseteq P(n)$ an~$S$-closed
  set of partitions. We define
  \[ S\setminus\Pi=S\setminus\bigcup_{\pi\in\Pi}S^\pi.\]
\end{df}

In the next definition we give a~combinatorial analog of viewing
a~stratum relative to a~substratum.

\begin{df} \label{dfrel} $\,$
  
  \nin (1) A {\bf relative $n$-cut} is a~pair $(S,\Pi)$,
where $S\subseteq\{-1,0,1\}^n$, $\Pi\subseteq P(n)$, such that the 
following two conditions are satisfied:
\begin{itemize}
\item $(\lspan S)\setminus\Pi=S$;
\item $\Pi$ is $(\lspan S)$-closed.
\end{itemize}
  
\nin (2) The permutation $\cs_n$-action on $\{-1,0,1\}^n$ induces
an~$\cs_n$-action on the relative $n$-cuts by
$(S,\Pi)\stackrel{\sigma}{\mapsto}(\sigma S,\Pi\sigma^{-1})$, for
$\sigma\in\cs_n$.  The {\bf relative $n$-resonances} are defined to be
the orbits of this $\cs_n$-action. We let $[S,\Pi]$ denote the
relative $n$-resonance represented by the relative $n$-cut~$(S,\Pi)$.
\end{df}

When $S\in\catr_n$ and $\Pi\subseteq P(n)$, $\Pi$ is $S$-closed, it is
convenient to use the notation $Q(S,\Pi)$ to denote the relative cut
$(S\setminus\Pi,\Pi)$. Clearly we have $(S,\Pi)=Q(\tspan S,\Pi)$.
Analogously, $[Q(S,\Pi)]$ denotes the relative resonance
$[S\setminus\Pi,\Pi]$. We use these two notations interchangeably
depending on which one is more natural in the current context.

The special case of the particular importance for our computations in
the later sections is that of $Q(S,\pi\dar S)$, where $\pi$ is
a~partition of $[n]$ with $m$ parts. In this case, we call
$(S\setminus(\pi\dar S),\pi\dar S)$ the relative $(n,m)$-cut
associated to $S$ and $\pi$.

By the Definition~\ref{dfrel}, the relative cut $(S,\Pi)=((\lspan
S)\setminus\Pi,\Pi)$ consists of two parts. We intuitively think of
$(\lspan S)\setminus\Pi$ as the set of all resonances which survive
the shrinking of the strata associated to the elements of $\Pi$, so it
is natural to call them {\it surviving elements}. We also think of
$\Pi$ as the set of all partitions whose associated strata are shrunk
to the infinity point, so, accordingly, we call them {\it partitions
  at infinity}.

\subsection{Direct products of relative resonances} 

\begin{df} \label{df3.3}$\,$
For relative resonances $(S,\Pi)$ and $(T,\Lambda)$ we define
$$(S,\Pi)\times (T,\Lambda)=
(S\times T,(\Pi\times P(m))\cup (P(n)\times\Lambda)).$$
Clearly the orbit $[(S,\Pi)\times (T,\Lambda)]$ does not depend
on the choice of representatives of the orbits $[S,\Pi]$ and
$[T,\Lambda]$, so we may define $[S,\Pi]\times[T,\Lambda]$ to be
$[(S,\Pi)\times (T,\Lambda)]$.
\end{df}

\nin
The following special cases are of particular importance for our
computation:

\vskip4pt

\nin
(1) {\bf A direct product of two resonances.}

\nin For an~$m$-cut~$S$, and an~$n$-cut~$T$, we have 
$S\times T=\{(x_1,\dots,x_m,y_1,\dots,y_n)\,|$ $(x_1,\dots,x_m)\in S,
(y_1,\dots,y_n)\in T\}\in\catr_{m+n}$, 
and $[S]\times[T]=[S\times T]$.

\vskip4pt

\nin
(2) {\bf A direct product of a relative resonance and a resonance.}

\nin For $S\in\catr_n$, $\Pi\subseteq P(n)$ an $S$-closed set of
partitions, and $T\in\catr_k$, we have $Q(S,\Pi)\times T=Q(S\times
T,\widetilde\Pi)$, where $\widetilde\Pi=\Pi\times
P(\{n+1,n+2,\dots,n+k\})$, and $[Q(S,\Pi)]\times [T]=[Q(S,\Pi)\times T]$.

\begin{exam}
\begin{multline*}
Q(\{(0,0,0),\pm(1,-1,-1),\pm(0,1,-1)\},\{1\}\{23\})=\\
\{(0)\}\times Q(\{(0,0),\pm(1,-1)\},\{12\}).
\end{multline*}
\end{exam}

\begin{rem}
One can define a~category, called {\bf relative resonance category}, 
whose set of objects is the set of all relative $n$-cuts. 
A~new structure which it has in comparison to~$\catr$ is provided by 
``shrinking morphisms'': $(S,\Pi)\rsa(T,\Lambda)$, for 
$S,T\subseteq\{-1,0,1\}^n$, $P(n)\supseteq\Lambda\supseteq\Pi$, such 
that $(\lspan S)\setminus\Lambda=T$. They correspond to shrinking 
strata to infinity. The full definition with relations on morphisms 
and the corresponding combinatorial analysis, will appear in~\cite{Ko01c}.
\end{rem}

\subsection{Resonance functors} 


Given a~functor $\cf:\catr\lra\ctop^*$, we introduce the following
notation:
$$\cf(Q(S,\Pi))=\cf(S)\bigg/\bigcup_{\un(\pi)\in\Pi}\im\cf(\pi S
\stackrel{\pi}{\hra}S).$$

\begin{df} \label{dfresf}
  A functor $\cf:\catr\lra\ctop^*$ is called a~{\bf resonance functor}
if it satisfies the following axioms:
\begin{enumerate}
\item [(A1)] {\bf Inclusion axiom.} 
  
  \nin If $S\in\catr_n$, and $\pi\in OP(n)$, then $\cf(\pi S
\stackrel{\pi}{\hra}S)$ is an~inclusion map, and 
$\cf(S)\big/\im\cf(\pi S\stackrel{\pi}{\hra}S)\simeq\cf(Q(S,\pi\dar S))$.

\item [(A2)] {\bf Relative resonance axiom.}

  \nin If, for some $S,T\in\catr_n$, and $\Pi,\Lambda\subseteq P(n)$,
  $[Q(S,\Pi)]=[Q(T,\Lambda)]$, then
  $\cf(Q(S,\Pi))\simeq\cf(Q(T,\Lambda))$.

\item [(A3)] {\bf Direct product axiom.}

\nin For two relative $n$-cuts $(S,\Pi)$ and $(T,\Lambda)$ we have
  $$\cf(S,\Pi)\times\cf(T,\Lambda)\simeq\cf((S,\Pi)\times(T,\Lambda)).$$

\end{enumerate}
\end{df}

\noindent
Given $S\in\catr_n$, and $\pi\in OP(n)$, let $i_{S,\pi}$ denote the
inclusion map $\cf(\pi S\stackrel{\pi}{\hra}S)$. There is a~canonical
homology long exact sequence associated to the triple
\begin{equation} \label{triple}
 \cf(\pi S)\stackrel{i_{S,\pi}}{\hra}\cf(S)\stackrel{p}{\lra}
\cf(Q(S,\pi\dar S)),
\end{equation}
namely
\begin{multline} \label{stls}
\dots\stackrel{\partial_*}{\lra}\widetilde H_n(\cf(\pi S))
\stackrel{(i_{S,\pi})_*}{\lra}\widetilde H_n(\cf(S))\stackrel{p_*}{\lra}
\widetilde H_n(\cf(Q(S,\pi\dar S)))\stackrel{\partial_*}{\lra}  \\
\widetilde H_{n-1}(\cf(\pi S))\stackrel{(i_{S,\pi})_*}{\lra}\dots
\end{multline}
We call \eqref{triple}, resp.~\eqref{stls}, the {\it standard triple},
resp.~the {\it standard long exact sequence} associated to the
morphism $\pi S\stackrel{\pi}{\hra}S$ and the functor $\cf$ (usually
$\cf$ is fixed, so its mentioning is omitted).


     \section{First applications} \label{s1.4}

\subsection{Resonance compatible stratifications} 
\label{ss4.1}

As mentioned in the Section~\ref{s1.1} we shall now look at the
natural strata of the spaces $X^{(n)}$. The strata are defined by
point coincidences and are indexed by number partitions of $n$. Let
$\Sigma_\lambda^X$ denote the stratum indexed by $\lambda$.

Let $\lambda$ be a number partition of $n$ and let $\tilde\lambda\in
OP(n)$ be $\lambda$ with some fixed order on the parts. Then
$\tilde\lambda$ can be thought of as a~vector with positive integer
coordinates in~$\br^n$. Let $S_{\tilde\lambda}$ be the set
$\{x\in\{-1,0,1\}^n\,|\,\langle x,\tilde\lambda\rangle =0\}$.
Obviously, $S_{\tilde\lambda}$ is an~$n$-cut, and the $n$-resonance
$S_\lambda$, which it defines, does not depend on the choice of
$\ti\lambda$, but only on the number partition~$\lambda$.

The crucial topological observation is that if $\nu$ is another
partition of $n$, such that $S_\lambda=S_\nu$, then the spaces $\slx$
and $\Sigma_\nu^X$ are homeomorphic. This is precisely the fact which
lead us to introduce resonances and the surrounding combinatorial
framework and to forget about the number partitions themselves.

That observation allows us to introduce a~functor $\cf$ mapping
$S_{\tilde\lambda}$ to $\Sigma_\lambda^X$; the morphisms map
accordingly. Clearly, $\cf(1^l)=X^{(l)}$. One can detect in this
example the justification for the names which we chose for the
morphisms of $\catr$: ``inclusions'' and ``gluings''. Furthermore, it
is easy, in this case, to verify the axioms of the
Definition~\ref{dfresf}, and hence to conclude that $\cf$ is
a~resonance functor. The only nontrivial point is the verification of
the second part of (A1), which we do in the next proposition.

\begin{prop} \label{prop4.1}
  Let $S$ be an $n$-cut and $\pi\in OP(n)$. Then $\cf(\nu
  S)\subseteq\cf(\pi S)$ if and only if $\un(\nu)\in\un(\pi)\dar S$.
\end{prop}

\pr It is obvious that all the steps of the definition of
$\un(\pi)\dar S$ which change the partition preserve the property
$\cf(\nu S)\subseteq\cf(\pi S)$, hence the {\it if} direction follows.

Assume now $\cf(\nu S)\subseteq\cf(\pi S)$. This means that there
exists $\tau\in OP(m)$, where $m$ is the number of parts of $\pi$,
such that $\cf(\tau\pi S)=\cf(\nu S)$. By definition,
$\tau\circ\pi\in\un(\pi)\dar S$. Now, we can reach $\un(\nu)$ from
$\un(\tau\circ\pi)$ by moves of type (3) from the definition of the
relative resonances. 

Indeed, if $\cf(\tau\pi S)=\cf(\nu S)=\Sigma_\lambda^X$, then the
sizes of the resulting blocks after gluing along $\tau\circ\pi$ and
along $\nu$ are the same. For every block $b$ of $\lambda$ we can go,
by means of moves of type (3), from the block of $\un(\tau\circ\pi)$
which glues to $b$ to the block of $\un(\nu)$ which glues to $b$.
Since we can do it for any block of $\lambda$, we can go from
$\un(\tau\circ\pi)$ to $\un(\nu)$, and hence 
$\un(\nu)\in\un(\pi)\dar S$. 
\qed

\vskip4pt

\nin In the context of this stratification the following central
question arises.

\vskip8pt

\nin {\bf The Main Problem.} {\it (Arnol'd, Shapiro,~\cite{Sh00}).
  Describe an algorithm which, for a~given resonance $\lambda$, would
  compute the Betti numbers of $\Sigma_\lambda^{S^1}$, or
  $\Sigma_\lambda^{S^2}$.}

\vskip8pt

The case of the strata $\Sigma_\lambda^{S^1}$ is simpler, essentially
because of the following elementary, but important property of smash
products: if $X$ and $Y$ are pointed spaces and $X$ is contractible,
then $X\wedge Y$ is also contractible.

In the subsequent subsections we shall look at a~few interesting
special cases, and also will be able to say a~few things about the
general problem.

\subsection{Resonances $(a^k,1^l)$}  \label{ss4.2}

Let $a,k,l$ be positive integers such that $a\geq 2$. Let $S$ be the
$(l+k)$-cut consisting of all the elements of
$\{-1,0,1\}^{l+k}$, which are orthogonal to the vector
$(\underbrace{1,\dots,1}_l,\underbrace{a,\dots,a}_k)$. Clearly,
the~$(l+k)$-resonance $[S]$ is equal to $(a^k,1^l)$. The case $l<a$ is
not very interesting, since then $(a^k,1^l)=(1^k)\times(1^l)$.
Therefore we may assume that $l\geq a$.

We would like to understand the topological properties of the space
$\cf(a^k,1^l)$. In general, this is rather hard. However, as the
following theorem shows, it is possible under some additional
conditions on $\cf$.

\begin{thm} \label{thm4.2}
  Let $\cf:\catr\lra\ctop^*$ be a~resonance functor, such that
  $\cf(1^l)$ is contractible for $l\geq 2$. Let $l=am+\epsilon$, where
  $0\leq\epsilon\leq a-1$.
  \begin{enumerate}
  \item [(a)] If $k\neq 1$, or $\epsilon\geq 2$, then $\cf(a^k,1^l)$
    is contractible.
  \item [(b)] If $k=1$, and $\epsilon\in\{0,1\}$, then
  \begin{equation} \label{eq:4.1}
   \cf(a^k,1^l)\simeq\susp^m(\cf(1)^{m+\epsilon+1}), 
  \end{equation}
  where $\cf(1)^{m+\epsilon+1}$ denotes the $(m+\epsilon+1)$-fold
  smash product.
  \end{enumerate}
\end{thm}

Since for the resonance functor $\cf$ described in the
subsection~\ref{ss4.1} we have $\cf(1^l)=X^{(l)}$, we have the
following corollary.

\begin{crl} \label{crl4.3}
  If $X^{(l)}$ is contractible for $l\geq 2$, then
  \begin{enumerate}
  \item [(a)] If $k\neq 1$, or $\epsilon\geq 2$ (again
    $l=am+\epsilon$), then $\Sigma_{(a^k,1^l)}^X$ is contractible.
  \item [(b)] If $k=1$, and $\epsilon\in\{0,1\}$, then
    $\Sigma_{(a^k,1^l)}^X\simeq\susp^m(X^{m+\epsilon+1})$, where
    $X^{m+\epsilon+1}$ denotes the $(m+\epsilon+1)$-fold smash
    product.
  \end{enumerate}
\end{crl}

\begin{rem} 
  Clearly, $(S^1)^{(l)}$ is contractible for $l\geq 2$, so the
  Corollary~\ref{crl4.3} is valid. In this situation, the case $k>1$
  was proved in~\cite{Ko00b}, and the case $k=1$
  in~\cite{BWa97,ShW98}.
\end{rem}

Before we proceed with proving the Theorem~\ref{thm4.2} we need
a~crucial lemma. Let $\pi\in P(k+l)$ be 
$(\{1,\dots,a\},\{a+1\},\{a+2\},\dots,\{k+l\})$. It is immediate
that $[\tilde\pi S]=(a^{k+1},1^{l-a})$, if $\un(\tilde\pi)=\pi$.

\begin{lm} \label{lm4.4}
  Let $S$ be as above, $T\in\catr_l$, such that $[T]=(1^l)$, and let
  $\nu$ be the partition
  $(\{1,\dots,a\},\{a+1\},\{a+2\},\dots,\{l\})$, then we have
  \begin{equation}
    \label{eq:4.2}
    [Q(S,\pi\dar S)]=[Q(T,\nu\dar T)]\times(1^k).
  \end{equation}
\end{lm}

\begin{rem} Lemma~\ref{lm4.4} is a~special case of the
  Lemma~\ref{lm4.6}, however we choose to include a~separate proof for
  it for two reasons: firstly, it is the first, still not too
  technical example of investigating the combinatorial structure of
  the resonance category, which is a~new object; secondly, the
  particular case of $(a^k,1^l)$ resonances was a~subject of
  substantial previous attention.
\end{rem}

\nin {\bf Proof of the Lemma~\ref{lm4.4}.} 

\vskip4pt

\nin Recall that by the definition of the direct product,
$$[Q(T,\nu\dar T)]\times(1^k)=[Q(T\times U,(\nu\dar T)\times
P(\{l+1,\dots,l+k\}))],$$ 
where $U\in\catr_k$ and $[U]=(1^k)$. Clearly, $(\nu\dar T)\times
P(\{l+1,\dots,l+k\})=\pi\dar S$, hence we just need to show that
$S\setminus(\pi\dar S)=(T\times U)\setminus((\nu\dar T)\times 
P(\{l+1,\dots,l+k\}))$. Note that $(T\times U)\setminus((\nu\dar T)\times 
P(\{l+1,\dots,l+k\}))=(T\setminus(\nu\dar T))\times U$.
Furthermore, 
\vskip-20pt 
$$S=\bigg\{(x_1,\dots,x_{l+k})\in\{-1,0,1\}^{l+k}\,\Bigm|\,
\sum_{j=l+1}^{l+k}x_j+a\sum_{i=1}^{l}x_i=0\bigg\},$$
\vskip-5pt
\noindent and
\vskip-20pt
\begin{multline*}
\bigcup_{\tau\in\pi\dar S}S^\tau=\bigg\{(x_1,\dots,x_{l+k})\in\{-1,0,1\}^{l+k}
\,\Bigm|\,\sum_{j=l+1}^{l+k}x_j+a\sum_{i=1}^{l}x_i=0, \\
\max(|\plus(x_1,\dots,x_l)|,|\minus(x_1,\dots,x_l)|)\geq a\bigg\}.
\end{multline*}
\vskip-5pt
\noindent
Therefore, by the definition of the relative resonances, we have
\vskip-20pt 
\begin{multline*}
S\setminus(\pi\dar S)=\bigg\{(x_1,\dots,x_{l+k})\in\{-1,0,1\}^{l+k}\,\Bigm|\,\\
\sum_{i=1}^l x_i=0,\,\,\sum_{j=l+1}^{l+k}x_j=0,\,\,
|\plus(x_1,\dots,x_l)|<a\bigg\}.
\end{multline*}
On the other hand, 
$(1^k)=[\{(y_1,\dots,y_k)\in\{-1,0,1\}^k\,|\,\sum_{j=1}^k y_j=0\}]$,
and
$$T\setminus(\nu\dar T)=\bigg\{(z_1,\dots,z_l)\in\{-1,0,1\}^l\,\Bigm|\,
\sum_{i=1}^l z_i=0,\,\,|\plus(z_1,\dots,z_l)|<a\bigg\},$$
which proves~\eqref{eq:4.2}. \qed

\vskip4pt
\nin {\bf Proof of the Theorem~\ref{thm4.2}.} 

\vskip4pt

\nin {\bf (a)} We use induction on~$l$. The case $l<a$ can be taken as
an~induction base, since then $(a^k,1^l)=(1^k)\times(1^l)$, hence, by
the axiom (A3), $\cf(a^k,1^l)=\cf(1^k)\wedge\cf(1^l)$, which is
contractible, since $\cf(1^k)$ is. Thus we assume that $l\geq a$, and
$\cf(a^k,1^{l'})$ is contractible for all $l'<l$.

Let $S$ and $\pi$ be as in the Lemma~\ref{lm4.4}. The standard triple
associated to the morphism $\pi S\stackrel{\pi}{\hra}S$ is
$\cf(a^{k+1},1^{l-a})\hra\cf(a^k,1^l)\ra\cf(a^k,1^l)/\cf(a^{k+1},1^{l-a})$.
Since, by the induction assumption, $\cf(a^{k+1},1^{l-a})$ is
contractible, we conclude that $\cf(a^k,1^l)\simeq\cf(a^k,1^l)/
\cf(a^{k+1},1^{l-a})$.

Basically by the definition, we have
$\cf(a^k,1^l)/\cf(a^{k+1},1^{l-a})= \cf(Q(S,\pi\dar S))$. On the other
hand, we have proved in the Lemma~\ref{lm4.4} that $[Q(S,\pi\dar
S)]=Q(T,\nu\dar T)\times (1^k)$, where $T$ and $\nu$ are described in
the formulation of that lemma. By axioms (A2) and (A3) we get that
$\cf(Q(S,\pi\dar S))\simeq\cf(Q(T,\nu\dar T))\wedge \cf(1^k)$, which
is contractible, since $\cf(1^k)$ is. Therefore, $\cf(a^k,1^l)$ is
also contractible.

\vskip4pt

\nin {\bf (b)} The argument is very similar to (a). We again assume
$l\geq a$, which implies $l\geq 2$. By the using the same ordered set
partition~$\pi$ as in (a), we get that $\cf(a,1^l)\simeq
\cf(a,1^l)/\cf(a^2,1^{l-a})$. Further, by Lemma~\ref{lm4.4} and the
axioms (A2) and (A3) we conclude that $\cf(a,1^l)\simeq
\cf(1)\wedge(\cf(1^l)/\cf(a,1^{l-a}))$. Since $\cf(1^l)$ is
contractible, we get 
\begin{equation}
  \label{eq:4.3}
  \cf(a,1^l)\simeq \cf(1)\wedge \susp\cf(a,1^{l-a}).
\end{equation}
Since $\cf(a)=\cf(1)$, $\cf(a,1)=\cf(1)\wedge\cf(1)$, and $\cf(a,1^l)$
is contractible if $2\leq l<a$, we obtain \eqref{eq:4.1} by the
repeated usage of~\eqref{eq:4.3}. \qed


\subsection{Resonances $(a^k,b^l)$} 

The algebraic invariants of these strata have not been computed
before, not even in the case $X=S^1$, and $\cf$ - the standard
resonance functor associated to the stratification of $X^{(n)}$.

We would like to apply a~technique similar to the one used in the
subsection~\ref{ss4.2}. A~problem is that, once one starts to ``glue''
$a$'s, one cannot get $b$'s in the same way as one could in the
previous section from 1's. Thus, we are forced to consider a~more 
general case of resonances, namely $(g^m,a^k,b^l)$, where $g$ is
the least common multiple of $a$ and $b$. Assume $g=a\cdot\bar a=
b\cdot\bar b$, and $b>a\geq 2$. Analogously with the Theorem~\ref{thm4.2}
we have the following result.


\begin{thm} \label{thm4.5}
  Let $\cf$ be as in the Theorem~\ref{thm4.2}. Let furthermore 
$k=x\cdot\bar a+\epsilon_1$, $l=y\cdot\bar b+\epsilon_2$, where 
$0\leq\epsilon_1<\bar a$, $0\leq\epsilon_2<\bar b$. Then

\begin{equation}
  \label{eq:4.4}
  \cf(g^m,a^k,b^l)\simeq
     \begin{cases}
        \susp^{x+y+m-1}(\cf(1)^{x+y+m+\epsilon_1+\epsilon_2}), & 
        \text{ if } m,\epsilon_1,\epsilon_2\in\{0,1\}; \\
        \text{ point, } & \text{ otherwise. } 
     \end{cases} 
\end{equation}
\end{thm}

Just as in the subsection~\ref{ss4.2} (Corollary~\ref{crl4.3}), the
Theorem~\ref{thm4.5} is true if one replaces $\cf(\lambda)$ with
$\Sigma_\lambda^{S^1}$.

The proof of the Theorem~\ref{thm4.5} follows the same general scheme
as that of the Theorem~\ref{thm4.2}, but the technical details are
more numerous. Again there is a~crucial combinatorial lemma.

Let $S$ be an~$(m+k+l)$-cut consisting of all the elements of
$\{-1,0,1\}^{m+k+l}$ which are orthogonal to the vector
$(\underbrace{a,\dots,a}_k,\underbrace{b,\dots,b}_l,
\underbrace{g,\dots,g}_m)$. Assume $k\geq\bar a$, and let an~unordered
set partition~$\pi$ be equal to $(\{1,\dots,\bar a\},\{\bar
a+1\},\dots,\{k+l+m\})$. We see that $[S]=(g^m,a^k,b^l)$, and
$[\tilde\pi S]=(g^{m+1},a^{k-\bar a},b^l)$, if $\pi=\un(\tilde\pi)$.

\begin{lm} \label{lm4.6}
Let $T\in\catr_k$, such that $[T]=(1^k)$, and 
$\nu=(\{1,\dots,\bar a\},\{\bar a+1\},\dots,\{k\})$, then
\begin{equation}
  \label{eq:4.5}
  [Q(S,\pi\dar S)]=[Q(T,\nu\dar T)]\times(\bar b^m,1^l).
\end{equation}
\end{lm}
\pr Again, it is easy to see that the sets of the partitions at infinity
on both sides of~\eqref{eq:4.5} coincide. Indeed,
$$[Q(T,\nu)]\times(\bar b^m,1^l)=[Q(T\times U,(\nu\dar T)\times
P(\{k+1,\dots,k+m+l\}))],$$
where $U\in\catr_{m+l}$, such that
$[U]=({\bar b}^m,1^l)$, and $(\nu\dar T)\times
P(\{k+1,\dots,k+m+l\})=\pi\dar S$. Also, we again have the equality
$$(T\times U)\setminus((\nu\dar T)\times
P(\{k+1,\dots,k+m+l\}))=T\setminus(\nu\dar T)\times U,$$
which greatly
helps to prove that the sets if the surviving elements on the two
sides of~\eqref{eq:4.5} coincide. 

By the definition
\begin{multline*}
S=\bigg\{(x_1,\dots,x_{k+l+m})\in\{-1,0,1\}^{k+l+m}\,\Bigm| \\
a\sum_{i=1}^{k}x_i+b\sum_{i=k+1}^{k+l}x_i+g\sum_{i=k+l+1}^{k+l+m}x_i=0\bigg\},
\end{multline*}
and
\begin{multline*}
\bigcup_{\tau\in\pi\dar S}S^\tau=\bigg\{(x_1,\dots,x_{k+l+m})
\in\{-1,0,1\}^{k+l+m}\,\Bigm|\,
a\sum_{i=1}^{k}x_i+b\sum_{i=k+1}^{k+l}x_i+ \\ g\sum_{i=k+l+1}^{k+l+m}x_i=0,\,\,
\max(|\plus(x_1,\dots,x_k)|,|\minus(x_1,\dots,x_k)|)\geq\bar a\bigg\}.$$
\end{multline*}
By the definition of the relative resonances and some elementary number 
theory we conclude that 
\begin{multline*}
S\setminus(\pi\dar S)=\bigg\{(x_1,\dots,x_{k+l+m})\in\{-1,0,1\}^{k+l+m}\,\Bigm|\,
 |\plus(x_1,\dots,x_k)|<\bar a,  \\
\sum_{i=1}^{k}x_i=0,\quad
b\sum_{i=k+1}^{k+l}x_i+g\sum_{i=k+l+1}^{k+l+m}x_i=0\bigg\}.$$
\end{multline*}
The number theory argument which we need is that if $ax+by+\lcm(a,b)z=0$,
then $\bar a\,\big|\,x$, where $\bar a\cdot a=\lcm(a,b)$. This can be 
seen by, for example, noticing that if $ax+by+\lcm(a,b)z=0$, then  
$b\,\big|\,ax$, but since also $a\,\big|\,ax$, we have 
$\lcm(a,b)\,\big|\,ax$, hence $\bar a\,\big|\,x$.

The equation~\eqref{eq:4.5} follows now from the earlier observations
together with the equalities \vskip-10pt
$$T\setminus(\nu\dar T)=\!\bigg\{(x_1,\dots,x_k)\in\{-1,0,1\}^k\Bigm|
|\plus(x_1,\dots,x_k)|<\bar a,\sum_{i=1}^{k}x_i=0\bigg\}$$
and \vskip-10pt
$$(\bar b^m,1^l)=\bigg[\bigg\{(y_1,\dots,y_{m+l})\in\{-1,0,1\}^{m+l}
\,\Bigm|\,\sum_{i=1}^l y_i+\bar b\sum_{i=l+1}^{l+m}x_i=0\bigg\}\bigg].
\,\,\mqed$$

\vskip4pt
\nin {\bf Proof of the Theorem~\ref{thm4.5}.}
  The cases $k<\bar a$ and $l<\bar b$ are easily reduced to the 
Theorem~\ref{thm4.2}. Assume therefore that $k\geq\bar a$ and 
$l\geq\bar b$. Recall also that $b>a\geq 2$, and hence $\bar a\geq 2$. 

Let $S$ and $\pi$ be as in the formulation of the Lemma~\ref{lm4.6}. 
The standard triple associated to the morphism 
$\pi S\stackrel{\pi}{\hra}S$ is
\begin{equation}
  \label{eq:4.6}
  \cf(g^{m+1},a^{k-\bar a},b^l)\hra\cf(g^m,a^k,b^l)\ra
  \cf(g^m,a^k,b^l)/\cf(g^{m+1},a^{k-\bar a},b^l).
\end{equation}
We break the rest of the proof into 3 cases.

\vskip4pt

\nin {\bf Case $m\geq 2$}. Again, we prove that $\cf(g^m,a^k,b^l)$ is
contractible by induction on~$k$. This is clear if $k<\bar a$. If
$k\geq\bar a$, it follows from~\eqref{eq:4.6} that $\cf(g^m,a^k,b^l)
\simeq\cf(g^m,a^k,b^l)/\cf(g^{m+1},a^{k-\bar a},b^l)=\cf(Q(S,\pi\dar
S))$. By Lemma~\ref{lm4.6} we conclude that $\cf(g^m,a^k,b^l)\simeq
\cf(Q(T,\nu\dar T))\wedge\cf(\bar b^m,1^l)$. By the
Theorem~\ref{thm4.2}, $\cf(\bar b^m,1^l)$ is contractible, hence so is
$\cf(g^m,a^k,b^l)$.

\vskip4pt

\nin {\bf Case $m=0$}. By Lemma~\ref{lm4.6} we get that
$\cf(Q(S,\pi\dar S))\simeq\cf(Q(T,\nu\dar T))\wedge\cf(1^l)$. Since
$l\geq 2$, we have that $\cf(1^l)$ is contractible, hence so is
$\cf(Q(S,\pi\dar S))= \cf(a^k,b^l)/\cf(g,a^{k-\bar a},b^l)$.
Therefore, by~\eqref{eq:4.6} $\cf(a^k,b^l)\simeq\cf(g,a^{k-\bar
  a},b^l)$.

\vskip4pt

\nin {\bf Case $m=1$}. Since $\cf(g^2,a^{k-\bar a},b^l)$ is
contractible, we conclude by~\eqref{eq:4.6} that
$\cf(g,a^k,b^l)\simeq\cf(g,a^k,b^l)/ \cf(g^2,a^{k-\bar
  a},b^l)=\cf(Q(S,\pi\dar S))$. By Lemma~\ref{lm4.6}, and the
properties of the resonance functors, we have
\begin{equation}
  \label{eq:4.7}
\cf(g,a^k,b^l)\simeq\cf(\bar b,1^l)\wedge(\cf(1^k)/\cf(\bar a,1^{k-\bar a}))
\simeq\cf(\bar b,1^l)\wedge\susp(\cf(\bar a,1^{k-\bar a})).
\end{equation} 
By the repeated usage of~\eqref{eq:4.7} we obtain~\eqref{eq:4.4}.
\qed

     \section{Sequential resonances} \label{s1.5}

\subsection{The structure theory of strata associated to sequential 
resonances} 

\begin{df}\label{dfseqres}
  Let $\lambda=(\lambda_1,\dots,\lambda_n)$,
  $\lambda_1\leq\dots\leq\lambda_n$, be a~number partition.  We call
  $\lambda$ {\bf sequential} if, whenever $\sum_{i\in
    I}\lambda_i=\sum_{j\in J}\lambda_j$, and $q\in I$, such that
  $q=\max(I\cup J)$, then there exists $\widetilde J\subseteq J$, such
  that $\lambda_q=\sum_{j\in\widetilde J}\lambda_j$.
  
  Correspondingly, we call a~resonance $S$ sequential, if it can be
  associated to a~sequential partition.
\end{df}

\nin Note that the set of sequential partitions is closed under
removing blocks.

\begin{exams} {\rm Sequential partitions.}
\begin{enumerate}
\item[(1)] All partitions whose blocks are equal to powers of some number;
\item[(2)] $(a^k,b^l,1^m)$, such that $a>bl$; more generally
  $(a_1^{k_1},\dots,a_t^{k_t},1^m)$, such that $a_i>\sum_{j=i+1}^t a_j
  k_j$, for all $i\in[t]$.
\end{enumerate}
\end{exams}

Through the rest of this subsection, we let $\lambda$ be as in the
Definition~\ref{dfseqres}. For such $\lambda$ we use the following
additional notations:
\begin{itemize}
\item $mm(\lambda)=|\{i\in[n]\,|\,\lambda_i=\lambda_n\}|$. 
In other words
$\lambda_{n-mm(\lambda)}\neq\lambda_{n-mm(\lambda)+1}=\dots=\lambda_n$.
\item $I(\lambda)\subseteq[n]$ is the lexicographically maximal set
  (see below the convention that we use to order lexicographically),
  such that $|I(\lambda)|\geq 2$, and $\lambda_n=\sum_{i\in
    I(\lambda)}\lambda_i$. Note that it may happen that $I(\lambda)$
  does not exist, in which case $\cf(\lambda)\simeq
  \cf(\lambda_1,\dots,\lambda_{n-mm(\lambda)})\wedge\cf(1^{mm(\lambda)})$,
  and can be dealt with by induction.
\end{itemize}

Let $n$ be a~positive integer. We use the following convention for the
lexicographic order on $[n]$. For $A=\{a_1,\dots,a_k\}$,
$B=\{b_1,\dots,b_m\}$, $A,B\subseteq[n]$, $a_1\leq\dots\leq a_k$,
$b_1\leq\dots\leq b_m$, we say that $A$ is lexicographically larger
than~$B$ if, either $A\supseteq B$, or there exists $q<\min(k,m)$,
such that $a_k=b_m$, $a_{k-1}=b_{m-1}$, $\dots$,
$a_{k-q+1}=b_{m-q+1}$, and $a_{k-q}>b_{m-q}$.

\begin{prop} \label{prop4.9}
  If $\lambda=(\lambda_1,\dots,\lambda_n)$,
  $\lambda_1\leq\dots\leq\lambda_n$, is a sequential partition, then
  so is $\bar\lambda=(\lambda_{j_1},\dots,\lambda_{j_t},\sum_{i\in
    I(\lambda)}\lambda_i)$, where $t=n-|I(\lambda)|$, and
  $\{j_1,\dots,j_t\}=[n]\setminus I(\lambda)$.
\end{prop}
\pr Let $\bar\lambda_1=\lambda_{j_1},\dots,\bar\lambda_t=\lambda_{j_t}$,
$\bar\lambda_{t+1}=\sum_{i\in I(\lambda)}\lambda_i$.
We need to check the condition of the Definition~\ref{dfseqres}
for the identity
\begin{equation}
  \label{eq:4.8a}
  \sum_{i\in I}\bar\lambda_i=\sum_{j\in J}\bar\lambda_j.
\end{equation}

If $t+1\not\in I\cup J$, then it follows from the assumption that
$\lambda$ is sequential. Assume $t+1\in I$. If
$\bar\lambda_j=\lambda_n$, for some $j\in J$, take $\widetilde
J=\{j\}$, and we are done. If $\bar\lambda_i=\lambda_n$, for some
$i\in I\setminus\{t+1\}$, then, since $\lambda$ is sequential, there
exists $\widetilde J\subseteq J$, such that $\sum_{j\in\widetilde
  J}\bar\lambda_j=\lambda_n=\bar\lambda_{t+1}$, and we are done again.

Finally, assume $\bar\lambda_i\neq\lambda_n$, for $i\in(I\cup
J)\setminus\{t+1\}$. Substituting $\lambda_n$ instead of
$\bar\lambda_{t+1}$ into the identity~\eqref{eq:4.8a} is allowed,
since $\lambda_n$ does not appear among $\{\bar\lambda_i\}_{i\in(I\cup
  J)\setminus\{t+1\}}$. This gives us an identity for $\lambda$, and
again, since $\lambda$ is sequential, we find the desired set
$\widetilde J\subseteq J$, such that $\sum_{j\in\widetilde
  J}\bar\lambda_j=\bar\lambda_{t+1}$. \qed

\vskip4pt

Let $S\in\catr_n$ be the set of all elements of $\{-1,0,1\}^n$, which
are orthogonal to the vector $\lambda=(\lambda_1,\dots,\lambda_n)$.
Clearly, $[S]=\lambda$. Let $\pi\in P(n)$ be the partition whose only
nonsingleton block is given by $I(\lambda)$. The next lemma expresses
the main combinatorial property of sequential partitions.

\begin{lm}\label{green}
  Let $\tau\in P(n)$ be a~partition which has only one nonsingleton
block~$B$, and assume $\lambda_n=\sum_{i\in B}\lambda_i$. 
Then $\tau\in\pi\dar S$.
\end{lm}
\pr Assume there exists partitions $\tau$ as in the formulation of the
lemma, such that $\tau\not\in\pi\dar S$. Choose one so that the
block~$B$ is lexicographically largest possible. Let $C=B\cap
I(\lambda)$. By the definition of $I(\lambda)$, and the choice of $B$,
we have $\sum_{i\in I(\lambda)\setminus C}\lambda_i=
\sum_{j\in B\setminus C}\lambda_j$, and $q\in I(\lambda)\setminus C$,
where $q=\max((I(\lambda)\cup B)\setminus C)$.

Since partition $\lambda$ is sequential, there exists $D\subseteq
B\setminus C$, such that $\lambda_q=\sum_{j\in D}\lambda_j$.  Let
$\gamma\in P(n)$ be the partition whose only nonsingleton block is
$G=(B\setminus D)\cup\{q\}$. Clearly, $\sum_{i\in G}\lambda_i=
\lambda_n$, and $|G|\geq 2$. By the choice of $q$, $G$ is 
lexicographically larger than $B$, hence $\gamma\in\pi\dar S$.

Let furthermore $\tilde\gamma\in P(n)$ be the partition having two
nonsingleton blocks: $D$ and $G$. By the Definition~\ref{df3.1}(2) if
$\gamma\in\pi\dar S$, then $\tilde\gamma\in\pi\dar S$. By the
Definition~\ref{df3.1}(3), if $\tilde\gamma\in\pi\dar S$, then
$\tau\in\pi\dar S$, which yields a~contradiction.
\qed

\vskip4pt

Let $T\in\catr_{n-mm(\lambda)}$ be the set of all elements of
$\{-1,0,1\}^{n-mm(\lambda)}$, which are orthogonal to the vector
$(\lambda_1,\dots,\lambda_{n-mm(\lambda)})$. Let $\nu\in
P(n-mm(\lambda))$ be the partition whose only nonsingleton block is
given by $I(\lambda)$. We are now ready to state the combinatorial
result which is crucial for our topological applications.

\begin{lm} \label{lm4.11}
  \begin{equation}
    \label{eq:4.9}
    [Q(S,\pi\dar S)]=[Q(T,\nu\dar T)]\times(1^{mm(\lambda)}).
  \end{equation}
\end{lm}
\pr By definition we must verify that the sets of partitions at
infinity and the surviving elements coincide on both sides of the
equation~\eqref{eq:4.9}. 

Let us start with the partitions at infinity. Once filtered through the
Pro\-po\-si\-tion~\ref{prop4.1}, the identity $\pi\dar S=(\nu\dar T)\times
P(\{n-mm(\lambda)+1,\dots,n\})$ becomes essentially tautological. Both
sides consist of the partitions $\tau=(\tau_1,\dots,\tau_k)\in P(n)$,
such that the number partition 
$(\sum_{i\in\tau_1}\lambda_i,\dots,\sum_{i\in\tau_k}\lambda_i)$ 
can be obtained from the number partition 
$(\lambda_{j_1},\dots,\lambda_{j_t},\sum_{i\in I(\lambda)}\lambda_i)$,
where $\{j_1,\dots,j_t\}=[n]\setminus I(\lambda)$, by summing parts.

Let us now look at the surviving elements. Obviously,
$S\setminus(\pi\dar S)\supseteq(T\setminus(\nu\dar T))\times U$, where
$U\in\catr_k$, such that $[U]=(1^{mm(\lambda)})$, and we need to show
the converse inclusion. Let $x=(x_1,\dots,x_n)\in S$, such that
$\sum_{i=n-mm(\lambda)+1}^n x_i\neq 0$ (otherwise
$x\in(T\setminus(\nu\dar T))\times U$), we can assume
$\sum_{i=n-mm(\lambda)+1}^n x_i>0$. Then, since $S$ is a~sequential
resonance, there exists $y=(y_1,\dots,y_n)\in S$, such that
\begin{itemize}
\item if $y_i\neq 0$, then $x_i=y_i$; 
\item $|\plus(y)|=1$, and $\plus(y)\subseteq\{n-mm(\lambda)+1,\dots,n\}$.
\end{itemize}
This means that $y\in S^\tau$, for some $\tau\in P(n)$, which 
satisfies the conditions of the Lemma~\ref{green},
which implies that $\tau\in\pi\dar S$. On the other hand, $y\in
S^\tau$ necessitates $x\in S^\tau$, and hence $x\not\in
S\setminus(\pi\dar S)$. This finishes the proof of the lemma. 
\qed

\vskip4pt

Just as before, this combinatorial fact about the resonances
translates into a~topological statement, which can be further
strengthened by requiring some additional properties from~$\lambda$.

\begin{df}
  Let $\lambda=(\lambda_1,\dots,\lambda_n)$,
  $\lambda_1\leq\dots\leq\lambda_n$, be a~sequential partition, and
  let $q=\max I(\lambda)$. $\lambda$ is called {\bf strongly
    sequential}, if there exists $J\subseteq
  I(\lambda)\setminus\{q\}$, such that $\lambda_q=\sum_{i\in
    J}\lambda_i$ (note that we do not require $|J|\geq 2$).
\end{df}

We are now in a~position to prove the main topological structure
theorem concerning the sequential resonances.

\begin{thm} \label{thm4.13}
  Let $\cf$ be as in the Theorem~\ref{thm4.2}. Let $\lambda$ be 
a~sequential partition, such that $I(\lambda)$ exists, then
\begin{enumerate}
\item[(1)] if $mm(\lambda)\geq 2$, then $\cf(\lambda)$ is contractible;
\item[(2)] if $mm(\lambda)=1$, then 
$\cf(\lambda)\simeq\cf(Q(T,\nu\dar T))\wedge\cf(1)$, and we have
the inclusion triple $\cf(\mu)\stackrel{i}{\hra}\cf(\lambda_1,\dots,
\lambda_{n-1})\ra\cf(Q(T,\nu\dar T))$, where $\mu=(\lambda_{j_1},\dots,
\lambda_{j_t})$, $\{j_1,\dots,j_t\}=[n]\setminus I(\lambda)$,
and $\nu\in P(n-mm(\lambda))$ is the partition whose only nonsingleton
block is given by $I(\lambda)$.

If moreover $\lambda$ is strongly sequential, then the map $i$
is homotopic to a~trivial map (mapping everything to a~point), hence
the triple splits and we conclude that
\begin{equation}
  \label{eq:4.10}
  \cf(\lambda)\simeq(\cf(1)\wedge\cf(\lambda_1,\dots,\lambda_{n-1}))
  \vee\susp(\cf(1)\wedge\cf(\mu)).
\end{equation}
\end{enumerate}
\end{thm}
\pr 

\vskip4pt

\nin 
{\bf (1)} We use induction on $\sum_{i=1}^{n-mm(\lambda)}\lambda_i$.
If $I(\lambda)$ does not exist, then $\lambda_n$ is independent,
i.e., $\cf(\lambda)\simeq\cf(\lambda_1,\dots,\lambda_{n-mm(\lambda)})
\times\cf(1^{mm(\lambda)})$, and hence $\cf(\lambda)$ is contractible.
Otherwise consider the inclusion triple
\begin{equation}
  \label{eq:4.11}
  \cf(\bar\lambda)\hra\cf(\lambda)\ra\cf(\lambda)/\cf(\bar\lambda)=
\cf(Q(S,\pi\dar S)),
\end{equation}
where $\bar\lambda=(\lambda_{j_1}, \dots,\lambda_{j_t},\sum_{i\in
I(\lambda)}\lambda_i)$, and $\pi\in P(n)$ is the partition whose 
only nonsingleton block is given by $I(\lambda)$. By the induction 
assumption $\cf(\bar\lambda)$ is contractible. On the other hand, 
by Lemma~\ref{lm4.11}, $\cf(Q(S,\pi\dar S))\simeq\cf(Q(T,\nu\dar
T))\wedge\cf(1^{mm(\lambda)})$, which is also contractible if
$mm(\lambda)\geq 2$.

\vskip4pt

\nin 
{\bf (2)} If $mm(\lambda)=1$, then we can conclude from~\eqref{eq:4.11}
that $\cf(\lambda)\simeq\cf(1)\wedge\cf(Q(T,\nu\dar T))$. Next,
consider the inclusion triple 
\begin{equation}
  \label{eq:4.12}
  \cf(\mu)\stackrel{i}{\hra}\cf(\lambda_1,\dots,\lambda_{n-1})\ra
  \cf(Q(T,\nu\dar T)).
\end{equation}
  If $\lambda$ is strongly sequential, then there exists
$J\subseteq I(\lambda)\setminus\{q\}$, such that 
$\lambda_q=\sum_{i\in J}\lambda_i$ (here $q=\max I(\lambda)$). 
The map $i$ factors:
\begin{equation}
  \label{eq:4.13}
  \cf(\mu)\stackrel{i_1}{\hra}\cf(\lambda_{p_1},\dots,
\lambda_{p_{n-1-|J|}},\sum_{i\in I(\lambda)}\lambda_i)
\stackrel{i_2}{\hra}\cf(\lambda_1,\dots,\lambda_{n-1}),
\end{equation}
where $\{p_1,\dots,p_{n-1-|J|}\}=[n-1]\setminus J$. Since
$(\lambda_{p_1},\dots,\lambda_{p_{n-1-|J|}},\sum_{i\in
  I(\lambda)}\lambda_i)$ is sequential, and $mm((\lambda_{p_1},
\dots,\lambda_{p_{n-1-|J|}},\sum_{i\in I(\lambda)}\lambda_i))\geq 2$, 
we can conclude that the middle space in~\eqref{eq:4.13} is 
contractible, and hence $i$ in~\eqref{eq:4.12} is homotopic to 
a~trivial map. This yields the conclusion. 
\qed

\subsection{Resonances $(a^k,b^l,1^m)$} 

We give here the first application of the structure theory 
described in the previous subsection.

\begin{thm} \label{thm4.14}
  Let $a,b,k,l,m,r$ be positive integers, such that $b>1$, $m\geq r$,
  and $a=bl+r$.  Then
  \begin{equation}
    \label{eq:4.14}
    \cf(a^k,b^l,1^m)\simeq\susp(\cf(1^k)\wedge\cf(a,1^{m-r}))\vee
    (\cf(1^k)\wedge\cf(b^l,1^m)).
  \end{equation}
\end{thm}

\begin{rem}
  The restriction $m\geq r$ is unimportant. Indeed, if $m<r$, then
  $a>bl+m$, hence $a$ is not involved in any resonance other than
  $a=a$. This implies that
  $\cf(a^k,b^l,1^m)=\cf(1^k)\times\cf(b^l,1^m)$, and we have
  determined the homotopy type of $\cf(a^k,b^l,1^m)$ by the previous
  computations.
\end{rem}

\nin {\bf Proof of the Theorem~\ref{thm4.14}.}

\vskip4pt

\nin Obviously, the condition $a>bl$ guarantees that the partition
$(a^k,b^l,1^m)$ is sequential, hence the Theorem~\ref{thm4.13} is
valid. It follows that if $k\geq 2$, then $\cf(a^k,b^l,1^m)$ is
contractible, hence~\eqref{eq:4.14} is true.

Furthermore, if $l\geq 2$, or, $l=1$ and $m\geq b$, then $(a,b^l,1^m)$
is strongly sequential, hence in this case~\eqref{eq:4.10} is valid,
which in new notations becomes precisely the equation~\eqref{eq:4.14}.

Finally, assume $l=1$ and $b>m\geq r\geq 1$. Let $a=b+d$. If
$\cf(a,1^{m-d})$ or $\cf(b,1^m)$ is contractible, then the map~$i$ in
the inclusion triple $\cf(a,1^{m-d})\stackrel{i}{\hra}\cf(b,1^m)\ra
\cf(b,1^m)/\cf(a,1^{m-d})$ is homotopic to a~trivial map, and we again
conclude~\eqref{eq:4.14}. If both of these spaces are not contractible
then $\cf(a,1^{m-d})\simeq S^{2y+\epsilon_2+1}$ and $\cf(b,1^m)\simeq
S^{2x+\epsilon_1+1}$, where nonnegative integers
$x,y,\epsilon_1,\epsilon_2$ are defined by
\begin{equation}
  \label{eq:4.15}
  m=bx+\epsilon_1,\,\,\,m-d=(b+d)y+\epsilon_2,\,\,\,
\epsilon_1,\epsilon_2\in\{0,1\}.
\end{equation}
Let us show that $2x+\epsilon_1>2y+\epsilon_2$. If $x>y$, then
$2x+\epsilon_1\geq 2x\geq 2y+2>2y+\epsilon_2$. From~\eqref{eq:4.15} we
have that $b(x-y)=d+dy+\epsilon_2-\epsilon_1$. If $x\leq y$, then the
left hand side is nonpositive. On the other hand, since $d\geq 1$, the
right hand side is nonnegative. Hence, both sides are equal to 0,
which implies $x=y$, $d=\epsilon_1=1$, $\epsilon_2=y=0$. This yields
$2x+\epsilon_1>2y+\epsilon_2$.

The homotopic triviality of the map $i$ follows now from the fact that
the homotopy groups of a~sphere are trivial up to the dimension of
that sphere, i.e., $\pi_k(S^n)=0$, for $0\leq k\leq n-1$.  
\qed

\subsection{Resonances consisting of powers} 

Let us fix an integer $a\geq 2$. In this subsection we will study the
topology of the strata indexed by the following class of partitions:
all number partitions whose blocks are powers of $a$. Let $\Lambda(a)$
denote the set of all such partitions. 

It is convenient to introduce a~different notation for the partitions
in this class. Let $\Lambda_a(\alpha_n,\alpha_{n-1},\dots,\alpha_0)$,
where $\alpha_n,\alpha_{n-1},\dots,\alpha_0$ are nonnegative integers,
denote $\lambda\in\Lambda(a)$, which consists of $\alpha_n$ parts
equal to $a^n$, $\alpha_{n-1}$ parts equal to $a^{n-1}$, $\dots$,
$\alpha_0$ parts equal to 1. For example
$(8,4,2,2,2,1,1,1,1,1,1)=\Lambda_2(1,1,3,6)$.

Obviously, all partitions from $\Lambda(a)$ are strongly sequential,
hence the Theorem~\ref{thm4.13} applies, and it yields:
\begin{enumerate}
\item[(1)] if $\alpha_n\geq 2$, then $\cf(\Lambda_a(\alpha_n,\alpha_{n-1},
  \dots,\alpha_0))$ is contractible;
\item[(2)] if $I(\Lambda_a(1,\alpha_{n-1},\dots,\alpha_0))$ exists, then
  \begin{multline}
    \label{eq:4.16}
  \cf(\Lambda_a(1,\alpha_{n-1},\dots,\alpha_0)) \simeq 
  (\cf(1)\wedge\cf(\Lambda_a(\alpha_{n-1},\dots,\alpha_0)))\vee  \\
  (S^1\wedge\cf(1)\wedge\cf(\Lambda_a(1,\beta_{n-1},\dots,\beta_0))),
  \end{multline}
where $\Lambda_a(1,\beta_{n-1},\dots,\beta_0)$ is obtained from
$\Lambda_a(1,\alpha_{n-1},\dots,\alpha_0)$ by removing the blocks
indexed by $I(\Lambda_a(1,\alpha_{n-1},\dots,\alpha_0))$.
\item[(3)] If $I(\Lambda_a(1,\alpha_{n-1},\dots,\alpha_0))$ does not exist,
then
\begin{equation}
  \label{eq:4.17}
  \cf(\Lambda_a(1,\alpha_{n-1},\dots,\alpha_0)) \simeq 
  \cf(1)\wedge\cf(\Lambda_a(\alpha_{n-1},\dots,\alpha_0)).
\end{equation}
\end{enumerate}

It is immediate from the formulae~\eqref{eq:4.16} and~\eqref{eq:4.17}
that each topological space $\cf(\Lambda_a(\alpha_n,\alpha_{n-1},
\dots,\alpha_0))$ is homotopy equivalent to a~wedge of spaces of the
form $\cf(1)^\alpha\wedge S^\beta$, where $\cf(1)^\alpha$ means
an~$\alpha$-fold smash product of $\cf(1)$. The natural combinatorial
question which arises is how to enumerate these spaces. We shall now
construct a~combinatorial model: a~weighted graph which yields such
an~enumeration.

\begin{df}
  Let $\lambda=\Lambda_a(\alpha_n,\alpha_{n-1},\dots,\alpha_0))$, and
  set $\alpha_{-1}=1$. $\Gamma_\lambda$ is a~directed weighted graph
  on the set of vertices $\{n,\dots,0,-1\}$ whose edges and weights
  are defined by the following rule. For $x,x+d\in\{-1,0,\dots,n\}$,
  $d\geq 1$, there exists an~edge $e(x,x+d)$ (the edge is directed
  {\it from} $x$ {\it to} $x+d$) if and only if $\alpha_x>0$,
  $\alpha_{x+d}>0$, and
$$a^d\,\big|\,a^{d-1}\alpha_{x+d-1}+a^{d-2}\alpha_{x+d-2}+\dots+
a\alpha_{x+1}+\alpha_x-1.$$
In this case the weight of the edge is
defined as
$$w(x,x+d)=(a^{d-1}\alpha_{x+d-1}+a^{d-2}\alpha_{x+d-2}+\dots+
a\alpha_{x+1}+\alpha_x-1)\cdot a^{-d}.$$
\end{df}

\nin Note that if $d\geq 2$ and there exists an edge $e(x,x+d)$, then
there exists an edge $e(x,x+d-1)$.

We call a~directed path in $\Gamma_\lambda$ {\it complete} if it
starts in $-1$ and ends in $n$. Let $\gamma$ be a~complete path in
$\Gamma_\lambda$ consisting of $t$ edges, $\gamma=(e(x_0,x_1),\dots,
e(x_{t-1},x_t))$, where $x_0=-1$, and $x_t=n$. The weight of $\gamma$
is defined to be the pair $(l(\gamma),w(\gamma))$, where
$l(\gamma)=t$, and $w(\gamma)=\sum_{i=1}^t w(x_{i-1},x_i)$.

\begin{thm}\label{thm4.15}
  Let $\lambda\in\Lambda(a)$,
  $\lambda=\Lambda_a(1,\alpha_{n-1},\dots,\alpha_0)$, then
  \begin{equation}
    \label{eq:4.18}
    \cf(\lambda)\simeq\bigvee_{\gamma}(\cf(1)^{l(\gamma)+w(\gamma)}
    \wedge S^{w(\gamma)}),
  \end{equation}
  where the wedge is taken over all complete paths of
  $\Gamma_\lambda$.
\end{thm}
\pr We use induction on $\sum_{i=0}^n \alpha_i$. As the base of the
induction we take the case
$\lambda=\Lambda_a(1,\underbrace{0,\dots,0}_n)$. In this case
$\Gamma_\lambda$ is a~graph with only one edge $e(-1,n)$, $w(-1,n)=0$.
Thus, there is only one complete path. It has weight $(1,0)$, and
$\cf(\lambda)\simeq\cf(1)$.

Next, we prove the induction step. We break up the proof in three
cases. Let $t\in[n-1]$ be the maximal index for which $\alpha_t\neq 0$.

\vskip4pt

\nin {\bf Case 1.} {\it $I(\lambda)$ does not exist.}

\nin By~(\ref{eq:4.17}) we have
\begin{equation}
  \label{eq:4.19}
  \cf(\lambda)\simeq\cf(1)\wedge\cf(\Lambda_a(\alpha_{n-1},\dots,\alpha_0)).
\end{equation}
On the other hand, $I(\lambda)$ does not exist if and only if
$a^n>\alpha_{n-1}a^{n-1}+\dots+\alpha_1 a+\alpha_0$. We also know that
$\lambda\neq\Lambda_a(1,0,\dots,0)$, i.e.,
$\alpha_{n-1}a^{n-1}+\dots+\alpha_1 a+\alpha_0>0$. This implies that
there is at most one edge of the type $e(x,n)$. This edge exists
if and only if $\alpha_x=1$, and $\alpha_{n-1}=\dots=\alpha_{x+1}=0$,
in which case $w(x,n)=0$. 

If this edge does not exist then there are no complete paths in
$\Gamma_\lambda$ and, at the same time
$\cf(\Lambda_a(\alpha_{n-1},\dots,\alpha_0))$ is contractible by the
previous observations. This agrees with~(\ref{eq:4.18}).

If, on the other hand, this edge does exist, then all complete paths
$\gamma$ must be of the type $\gamma=(\tilde\gamma,e(x,n))$, where
$\tilde\gamma$ is a~complete path from $-1$ to $x$. Also in this 
case~(\ref{eq:4.19}) agrees with~(\ref{eq:4.18}).

\vskip4pt

\nin {\bf Case 2.} {\it $I(\lambda)$ exists and $\alpha_t\geq 2$.}

\nin In this case $\cf(\Lambda_a(\alpha_{n-1},\dots,\alpha_0))$ is
contractible, and 
\begin{equation}
  \label{eq:4.20}
  \cf(\lambda)\simeq S^1\wedge\cf(1)\wedge\cf(\Lambda_a(1,\beta_{n-1},
\dots,\beta_0)),
\end{equation}
where $\beta_{n-1}, \dots,\beta_0$ are as in~(\ref{eq:4.16}).

Let $q\in\{n-1,\dots,0,-1\}$ be the maximal index for which
$\beta_q\neq 0$ (we assume $\beta_{-1}=1$). Let $\tilde\lambda=
(1,\beta_{n-1},\dots,\beta_0)$. We can describe the graph
$\Gamma_{\tilde\lambda}$: it is obtained from $\Gamma_\lambda$ by
\begin{enumerate}
\item[(1)] removing the edges which have one of the endpoints 
in the set  $\{n-1,\dots,$ $q+1\}$;
\item[(2)] decreasing the weight of every existing edge $e(x,n)$ by $1$;
\item[(3)] keeping all existing edges with the old weights on the set
  $\{q,q-1,\dots,-1\}$.
\end{enumerate}

This operation on $\Gamma_\lambda$ is well-defined, since there can be
no edges in $\Gamma_\Lambda$ of the type $e(x,n)$, for
$x\in\{n-1,\dots,q+1\}$, and since the weight of edges $e(x,n)$, for
$x\in\{q,\dots,0,-1\}$ must be at least 1, as $b_q\neq 0$.
Furthermore, it is clear from the above combinatorial description of
$\Gamma_{\tilde\lambda}$, that the set of the complete paths of
$\Gamma_{\tilde\lambda}$ is the same as that of $\Gamma_\lambda$, and
that the weights of the edges in these paths are also the same except
for the edge with the endpoint $n$, whose weight has been decreased
by~1. Thus,~(\ref{eq:4.20}) agrees with~(\ref{eq:4.18}) in this case.

\vskip4pt

\nin {\bf Case 3.} {\it $I(\lambda)$ exists and $\alpha_t=1$.}

\nin This case is rather similar to the case 2, except that there is
an~edge $e(t,n)$ of weight $0$. Thus, $\Gamma_{\tilde\lambda}$
bookkeeps all the complete paths of $\Gamma_\lambda$, except for the
ones which have this edge $e(t,n)$.

However, the first term of the right hand side of~(\ref{eq:4.16})
bookkeeps the paths $(\tilde\gamma,e(t,n))$, just like in the case~1.
Since the set of all complete paths of $\Gamma_\lambda$ is the
disjoint union of the sets of those paths which contain $e(t,n)$, and
those which do not, we again get that~(\ref{eq:4.16}) provides the
inductive step for~(\ref{eq:4.18}).
\qed

\begin{exams} $\,$

\vskip4pt

\nin (1) Let $\lambda=(a,1^l)$, for $a\geq 2$. Then $\Gamma_\lambda$
is a~graph on the vertex set $\{1,0,-1\}$ having either one or two
edges:
\begin{enumerate}
\item it has in any case the edge $e(-1,0)$, $w(-1,0)=0$; 
\item if $a$ divides $l$, then it has the edge $e(-1,1)$, in which case
  $w(-1,1)=l/a$;
\item if $a$ divides $l-1$, then it has the edge $e(0,1)$, in which case
  $w(0,1)=(l-1)/a$.
\end{enumerate}
Clearly the Theorem~\ref{thm4.15} agrees with the
Theorem~\ref{thm4.2}. Indeed, if $\epsilon\not\in\{0,1\}$ (where
$\epsilon$ is taken from the formulation of the Theorem~\ref{thm4.2}),
then there are no complete paths in $\Gamma_\lambda$. If $\epsilon=0$,
then there is one path $(-1,1)$ of weight $(1,l/a)$; and if
$\epsilon=1$, then there is one path $((-1,0),(0,1))$ of weight
$(2,(l-1)/a)$. Thus,~(\ref{eq:4.18}) and~(\ref{eq:4.1}) are equivalent
in this case.


$$\begin{array}{c}
\epsffile{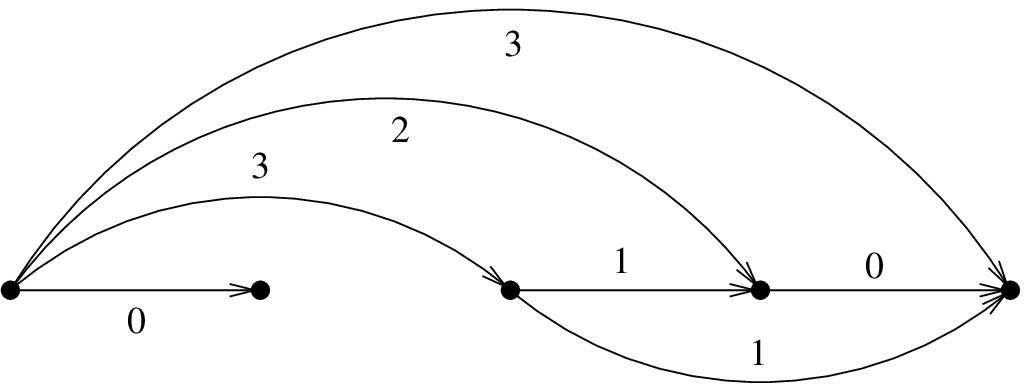}\\
\text{Figure 1.1.}\,\, \Gamma_{(8,4,2,2,2,1,1,1,1,1,1)}. 
\end{array}$$

\nin (2) Let $\lambda=(8,4,2,2,2,1,1,1,1,1,1)=\Lambda_2(1,1,3,6)$.
The graph $\Gamma_\lambda$ is shown on the Figure~1.1.
It has 4 directed paths from -1 to 3 and, by the Theorem~~\ref{thm4.15}, 
we have 
$$\cf(\lambda)\simeq
(\cf(1)^3\wedge S^2)\vee
(\cf(1)^5\wedge S^3)\vee
(\cf(1)^6\wedge S^4)\vee
(\cf(1)^7\wedge S^4),$$ 
in particular $\Sigma_\lambda^{\Bbb R}\simeq S^5\vee
S^8\vee S^{10}\vee S^{11}$.
\end{exams}

\section{Complexity of resonances}  \label{s1.6}

The main idea of all our previous computations was to find, for
a~given $n$-cut $S$, a~partition $\pi\in P(n)$, such that
$\lspan(S\setminus(\pi\dar S))\neq S$. Intuitively speaking, shrinking
the substratum corresponding to $\tilde\pi S$, where
$\un(\tilde\pi)=\pi$, {\it essentially} reduces the set of linear
identities in~$S$. It is easy to construct examples when such $\pi$
does not exist, e.g., Example~\ref{exsmres}(4).

These observations lead us to introduce a~formal notion of complexity
of a~resonance.

\begin{df} \label{df6.1}  $\,$

\nin
1) For $S\in\catr_n$, the {\bf complexity} of $S$ is denoted $c(S)$
and is defined by:
\begin{equation}
  \label{eq:6.1}
c(S)=\min\{|\Pi|\,|\,\Pi\subseteq P(n),\lspan(S\setminus(\Pi\dar S))\neq S\}.
\end{equation}

\nin
2) We define the complexity of an $n$-resonance to be the complexity
of one of its representing cuts. Clearly, it does not depend
on the choice of the representative.
\end{df}

\begin{rem} The number $c(S)$ would not change if we required the
  partitions in $\Pi$ to have one block of size 2, and all other
  blocks of size~1.
\end{rem}


The higher is the complexity of a resonance $[S]$, the less it is
likely that one can succeed with analyzing its topological structure
using the method decribed in this chapter. This is because one would
need to take a~quotient by a~union of $c([S])$ strata and it might be
difficult to get a~hold on the topology of that union.

\vskip4pt

We finish by constructing for an~arbitrary $n\in{\Bbb N}$, a~resonance
of complexity~$n$. Let $\lambda_n=(a_1,\dots,a_n,b_1,\dots,b_n)$, such
that $a_i,b_i\in{\Bbb N}$, for $i,j\in[n]$, and all other linear
identities among $a_i$'s and $b_i$'s with coefficients $\pm 1,0$ are
generated by such identities. In other words, the cut $S$
associated to $\lambda$ is equal to the set
\begin{equation}
  \label{eq:6.2}
  \Bigl\{(x_1,\dots,x_n,y_1,\dots,y_n)\in\{-1,0,1\}^{2n}\,\Big |\,
\sum_{i=1}^n y_i=0,\,\,x_i+y_i=0,\forall i\in[n]\Bigr\}.
\end{equation}

\nin
It is not difficult to construct such $\lambda_n$ directly:

\vskip4pt

\nin 1) Choose $a_1,\dots,a_n$, such that the only linear identities
with coefficients $\pm 1,0$ on the set $a_1,a_1,a_2,a_2,\dots,a_n,a_n$
are of the form $a_i=a_i$; in other words, there are no linear
identities with coefficients $\pm 2,\pm 1,0$ on the set
$a_1,\dots,a_n$. One example is provided by the choice $a_1=1$,
$a_2=3$, $\dots$, $a_n=3^{n-1}$.

\vskip4pt

\nin 2) Let $b_i=N+a_i$, for $i\in[n]$, where $N$ is sufficiently
large. As the proof of the Proposition~\ref{prop6.2} will show, it is
enough to choose $N>2\sum_{i=1}^n\lambda_i$. This bound is far from 
sharp, but it is sufficient for our purposes.

\begin{prop} \label{prop6.2}
  Let $S_n$ be the $n$-cut associated to the ordered sequence
of natural numbers $\lambda_n$ described above, then $c(S_n)=n$.
\end{prop}

\pr First, let us verify that the cut $S_n$ associated to
$\lambda_n$ is equal to the one described in~\eqref{eq:6.2}. Take
$(x_1,\dots,x_n,y_1,\dots,y_n)\in S_n$. 

Assume first that $\sum_{i=1}^n y_i\neq 0$. Then,
$(x_1,\dots,x_n,y_1,\dots,y_n)$ stands for the identity
\begin{equation}
  \label{eq:6.3}
  \sum_{i\in I_1}a_i+\sum_{j\in J_1}b_j=
\sum_{i\in I_2}a_i+\sum_{j\in J_2}b_j,
\end{equation}
such that $|J_1|\geq|J_2|+1$. This implies that $N$ is equal to some 
linear combination of $a_i$'s with coefficients $\pm 2,\pm 1,0$.
This leads to contradiction, since $N>2\sum_{i=1}^n\lambda_i$.

Thus, we know that $\sum_{i=1}^n y_i=0$. Cancelling $N\cdot|J_1|$ out
of~\eqref{eq:6.3} we get an~identity with coefficients $\pm 2,\pm 1,0$
on the set $a_1,\dots,a_n$. By the choice of $a_i$'s, this identity
must be trivial, which amounts exactly to saying that $x_i+y_i=0$, for
$i\in[n]$. 

Second, it is a~trivial observation that $c(S_n)\leq n$. Indeed, let
$\pi_i\in P(n)$ be a~partition with only one nonsingleton block
$(1,n+i)$, for $i\in[n]$. Then
$\lspan(S_n\setminus(\{\pi_1,\dots,\pi_n\}\dar S_n))\neq S_n$, since for any
$(x_1,\dots,x_n,y_1,\dots,y_n)\in S_n\setminus(\{\pi_1,\dots,\pi_n\}\dar S_n)$,
we have $x_1=0$.

Finally, let us see that $c(S_n)>n-1$. As we have remarked after the
Definition~ref{df6.1}, it is enough to consider the case when the
partitions of $\Pi$ have one block of size 2, and the rest are
singletons. Let us call the identity $a_i+b_j=a_j+b_i$ {\it the 
elementary identity indexed $(i,j)$}. 

From the definition of the closure operation $\dar$ it is clear that
an~elementary identity indexed $(i,j)$ is not in $S_n\setminus(\Pi\dar
S_n)$ if and only if the partition whose only nonsingleton block is
$(i,n+j)$ belongs to $\Pi$, or the partition whose only nonsingleton
block is $(j,n+i)$ belongs to $\Pi$. That is because the only reason
this identity would not be in $S_n\setminus(\Pi\dar S_n)$ would be
that one of these two partitions is in $\Pi\dar S_n$. But, if such
a~partition is in $\Pi\dar S_n$, then it must be in $\Pi$: moves (2)
of the Definition~\ref{df3.1} can never produce a~partition whose only
nonsingleton block has size~2, while the moves (3) of the
Definition~\ref{df3.1} may only interchange between partitions
$(i,n+j)$ and $(j,n+i)$ in our specific situation. Thus, we can
conclude that if $|\Pi|\leq n-1$, then at most $n-1$ elementary
identities are not in $S_n\setminus(\Pi\dar S_n)$.

Next, we note that for any distinct $i,j,k\in[n]$, the elementary
identities $(i,j)$ and $(j,k)$ imply the elementary identity $(i,k)$.
Let us now think of elementary identities as edges in a~complete graph
on $n$~vertices, $K_n$. Then, any set $M$ of elementary identities
corresponds to a~graph $G$ on $n$ vertices, and the collection of the
elementary identities which lie in the $\lspan M$ is encoded by the
{\it transitive closure} of $G$. It is a~well known combinatorial fact
that $K_n$ is $(n-1)$-connected, which means that removal of at most
$n-1$ edges from it leaves a~connected graph. Hence, if we remove at
most $n-1$ edges from $K_n$ and then take the transitive closure, we
get $K_n$ again. Thus, if $|\Pi|\leq n-1$, all elementary identities
lie in $\lspan(S_n\setminus(\Pi\dar S_n))$. Since the elementary
identities generate the whole $S_n$, we conclude that
$S_n=\lspan(S_n\setminus(\Pi\dar S_n))$, hence $c(S_n)>n-1$.
\qed


\clearemptydoublepage

\chapter{Incidence Combinatorics of Resolutions}
           
\section{The motivation for the abstract framework}\label{sect_intr}

In this chapter we introduce notions of {\it combinatorial blowups,
  building sets\/}, and {\it nested sets\/}, for an arbitrary
meet-semilattice.  The definitions are given on a~purely
order-theoretic level without any reference to geometry. This provides
a common abstract framework for the incidence combinatorics occurring
in at least two different situations in algebraic geometry: the
construction of De~Concini-Procesi models of subspace
arrangements~\cite{DP95}, and the resolution of singularities in toric
varieties.

  The various parts of this abstract framework have received 
different emphasis within different situations: while the notion of
combinatorial blowups clearly specializes to stellar subdivisions of 
defining fans in the context of toric varieties, building sets and 
nested sets were introduced in the context of model constructions 
by De~Concini  \& Procesi~\cite{DP95} (earlier and in a more special 
setting by Fulton~ \& MacPherson~\cite{FM94}), from where we adopt our 
terminology. This correspondence however is not complete: 
the building sets in \cite{DP95, FM94} are not canonical, they depend 
on the geometry, while ours do not. See Section~\ref{ssect_DPmodels} 
for further details.

  It was proved in \cite{DP95} that a sequence of blowups within
an arrangement  of complex linear subspaces leads from
the~intersection stratification of complex space given by the 
maximal subspaces of the arrangement to an arrangement model 
stratified by divisors with normal crossings. 
In the context of toric varieties, there exist many different
procedures for stellar subdivisions of a defining fan that result
in a simplicial fan, so-called simplicial resolutions.  

  The purpose of our Main Theorem \ref{thm_main} is to unify 
these two situations on the combinatorial level: a sequence of 
combinatorial blowups, performed on a~(combinatorial) building 
set in linear extension compatible order, transforms the initial
semilattice to a semilattice  where all intervals are boolean 
algebras, more precisely to the face poset of the corresponding
simplicial complex of nested sets. In particular, the structure of 
the resulting semilattice can be fully described by the initial 
data of nested sets. Both the formulation and the proof of our 
main theorem are purely combinatorial.

\vskip4pt

\nin
We sketch the content of this chapter:

\vskip4pt

\noindent
{\bf Section \ref{sect_buidg_nest_sets}.}
After providing some basic poset terminology, we define building sets
and nested sets for meet-semilattices in purely order-theoretic terms
and develop general structure theory for these notions.
 
\vskip4pt

\noindent  
{\bf Section \ref{sect_comb_blowups}.} We define combinatorial blowups
of meet-semilattices, and study their effect on building sets and
nested sets. The section contains our Main Theorem~\ref{thm_main}
which describes the result of blowing up the elements of a building
set in terms of the initial nested set complex.
   
\vskip4pt

\noindent  
The next two sections are devoted to relating our abstract framework
to two different contexts in algebraic geometry. 
 
\vskip4pt

\noindent 
{\bf Section~\ref{ssect_DPmodels}.} We briefly review the construction
of De~Concini-Procesi models for subspace arrangements. After that, we
show that the change of the incidence combinatorics of the
stratification in a single construction step is described by
a~combinatorial blowup of the semilattice of strata, and trace their
resolution procedure step-by-step.
 
\vskip4pt

\noindent 
{\bf Section~\ref{ssect_tv}.} We draw the connection to simplicial
resolutions of toric varieties: we recognize stellar subdivisions as
combinatorial blowups of the face posets of defining fans and discuss
the notions of building and nested sets in this context.


\section{Building sets and nested sets of meet-semilattices} 
\label{sect_buidg_nest_sets}

\subsection{Irreducible elements in posets} \label{ssect_pos_term}

We assume known the parts of the terminology of posets described in
Appendix~B. A poset is called {\it irreducible\/} if it is not
a~direct product of two other posets, both consisting of at least two
elements. For a poset $P$ with a unique minimal element~$\hat 0$, we
call $I(P)=\{x\,{\in}\,P\,|\,[\hat 0,x]\,\, \text{is irreducible}\,\}$
the {\em set of irreducible elements\/} in~$P$. In particular, the
minimal element~$\hat 0$ and all atoms of $P$ are irreducible elements
in~$P$. For $x\,{\in}\,P$, we call $D(x)= {\rm max}\,(I(P)_{\leq x})$
the {\em set of elementary divisors\/} of~$x$ -- a~term which is
explained by the following proposition:

\begin{prop} \label{prop_elem_div} 
Let $P$ be a poset with a unique minimal element~$\hat 0$.
For $x\,{\in}\,P$ there exists a unique finest decomposition of the 
interval~$[\hat 0,x]$ in~$P$ as a direct product, which is given 
by an isomorphism 
$\varphi_x^{\rm el}:\,  \prod_{j=1}^l\, [\hat 0,y_j] \,  
                         \stackrel{\cong}{\lra} \,  [\hat 0,x]
$,
with $\varphi_x^{\rm el}(\hat 0, \ldots, y_j, \ldots, \hat 0)\,{=}\, y_j$ 
for $j=1,\ldots,l$. The factors of this decomposition are the
intervals below the elementary divisors of $x$: $\{y_1,\ldots,y_l\}=D(x)$. 
\end{prop}

\pr 
Whenever a poset with a minimal element~$\hat 0$ is represented as a~direct 
product, all elements which have more than one coordinate different from 
$\hat 0$ are reducible. 
Hence, if $\prod_{j=1}^l[\hat 0,y_j]\,{\cong}\,[\hat 0,x]$, 
and the $y_j$ are irreducible for $j\,{=}\,1,\ldots,l$, then 
$\{y_1,\ldots,y_l\}\,{=}\, D(x)$.
\qed


\subsection{Building sets}\label{ssect_buildg}

\noindent 
 In this subsection we define the notion of building sets of
 a~semilattice and develop their structure theory.

\begin{df} \label{df_buildg}
Let $\cl$ be a semilattice. A subset $\cg$ in~$\cl$ is called 
a {\bf building set} of~$\cl$ if for any $x\,{\in}\,\cl$  and
{\rm max}$\, \cg_{\leq x}=\{x_1,\ldots,x_k\}$ there is an isomorphism
of posets
\begin{equation}\label{eq_buildg}
\varphi_x:\,\,\, \prod_{j=1}^k\,\,\, [\hat 0,x_j] \,\, 
                               \stackrel{\cong}{\lra}
                                \,\, [\hat 0,x]
\end{equation}
with $\varphi_x(\hat 0, \ldots, x_j, \ldots, \hat 0)\, = \, x_j$ 
for $j=1,\ldots, k$. 
We call $F(x)=\max \cg_{\leq x}$  the {\bf set of factors} of $x$ in $\cg$.
\end{df}

 The next proposition provides several equivalent conditions for a subset
 of $\cl$ to be a~building set.

\begin{prop} \label{prop_buildingsets}
For a semilattice $\cl$ and a subset $\cg$ of~$\cl$ the following 
are equivalent:
\begin{itemize}

\item[(1)] $\cg$ is a building set of $\cl$;
  
\item[(2)] $\cg\/{\supseteq}\,I(\cl)$, and for every $x\,{\in}\,\cl$
  with $D(x)\,{=}\,\{y_1,\ldots, y_l\}$ the elementary divisors
  of~$x$, there exists a~partition $\pi_x\,{=}\,\pi_1|\ldots|\pi_k$, of
  the set $[l]$, with blocks $\pi_t\,{=}\,\{i_1,\ldots, i_{|\pi_t|}\}$,
  for $t\in[k]$, such that the elements in the set ${\rm max}\,
  \cg_{\leq x}=$ $\{x_1,\ldots, x_k\}$ are of the form $ x_t\, =
  \, \varphi_x^{{\rm el}}(\hat 0, \ldots, \hat 0, y_{i_1},\hat 0,
  \ldots, \hat 0,y_{i_2},\hat 0, \ldots, \hat 0,$ 
$y_{i_{|\pi_t|}},\hat 0) $.

\noindent
Informally speaking, the factors of $x$ in $\cg$ are products of disjoint
sets of elementary divisors of~$x$.
\item[(3)] $\cg$ generates $\cl$ by $\vee$, and
for any $x\,{\in}\,\cl$, any $\{y,y_1,\dots,y_t\}\subseteq 
\max \cg_{\leq x}$, and  $z\,{\in}\,\cl$ with $z\,{<}\,y$, we have 
$
  \cg_{\leq y}\cap \cg_{\leq z\vee y_1\vee\dots\vee y_t}= \cg_{\leq z}
$.

\item[(4)] $\cg$ generates $\cl$ by $\vee$, and for any $x\,{\in}\,\cl$, 
any $\{y,y_1,\dots,y_t\}\subseteq\max \cg_{\leq x}$,  and $z\,{\in}\,\cl$ with 
$z\,{<}\,y$,  the following two conditions are satisfied:
\[ 
\begin{array}{rll}
i)  & \cg_{\leq y}\cap \cg_{\leq y_1\vee\dots\vee y_t}=\{\hat 0\}
&\text{``disjointness,''}\\
i\!i) &z\vee y_1\vee\dots\vee y_t<y\vee y_1\vee\dots\vee y_t
&\text{``necessity.''} 
\end{array} 
\]
\end{itemize}
\end{prop}

\pr 

\noindent
\underline{(1)$\Rightarrow$(2)}: That $\cg$ contains $I(\cl)$ follows directly
from the definition of building sets. We have the following isomorphisms:
$\varphi_x:\,\prod_{j=1}^k\, [\hat 0,x_j]\,{\lra}\,[\hat 0,x]$
by the building set property, and
$\varphi_{x_j}^{\rm el}:\,\prod_{y\in D(x_j)}\, [\hat 0,y]
   \,{\lra}\,
[\hat 0,x_j]$
for $j\,{=}\, 1,\ldots, k$ by Proposition~\ref{prop_elem_div}.
The composition 
$\varphi_x\,\circ\, (\prod_{j=1}^k\,\varphi_{x_j}^{\rm el})$
yields the finest decomposition $\varphi_{x}^{\rm el}$ of~$[\hat 0,x]$.
Thus, $D(x)\,{=}\,\uplus_{j=1}^k\, D(x_j)$, which gives the partition 
described in~(2).


\noindent
\underline{(2)$\Rightarrow$(1)}: The decomposition of the interval
$[\hat 0,x]$ into intervals below the elements in $\max\cg_{\leq x}$
follows from Proposition~\ref{prop_elem_div} by assembling
factors~$[\hat 0, y_j]$ with maximal elements indexed by elements from
the same block of the partition~$\pi_x$ into one factor.

\noindent
\underline{(1)$\Rightarrow$(3)}: (3) is a direct consequence of~$[\hat 0,x]$
decomposing into a direct product of the form described in the definition of
building sets. 

\noindent
\underline{(3)$\Rightarrow$(4)}: $i)$ follows by setting $z=\hat 0$ 
in~(3).
Equality in $i\!i)$ implies with~(3) that $\cg_{\leq y}=\cg_{\leq z}$,
in particular, $y\in \cg_{\leq z}$ -- a contradiction to $z<y$.

\noindent
\underline{(4)$\Rightarrow$(1)}: For $x\in \cl$ and max$\,\cg_{\leq x}\, = \, 
\{x_1,\ldots, x_k\}$ consider the poset map
\[
\phi:\,\,\, \prod_{j=1}^k\,\,\, [\hat 0,x_j] \,\, \lra \,\, [\hat 0,x] \, ,
\quad
      (\alpha_1, \ldots, \alpha_k)\,\,  \lmt \,\,  
                               \alpha_1 \vee \ldots  \vee \alpha_k\, .
\] 

\noindent
i) $\phi$ {\em is surjective\/}: For $\hat 0\,{\neq}\, y \,{\leq}\,
x$, let max$\,\cg_{\leq y}=\{y_1,\ldots,y_t\}$. First, we have
$\bigvee_{i=1}^t y_i\,{=}\,y$, since $\cg$ generates $\cl$ by $\vee$.
Second, define $\gamma_j\,{=}\, \bigvee_{y_i\in S_j}y_i$, with
$S_j=$ $({\rm max}\,\cg_{\leq y})\,{\cap}\,\cg_{\leq x_j}$, for
$j\,{=}\,1,\ldots,k$.  Clearly, $\gamma_j\,{\in}\,[\hat 0,x_j]$, and
$\cup_{j=1}^k S_j\,{=}\, {\rm max}\,\cg_{\leq y}$, since $\cg_{\leq
  y}\,{\subseteq}\, \cg_{\leq x}$. Hence, $\phi(\gamma_1,\ldots,
\gamma_k)=\bigvee_{i=1}^t y_i=y$.


\noindent 
ii) $\phi$ {\em is injective\/}: a) Assume
$\phi(\alpha_1,\ldots,\alpha_k)\, {=}\, \phi(\beta_1,\ldots,\beta_k)=
y \neq x$, and let max$\,\cg_{\leq y}=\{y_1,\ldots,y_t\}$. By
induction on the number of elements in~$[\hat 0,x]$ we can assume that
$[\hat 0,y]$ decomposes as a direct product $ [\hat 0,y] \, \cong \,
\prod_{i=1}^t\, [\hat 0,y_i] $. Moreover, the subsets $S_j$ of
max$\,\cg_{\leq y}$ defined in~i) actually partition max$\,\cg_{\leq
  y}$ as follows from the disjointness property applied to pairwise
intersections of the $\,\cg_{\leq x_j}$. Thus, $ [\hat 0,y] \, \cong
\, \prod_{j=1}^k\, [\hat 0,\gamma_j] $, with elements
$\gamma_j\,{\in}\, [\hat 0,x_j]$ as above, and it follows that
$\alpha_j=\beta_j=\gamma_j$ for $j=1,\ldots,k$. \newline b) Assume
that $\phi(\alpha_1,\ldots,\alpha_k)\, {=}\,
\phi(\beta_1,\ldots,\beta_k)\,{=}\,x$. By the necessity property it
follows that $\alpha_j=\beta_j=x_j$ for $j=1,\ldots,k$.  \qed

\begin{rem} \label{rem_indpce_of_joins}
The definition of building sets and of irreducible elements, 
as well as the characterization of building sets
in Proposition~\ref{prop_buildingsets}~(2), are independent of the existence 
of a join operation and can be formulated for any poset with a unique minimal 
element. 
\end{rem}

\noindent
We gather a few important properties 
of building sets. 

\begin{prop}\label{prop_buildingsets_properties}
For a building set $\cg$ of $\cl$, the following holds:
\begin{itemize}
\item [(1)] Let $x\in\cl$, $F(x)=\{x_1,\dots,x_k\}$ the set of factors of~$x$ 
in~$\cg$, and $\hat 0\,{\neq}\,y\,{\in}\,\cg$ with $y\,{\leq}\,x$. Then 
there exists a unique  $j\,{\in}\, \{1,\ldots,k\}$ such that 
$y\,{\leq}\,x_j$; i.e., $F(x)={\rm max}\,\cg_{\leq x}$ induces a partition of 
$\cg_{\leq x}\,{\setminus}\,\{\hat 0\}$.
\item[(2)] For $x\,{\in}\,\cl$ and $x_0\,{\in}\,F(x)$,
\[
\bigvee\, (F(x)\,{\setminus}\,\{x_0\}) \,\, < \,\, \bigvee\, F(x)\, =x\, ,
\]
i.e., each factor of~$x$ in $\cg$ is needed to generate~$x$.
\item [(3)] If $h_1,\dots,h_k$ in $\cg$ are such that 
$(h_i,\bigvee_{j=1}^k h_j] \cap \cg=\emptyset$ for $i=1,\dots,k$, 
then $F(\bigvee_{j=1}^k h_j)\, = \, \{h_1,\ldots, h_k\}$.
\end{itemize}
\end{prop}

\pr
(1) is a consequence of Proposition~\ref{prop_buildingsets}~(4){\em i}), 
as was noted already in the proof of (4)$\Rightarrow$(1), part ii)~a),
in the previous proposition. 
Taking the full set of factors and setting~$z\,{=}\,\hat 0$ in 
Proposition~\ref{prop_buildingsets}~(4){\em i$\!$i}), yields~(2). For~(3) note
that $\{h_1,\ldots,h_k\}\,{\subseteq}\,F(\bigvee_{j=1}^k h_j)$ by assumption.
If $\{h_1,\ldots, h_k\}$ were not the complete set of factors, 
we would obtain a contradiction to~(2).
\qed


\subsection{Nested sets}

\noindent 
  In this subsection we define the notion of nested subsets of
a~building set of a~semilattice and prove some of their properties.

\begin{df} \label{df_nested}
Let $\cl$ be a semilattice and $\cg$ a building set of $\cl$.
A subset $N$ in $\cg$ is called~{\bf nested} if, for any 
set of incomparable elements 
$x_1,\dots,x_t$ in $ N$ of cardinality at least two, 
the join $x_1\vee\dots\vee x_t$ exists and does not belong to $\cg$.
The nested sets in $\cg$ form an abstract simplicial complex, denoted 
$\cn(\cg)$. 
\end{df}

Note that the elements of $\cg$ are the vertices of the complex 
of nested sets $\cn(\cg)$. Moreover, 
the order complex of $\cg$ is a subcomplex of $\cn(\cg)$, since
linearly ordered subsets of~$\cg$ are nested.

\begin{prop}\label{prop_nested}
For a given semilattice $\cl$ and a subset $N$ of a building set $\cg$ 
of~$\cl$, the following are equivalent:
  \begin{enumerate}
  \item [(1)] $N$ is nested.
  \item [(2)] Whenever $x_1,\dots,x_t$ are noncomparable elements in~$N$, 
the join $x_1\vee\dots\vee x_t$ exists, and 
$F(x_1\vee\dots\vee x_t)=\{x_1,\dots,x_t\}$.
  \item [(3)] There exists a chain $C\subseteq\cl$, such that
$N=\bigcup_{x\in C}F(x)$. 
  \item [(4)] $N\in\Lambda$, where $\Lambda$ is the maximal subset of
  $2^{\cg}$, for which the following three conditions are satisfied:
    \begin{enumerate}
    \item [(o)] $\emptyset\in\Lambda$, and $\{g\}\in\Lambda$, for $g\in \cg$;
    \item [(i)] if $N\in\Lambda$ and $x\in\max N$, then $N_{<x}\in\Lambda$;
    \item [(ii)] if $N\in\Lambda$, then $\max N=F(\bigvee \, \max N)$.
    \end{enumerate}
  \end{enumerate}
\end{prop}
\pr 

\vskip4pt

\noindent
\underline{(1)$\Rightarrow$(2)}: Let $N$ be a nested set, and let
$M\,{=}\,\{x_1,\dots,x_t\}\,{\subseteq}\,N$ be a~set of incomparable
elements with $\bigvee_{i=1}^t x_i\,{\not \in}\,\cg$.  We can assume
that for some $x_j$ we have $(x_j,\bigvee_{i=1}^t
x_i]\,{\cap}\,\cg\,{\neq}\,\emptyset$, otherwise the claim follows by
Proposition~\ref{prop_buildingsets_properties}~(3). Without loss of
generality, we may assume that there exists an element $y$, such that
$y\,{\in}\,(x_1,\bigvee_{i=1}^t x_i]\,{\cap}\,\cg$, and
$y\,{\in}\,\max \cg_{\leq \bigvee M}$.  Define
$M'\,{=}\,\{x_1,\dots,x_t\}\,{\cap}\,\cg_{\leq y}\,{=}$ \linebreak
$\{x_1\,{=}\,x_{j_0},x_{j_1},\dots,x_{j_k}\}$, and
$z\,{=}\,\bigvee_{l=0}^k x_{j_l}$.  Since
$M'=\{x_{j_0},x_{j_1},\dots,x_{j_k}\}$ is nested (it is a subset
of $N$), we have the strict inequality $z\,{<}\,y$. Furthermore,
\[
\bigvee_{i=1}^t x_i  \, = \, 
z\vee\bigvee (M\setminus M')\, \leq \, 
z\vee\bigvee(\max G_{\leq\bigvee M}\setminus\{y\})  \, < \, 
\bigvee_{i=1}^t x_i\,,
\]
where the first inequality follows from 
Proposition~\ref{prop_buildingsets_properties}~(1) and the second 
inequality from Proposition~\ref{prop_buildingsets_properties}~(2). 
We thus arrive to a contradiction, which finishes the proof.

\vskip4pt

\noindent
\underline{(2)$\Rightarrow$(1)}: Obvious.

\vskip4pt

\noindent
\underline{(2)$\Rightarrow$(3)}: Let $N$ be a~set satisfying condition~(2). 
Fix a~particular linear extension $\{x_1,\dots,x_k\}$ on the partial order 
of~$N$, and define
$\alpha_j\,{=}\,x_1\,{\vee}\,\dots\,{\vee}\, x_j$, for $j\,{=}\,1,\dots,k$.
By (2) we have $F(\alpha_j)\,{=}\,\max\{x_1,\dots,x_j\}$, and therefore
$x_j\,{\in}\, F(\alpha_j)$ and $x_{j+1}\,{\not\in}\, F(\alpha_j)$ for
$j\,{=}\,1,\dots,k$. Hence, the $\alpha_j$'s are different and form
a~chain $C=\alpha_1<\alpha_2<\dots<\alpha_k$. By construction, 
$N=\bigcup_{x\in C}F(x)$.

\vskip4pt

\noindent
\underline{(1),(2)$\Rightarrow$(4)}: 
Let $N$ be a nested set, we shall prove that $N\in\Lambda$ 
by induction on the size of $N$:
\begin{enumerate}
\item if $|N|=0$, then $N\in\Lambda$ by condition (o);
\item if $|N|\geq 1$, then $\max N=F(\bigvee\,\max N)$ by 
condition~(2). Furthermore, since $|N_{<x}|<|N|$, and $N_{<x}$ 
is nested (it is a subset of $N$), $N_{<x}\in\Lambda$ by induction. 
Hence $N\in\Lambda$.
\end{enumerate}


\noindent
\underline{(3)$\Rightarrow$(1)}: 
Let $C\,{=}\,(\alpha_1\,{<}\,\dots\,{<}\,\alpha_k)$ be a chain in~$\cl$ and 
$N\,{=}\,\bigcup_{x\in C}F(x)$. 
Let $N'\,{=}\,\{x_1,\ldots,x_t\}\,{\subseteq}\, N$, $t\,{\geq}\,2$, be an
antichain in~$N$, and~$s$ the maximal index in~$C$ such that 
$N'\,{\cap}\,F(\alpha_s)\,{\neq}\,\emptyset$. In particular, 
$N'\,{\cap}\,F(\alpha_s)\,{\neq}\,\{\alpha_s\}$ due to $|N'|>1$ and
$N'$ being an antichain.

Let $y\,{\in}\,N'\,{\cap}\,F(\alpha_s)$. 
If $|N'\,{\cap}\,F(\alpha_s)|\,{>}\,1$,
\[
     y\, < \, \bigvee (N'\,{\cap}\,F(\alpha_s)) \, \leq \, 
              \bigvee \, N'      \, \leq \, \alpha_s\, ,  
\]
where the strict inequality is a consequence of the necessity property
for building sets. Thus, $ \bigvee  N'\,{\not \in}\, \cg$. 
If $|N'\,{\cap}\,F(\alpha_s)|\,{=}\,1$, we have 
$y\,{<}\,\bigvee N'\,{\leq}\,\alpha_s$,
due to $N'$ being an antichain with $|N'|\,{>}\,1$, and
again $\bigvee  N'\,{\not \in}\, \cg$.  

\vskip4pt

\noindent
\underline{(4)$\Rightarrow$(3)}: 
We need the following fact:

\vskip4pt
\noindent
{\bf Fact.} {\it If there exist elements $x_1,\dots,x_t$ and
  $y_1,\dots,y_k$ in~$\cl$, such that $x_t\,{>}\,y_j$, for
  $j\,{=}\,1,\dots,k$, and $F(\bigvee_{i=1}^t
  x_i)\,{=}\,\{x_1,\dots,x_t\}$, and $F(\bigvee_{j=1}^k
  y_j)\,{=}\,\{y_1,\dots,y_k\}$, then $F(x_1\vee\dots\vee x_{t-1}\vee
  y_1\vee\dots\vee y_k)\,{=}\, \{x_1,\dots,x_{t-1},y_1,\dots,y_k\}.$}

\vskip4pt Once the fact above is proved, one can derive (3) as
follows. For $N\in\Lambda$ we shall form a chain
$C\,{=}\,(\alpha_1\,{<}\,\dots\,{<}\,\alpha_{|N|})$ such that
$N\,{=}\,\bigcup_{i=1}^{|N|} F(\alpha_i)$. First, choose a~linear
extension $\{x_1,\dots,x_t\}$ of $N$. Then, set
$\alpha_t\,{=}\,\bigvee \, \max N$, $\alpha_{t-1}=$ $\bigvee \, \max
(N\,{\setminus}\,\{x_t\})$, $\alpha_{t-2}\,{=}\,\bigvee \, \max
(N\,{\setminus}\,\{x_t, x_{t-1}\})$, and so on. By (4)(ii), we have
$F(\alpha_t)\,{=}\,\max N$. Applying (4)(i) to $x_t\,{\in}\,\max N$,
and (4)(ii) to $N_{<x_t}$, we obtain $F(\bigvee\, \max
N_{<x_t})\,{=}\,\max N_{<x_t}$. Taking into account the fact above, we
conclude that $F(\alpha_{t-1})\,{=}\,\max (N\,{\setminus}\,\{x_t\})$,
and, using the same argument iteratively, we arrive to
$N\,{=}\,\bigcup_{i=1}^t F(\alpha_i)$.

\vskip3pt
\noindent
{\bf Proof of the fact.}
Set $\alpha\,{=}\,x_1\,{\vee}\,\dots\,{\vee}\, x_{t-1}\,{\vee}\, 
                   y_1\,{\vee}\,\dots\,{\vee}\, y_k$.
Since~$\alpha\,{\leq}\,\bigvee_{i=1}^t\,x_i$, the factors of~$\alpha$
can be partitioned into groups of elements below the $x_i$ for $i\,{=}\,
1,\ldots,t$, by Proposition~\ref{prop_buildingsets_properties}~(1).
Since $x_i\,{\leq}\,\alpha$ for $i\,{=}\,1,\ldots,t{-}1$, we obtain
$F(\alpha)\,{=}\,\{x_1,\ldots, x_{t-1}, \gamma_1,\ldots, \gamma_{m}\}$
with $\gamma_j\,{\leq}\,x_t$ for $j=1,\ldots,m$.

Again using Proposition~\ref{prop_buildingsets_properties}~(1), the
$y_1,\dots, y_k$ can be partitioned into groups below the factors
$\gamma_j$ for $j\,{=}\,1,\ldots,m$. The occurrence of one strict
inequality $\bigvee\,\{y_l\,|\,y_l\,{\leq}\,\gamma_j\}\,{<} \,
\gamma_j$, for some $j\in[m]$, yields a contradiction to
$$\alpha\,{=}\, \bigvee_{i=1}^{t-1} x_i \,{\vee}\bigvee_{j=1}^{k}
y_j{=}\, \bigvee_{i=1}^{t-1} x_i \,{\vee}\,\bigvee_{j=1}^{m}
\gamma_j ,$$ due to the necessity property of building sets. Moreover,
since the $y_i$ are factors themselves, joins of more than two of the
$y_i$'s are not elements of~$\cg$. Thus, $y_i\,{=}\,\gamma_i$, for
$i\,{=}\,1,\ldots,k{=}m$, as claimed. \qed

\section{Sequences of combinatorial blowups} 
\label{sect_comb_blowups}

\noindent 
We introduce the notion of a combinatorial blowup of an~element
in a~semilattice and prove that the set of semilattices is closed
under this operation.

\subsection{Combinatorial blowups}

\begin{df}
For a semilattice $\cl$ and an element $\al\in\cl$ we define 
a poset $\bl_\al\cl$, the {\bf combinatorial blowup of $\cl$ at $\al$},  
as follows:
\begin{itemize}
\item[$\circ$] elements of $\bl_\al\cl$:
  \begin{enumerate}
  \item[(1)] $y\in\cl$, such that $y\not\geq\al$;
  \item[(2)] $[\al,y]$, for $y\in\cl$, such that $y\not\geq\al$
             and $(y\vee\al)_\cl$ exists \newline
             (in particular, $[\al,\hat 0]$ can be
             thought of as the result of blowing up $\al$);
  \end{enumerate}
\item[$\circ$] order relations in $\bl_\al\cl$: 
  \begin{enumerate}
  \item[(1)] $y>z$ in $\bl_\al\cl$ if $y>z$ in $\cl$;
  \item[(2)] $[\al,y]>[\al,z]$ in $\bl_\al\cl$ if $y>z$ in $\cl$;
  \item[(3)] $[\al,y]>z$ in $\bl_\al\cl$ if $y\geq z$ in $\cl$;
  \end{enumerate}\nopagebreak[4]
  where in all three cases $y,z\not\geq\al$. 
\end{itemize}
\end{df}

Note that the atoms in $\bl_\al\cl$ are the atoms of~$\cl$
together with the element $[\alpha,\hat 0]$.
It is easy, albeit tedious, to check that the class of 
(meet-)semilattices is closed under combinatorial blowups.

\begin{lm}
Let $\cl$ be a semilattice and $\al\in\cl$, then
$\bl_\al\cl$ is a semilattice.
\end{lm}
\pr The joins in $\bl_\al\cl$ are defined by the rule
\begin{align*}
([\al,y]\vee[\al,z])_{\bl_\al\cl}&=[\al,(y\vee z)_\cl],\\
  ([\al,y]\vee z)_{\bl_\al\cl}   &=[\al,(y\vee z)_\cl],\\
  (y\vee z)_{\bl_\al\cl}         &=(y\vee z)_\cl,
\end{align*}
which is applicable only if $(y\vee z)_\cl$ exists, otherwise
the corresponding joins in $\bl_\al\cl$ do not exist. Also,
the first and second formulae are applicable only in the case 
$(y\vee z)_\cl\not\geq\al$, otherwise the corresponding
joins do not exist. 

We check this by considering four possible cases
separately:
$$
\begin{cases}
  [\al,x]\geq[\al,y]\\
 {}[\al,x]\geq[\al,z]
\end{cases}
\Leftrightarrow
\begin{cases}
  x\geq y\\
  x\geq z
\end{cases}
\Leftrightarrow
x\geq(y\vee z)_\cl
\Leftrightarrow
[\al,x]\geq[\al,(y\vee z)_\cl].
$$
$$
\begin{cases}
  [\al,x]\geq[\al,y]\\
 {}[\al,x]\geq z
\end{cases}
\Leftrightarrow
\begin{cases}
  x\geq y\\
  x\geq z
\end{cases}
\Leftrightarrow
x\geq(y\vee z)_\cl
\Leftrightarrow
[\al,x]\geq[\al,(y\vee z)_\cl].
$$
$$
\begin{cases}
  x\geq_{\bl_\al\cl} y\\
  x\geq_{\bl_\al\cl} z\\
  x\not\geq\al
\end{cases}
\Leftrightarrow
\begin{cases}
  x\geq_\cl y\\
  x\geq_\cl z\\
  x\not\geq\al
\end{cases}
\Leftrightarrow
\begin{cases}
x\geq(y\vee z)_\cl\\
x\not\geq\al
\end{cases}
\Rightarrow
\begin{cases}
y\vee z\not\geq\al\\
x\geq(y\vee z)_\cl.
\end{cases}
$$
$$
\begin{cases}
  [\al,x]\geq y\\
{}[\al,x]\geq z\\
  y,z\not\geq\al
\end{cases}
\hskip-3pt
\Leftrightarrow
\begin{cases}
  x\geq_\cl y\\
  x\geq_\cl z\\
  x,y,z\not\geq\al
\end{cases}
\hskip-3pt\Leftrightarrow
\begin{cases}
x\geq(y\vee z)_\cl\\
x,(y\vee z)_\cl\not\geq\al
\end{cases}
\hskip-5pt\Rightarrow
[\al,x]\geq(y\vee z)_{\bl_\al\cl}.
$$

 \noindent
Observe that it is possible that $(x\vee y)_\cl$ exists, while 
$(x\vee y)_{\bl_\al\cl}$ does not. \qed


\subsection{Blowing up building sets} %

\noindent
In this subsection we prove that if one combinatorially blows up 
a~building set of a~semilattice in any chosen linear extension order,
then one ends up with the face poset of the simplicial complex of nested
sets of this building set.
The following proposition provides the essential step for the proof. 

\begin{prop}\label{prop_singleblow}
  Let $\cl$ be a semilattice, $\cg$ a building set of~$\cl$, and 
$\al\,{\in}\,\max \cg$. Then, $\wti \cg=(\cg\sm\{\al\})\cup\{[\al,\hat 0]\}$
is a building set of $\bl_\al\cl$.
Furthermore, the nested subsets of $\wti \cg$ are precisely
the nested subsets of $\cg$ with $\al$ replaced by $[\al,\hat 0]$.
\end{prop}

\pr 
It is easy to see that $\wti \cg$ is a building set of $\bl_\al\cl$.
Indeed, given $x\in\cl\sm\cl_{\geq\al}$, \eqref{eq_buildg} is obvious for 
$x\in\bl_\al\cl$, and, if $(x\vee \al)_\cl$ exists, it follows for 
$[\al,x]\in\bl_\al\cl$ from the identity
\[
[\hat 0,[\al,x]]_{\bl_\al\cl}\, =\, [\hat 0,x]_{\bl_\al\cl}\times B_1\, ,
\]
where $B_1$ is the subposet consisting of the two comparable elements
$\hat 0\,{<}\,[\al,\hat 0]$.

Let us now see that the sets of nested subsets of $\cg$ and 
$\wti \cg$ are the same when replacing $\al$ by $[\al,\hat 0]$:
 
Let $N$ be a nested set in $\cg$, not containing $\al$. For incomparable
elements $x_1,\ldots,x_t$ in~$N$, 
$\bigvee_{i=1}^t\,x_i\,{\not \geq}\,\al$,
since otherwise we had 
$$\al\,{\in}\, \max \cg_{\leq \bigvee x_i} =
F(\bigvee_{i=1}^tx_i)\,{=}\,
\{x_1,\ldots,x_t\}$$
by Proposition~\ref{prop_nested}(2). Thus, $\bigvee_{i=1}^t x_i$ exists in
$\bl_{\al} \cl$ and $\bigvee_{i=1}^t x_i \,{\not\in}\,  
\wti \cg$. Hence, $N$ is nested in $\wti \cg$.
A nested subset in  $\wti \cg$ not containing $[\al,\hat 0]$ is obviously
nested in $\cg$.

Let now $N$ be nested in~$\cg$ containing $\alpha$, and set 
$\wti N\,{=}\,(N\sm\{\al\})\cup\{[\al,\hat 0]\}$. Subsets of incomparable 
elements in $\wti N$ not containing $[\al,\hat 0]$ can be dealt with as above.
Thus assume that $[\al,\hat 0],x_1,\ldots,x_t$ are incomparable in $\wti N$.
Then, $x_1,\ldots,x_t$ are incomparable in the nested set~$N$, and, as above,
we conclude that $\bigvee_{i=1}^t\,x_i$ exists and 
$\bigvee_{i=1}^t\,x_i\,{\not \geq}\,\al$. Moreover, $\alpha \vee 
\bigvee_{i=1}^t\,x_i$ exists in $\cl$ (joins of nested sets always exist!),
thus, 
$[\al,\bigvee_{i=1}^t\,x_i]\,{=}\,[\al,\hat 0]\vee \bigvee_{i=1}^t\,x_i$
exists in $\bl_{\al}\cl$ and is obviously not contained in $\wti \cg$.
We conclude that $\wti N$ is nested in $\wti \cg$. 

Vice versa, let $\wti N$ be nested in~$\wti \cg$ containing $[\al,\hat 0]$,
and set $N\,{=}\,(\wti N\sm \{[\al,\hat 0]\})\cup\{\al\}$.
Again it suffices to consider subsets of incomparable elements
$\al,x_1,\ldots,x_t$ in $N$. With  $[\al,\hat 0],x_1,\ldots,x_t$ 
incomparable in $\wti N$, $[\al,\hat 0]\vee \bigvee_{i=1}^t\,x_i\,{=}\,
[\al,\bigvee_{i=1}^t\,x_i]$ exists in $\bl_{\al}\cl$, thus 
$\al \vee \bigvee_{i=1}^t\,x_i$ exists in~$\cl$. Incomparability
implies that $\al \vee \bigvee_{i=1}^t\,x_i>\al$, and thus
$\al \vee \bigvee_{i=1}^t\,x_i \not \in \cg$. We conclude that~$N$ 
is nested in~$\cg$.
\qed




\medskip
By iterating the combinatorial blowup described in 
Proposition~\ref{prop_singleblow} through all of~$\cg$, 
we obtain the following theorem, which serves as a motivation 
for the entire development.

\begin{thm} \label{thm_main}
Let $\cl$ be a semilattice and $\cg$ a building set of~$\cl$
with some chosen linear extension: $\cg=\{G_1,\dots,G_t\}$, with  
$G_i>G_j$ implying $i<j$. Let $\bl_k\cl$ denote the result 
of subsequent blowups $\bl_{G_k}(\bl_{G_{k-1}}(\dots\bl_{G_1}\cl))$.
Then the final semilattice $\bl_t\cl$ is equal to the face poset of 
the simplicial complex $\cn(\cg)$.
\end{thm}

\pr
The building set~$\cg_t$ of $\bl_t\cl$ that results from iterated 
application of Proposition~\ref{prop_singleblow}  obviously is the set
of atoms ${\mathfrak A}$ in $\bl_t\cl$. Every element $x\in\bl_t\cl$ is the join 
of atoms below it: $x=\bigvee {\mathfrak A}_{\leq x}$. The subset 
${\mathfrak A}_{\leq x}$
of $\cg_t$ is nested, in particular, it is the set of factors of~$x$ 
in $\bl_t\cl$ with respect to $\cg_t$ (Proposition~\ref{prop_nested}(2)). 
Proposition~\ref{prop_buildingsets_properties}(2) implies that the
interval~$[\hat 0,x]$ in $\bl_t\cl$ is boolean. We conclude that
$\bl_t\cl$ is the face poset of a simplicial complex with faces in 
one-to-one correspondence with the nested sets in $\cg_t$, which in turn
correspond to the nested sets in~$\cg$ by Proposition~\ref{prop_singleblow}.
\qed



\section{De~Concini-Procesi models of subspace arrangements}
\label{ssect_DPmodels}

\subsection{Building sets for subspace arrangements}

Let $\ca=\{A_1,\ldots,A_n\}$ be an arrangement of linear subspaces in
$\dc^d$. Much effort has been spent on describing the cohomology of
the complement $\cm(\ca)=\dc^d\,{\setminus}\,\bigcup \ca$ of such an
arrangement and, in particular, on answering the question whether the
cohomology algebra is completely determined by the combinatorial data
of the arrangement. Here, combinatorial data is understood as the
lattice~$\cl(\ca)$ of intersections of subspaces of~$\ca$ ordered by
reverse inclusion together with the complex codimensions of the
intersections. A major step towards the solution of this problem (for
a complete answer see~\cite{DGM00, dLS01}) was the construction of
smooth models for the complement~$\cm(\ca)$ by De~Concini~ \&
Procesi~\cite{DP95} that allowed for an explicit description of
rational models for~$\cm(\ca)$ following~\cite{Mo78}.  The
De~Concini-Procesi models for arrangements in turn are one instance in
a sequence of model constructions reaching from compactifications of
symmetric spaces~\cite{DP83,DP85}, over the Fulton-MacPherson
compactifications of configuration spaces~\cite{FM94} to the general
framework of wonderful conical compactifications proposed by
MacPherson~ \& Procesi~\cite{MP98}.

Given a complex subspace arrangement~$\ca$ in~$\dc^d$, De~Concini  \& Procesi
describe a smooth irreducible variety~$Y$ together with a proper map 
$\pi:\, Y\, \lra \, \dc^d$ such that $\pi$ is isomorphism over $\cm(\ca)$,
and the complement of the preimage of $\cm(\ca)$ is a union of irreducible 
divisors with normal crossings in~$Y$. The model~$Y$ can be
constructed by a~sequence of blowups of smooth subvarieties that is 
prescribed by the stratification of complex space induced by the arrangement. 



In order to enumerate the strata in the intersection stratification of~$Y$
given by the irreducible divisors, De~Concini  \& Procesi introduced the 
notions of building sets, nested sets and irreducible elements as follows:  

\begin{df} \label{df_DPnotions}
{\rm (\cite[\S 2]{DP95})}
  Let $\cl(\ca)$ be the intersection lattice of an arrangement~$\ca$ 
of linear subspaces in a~finite dimensional complex vector space. 
Consider the lattice~$\cl(\ca)^*$ formed by the orthogonal 
complements of intersections ordered by inclusion. 
\begin{itemize}
\item[(1)] For $U\,{\in}\,\cl(\ca)^*$,
$U\, = \oplus_{i=1}^k\, U_i$ with $U_i\,{\in}\,\cl(\ca)^*$, is called a
{\bf decomposition\/} of~$U$ if for any $V\,{\subseteq}\,U$, 
$V\in \cl(\ca)^*$, $V\, = \oplus_{i=1}^k\, (U_i\cap V)$ and 
$U_i\cap V\,{\in}\,\cl(\ca)^*$ for $i=1,\ldots,k$. 
\item[(2)] Call  $U\,{\in}\,\cl(\ca)^*$ {\bf irreducible\/} if it does not
admit a non-trivial decomposition.
\item[(3)] A subset $\cg\,{\subseteq}\,\cl(\ca)^*$ is called a 
{\bf building set\/}
for~$\ca$ if for any $U\,{\in}\,\cl(\ca)^*$ and $G_1,\ldots, G_k$
maximal in~$\cg$ below~$U$,
$U\, = \oplus_{i=1}^k\, G_i$ is a decomposition (the $\cg$-decomposition)
of~$U$.
\item[(4)] A subset $\cs\,{\subseteq}\,\cg$ is called {\bf nested\/} if for
any set of non-comparable elements $U_1,\ldots,U_k$ in~$\cs$, 
$U\, = \oplus_{i=1}^k\, U_i$ is the $\cg$-decomposition of~$U$.
\end{itemize}
\end{df}

Note that $\cl(\ca)^*$ coincides with $\cl(\ca)$ as abstract lattices.
We will therefore talk about irreducible elements, building sets and 
nested sets in $\cl(\ca)$ without explicitly referring to the dual setting
of the preceding definition.

The notions of Definition~\ref{df_DPnotions} are in part based on the 
earlier notions introduced by 
Fulton~ \& MacPherson in~\cite{FM94} to study compactifications 
of configuration spaces. Our terminology is naturally adopted 
from~\cite{FM94,DP95}. Building sets and nested sets in the sense of 
De~Concini~ \& Procesi are building and nested sets for the intersection 
lattices of subspace arrangements in our combinatorial sense
(see Proposition~\ref{combgeomprop}~(1) below). However, there
are differences. The opposite is not true: A combinatorial building set 
for the intersection lattice of a subspace arrangement is not
necessarily a building set for this arrangement in the sense of 
De~Concini~ \& Procesi, neither
are irreducible elements in the sense of De~Concini~ \& Procesi 
irreducible in our sense. 

\begin{exam}
Consider the following arrangement of 3 subspaces 
in~$\dc^4$:
\[
A_1:\,\, z_4  =  0\, , \quad
A_2:\,\, z_1= z_2 =  0\, ,\quad
A_3:\,\, z_1= z_3 =  0 \,.
\]

\noindent
The intersection lattice is a boolean algebra on 3 elements.
$\{A_1,A_2,A_3\}\,{\subseteq}\,\cl(\ca)$  is a combinatorial building set, 
in fact, it is the set of irreducibles in the abstract lattice. However, the
minimal building set in the sense of De~Concini~ \& Procesi is 
$\{A_1,A_2,A_3, A_2\,{\cap}\,A_3\}$.
\end{exam}

The main difference can be formulated in the following way:
our constructions are order-theoretically canonical for a given 
semilattice. The set of combinatorial building sets, in particular 
the set of irreducible elements, depends only on the semilattice
itself and not on the geometry of the subspace arrangement which it 
encodes. See Proposition~\ref{combgeomprop} for the complete explanation.

\subsection{Local subspace arrangements}

 In order to trace the De~Concini-Procesi construction step by step
we need the more general notion of a~local subspace arrangement.

\begin{df} \label{lsadf}
  Let $M$ be a smooth complex $d$-dimensional manifold, and let $\ca$
  be a~union of finitely many smooth complex submanifolds of $M$ such
  that all non-empty intersections of submanifolds in~$\ca$ are
  connected smooth complex submanifolds. $\ca$ is called a~{\bf local
    subspace arrangement}~if for any $x\in\ca$ there exists an~open
  set $N$ in~$M$ with $x\,{\in}\,N$, a~subspace arrangement $\wti\ca$
  in ${\dc}^d$, and a~biholomorphic map $\phi:N\ra{\dc}^d$, such that
  $\phi(N\cap\ca)=\wti\ca$.
\end{df}

  Given a subspace arrangement $\ca$, the initial ambient space~$\dc^d$ 
of $\cm(\ca)$ carries a natural stratification  by the subspaces
of~$\ca$ and their intersections, the poset of strata being the 
intersection lattice $\cl(\ca)$ of the arrangement. For a~local
subspace arrangement $\ca=\{A_1,\dots,A_n\}$ in $M$ we again
consider the stratification of $M$ by all possible intersections
of the $A_i$'s, just like in the global case. The poset of strata
is also denoted by $\cl(\ca)$ and is called the intersection 
semilattice (it is a lattice if the intersection of all maximal
strata is nonempty).

\begin{df} \label{lsabsdf}
  Let $\ca$ be a local subspace arrangement and $\cl(\ca)$ its 
intersection semilattice. For $U\in\cl(\ca)$, $U_1,\dots,U_k\in\cl(\ca)$
are said to form a~{\bf decomposition} of $U$ if for any $x\in U$ there
exists an~open set $N$ with $x\,{\in}\,N$ and a~biholomorphic map 
$\phi:N\ra{\dc}^d$, such that $\phi(N\cap U_1),\dots,\phi(N\cap U_k)$
form a~decomposition of $\phi(N\cap U)$ in the sense of 
{\rm Definition~\ref{df_DPnotions}(1)}. 

  As in the global case, $\cg\subseteq\cl(\ca)$ is a~{\bf building 
set} for $\ca$ if for any $U\in\cl(\ca)$, the set of strata $\max\cg_{\leq U}$
gives a~decomposition of $U$.
\end{df}
 We shall refer to these building sets as {\em geometric\/} building sets.
The difference between combinatorial building sets and geometric
ones is contained in the dimension function as is explained in the 
following proposition.

\begin{prop} \label{combgeomprop} $\,$
  Let $\ca$ be a~local subspace arrangement with intersection 
semilattice~$\cl(\ca)$.
\begin{enumerate}
\item[(1)] If $\cg\subseteq\cl(\ca)$ is a~geometric building set of
$\ca$, then it is a~combinatorial building set.

\item[(2)] If $\cg\subseteq\cl(\ca)$ is a~combinatorial building set 
of $\cl(\ca)$, and for any $x\in\cl(\ca)$ the sum of codimensions of its
factors is equal to the codimension of $x$, then $\cg$ is
a~geometric building set.
\end{enumerate}
\end{prop}

\pr 
   In both cases it is enough to consider the case when $\ca$
is a~subspace arrangement.

(1) Consider~$\cg$ as a subset of~$\cl(\ca)^*$, then, for 
$U\,{\in}\,\cg$, the isomorphism $\varphi_U$ requested in 
Definition~\ref{df_buildg} is given by taking direct sums:
$$\varphi_U:\,\,\, \prod_{j=1}^k\,\,\, [\hat 0,G_j] \,\, 
\stackrel{\oplus_{j=1}^k}{\lra} \,\, [\hat 0,U]\, ,$$
where $G_1,\ldots,G_k$ are maximal in $\cg$ below~$U$.

(2) For $U\in\cl(\ca)^*$, the set $\{U_1,\dots,U_k\}=\max\cg_{\leq U}$
gives a~decomposition of $U$ because:
\begin{enumerate}
\item [a)] By the definition of $\cl(\ca)^*$ and the definition of
combinatorial building sets, we have $U=\Span(U_1,\dots,U_k)$, and, 
since $\sum_{i=1}^k\dim U_i=\dim U$, we have $U=\bigoplus_{i=1}^k U_i$;
\item [b)] for any $V\subseteq U$, $\bigoplus_{i=1}^k(U_i\cap V)
\subseteq V=\Span(U_1\wedge V,\dots,U_k\wedge V)\subseteq
\bigoplus_{i=1}^k(U_i\cap V)$, where "$\wedge$" denotes the meet 
operation in $\cl(\ca)^*$, hence $V=\bigoplus_{i=1}^k(U_i\cap V)$.
\qed
\end{enumerate}


\subsection{Intersection stratification of local arrangements 
               after blowup}

  Let a space~$X$ be given with an intersection stratification 
induced by a~local subspace arrangement, and let $G$ be a stratum in~$X$.
In the blowup of~$X$ at~$G$, $\bl_{G}X$, we find the following maximal strata: 
\begin{itemize}
\item[$\circ$] maximal strata in~$X$ that do not intersect with~$G$,
\item[$\circ$] blowups of maximal strata~$V$ at $G\cap V$, 
               $\bl_{G\cap V}V$,
               where $V$ is maximal in~$X$ and intersects~$G$,
\item[$\circ$] the exceptional divisor $\widetilde G$ replacing~$G$.
\end{itemize}
We consider the intersection stratification of $\bl_{G}X$ induced by
these maximal strata. We will later see (proof of 
Proposition~\ref{prop_compDP_cbl}) that in case $G$ is maximal
in a building set for the local arrangement in~$X$, then the union of
maximal strata in $\bl_GX$ is again a local arrangement with induced 
intersection stratification. In general, this is not the case,
see Remark~\ref{rem4.5}.

For ease of notation, let us agree here that formally blowing up an
empty (non-existing) stratum has no effect on the space.
We think about a stratum~$Y$ in~$X$, intersection of all maximal
strata $V_1,\ldots,V_t$ that contain~$Y$, as being replaced by the
intersection of corresponding maximal strata in~$\bl_{G}X$:
\begin{equation} \label{eq_blownupstratum}
\bl_{G\cap V_1} V_1 \,\,  \cap \, \,
             \ldots \,\,  \cap \, \,          
\bl_{G\cap V_t} V_t\, ,
\end{equation}
(recall that $\bl_{G\cap V_j} V_j\,{=}\,V_j$ for
$G\,{\cap}\,V_j\,{=}\,\emptyset$). The
intersection~(\ref{eq_blownupstratum}) being empty means that the 
stratum~$Y$ vanishes under blowup of~$G$. For notational convenience,
we most often retain names of strata under blowups, thereby
referring to the replacement of strata described above. 

\begin{rem} \label{rem4.5}
Let us mention here that blowing up a stratum in a local 
subspace arrangement does not necessarily result in a local subspace 
arrangement again. Consider the following arrangement of 2~planes and 
1~line in~$\dc^3$:
\[
A_1:\, \, y-z=0\, ,\quad
A_2:\, \, y+z=0\, ,\quad
L:\, \, x=y=0\, . 
\]
After blowing up~$L$, the planes $A_1$ and $A_2$ are replaced by 
complex line bundles over~$\dc {\rm P}^1$, which have in common their 
zero section~$Z$ and a complex line~$Y$; $L$ is replaced by a direct 
product of~$\dc$ and $\dc {\rm P}^1$, which intersects both line
bundles in~$Z$. The new maximal strata fail to form a local subspace 
arrangement in the point $Z\cap Y$.
\end{rem}


\subsection{Tracing incidence structure during arrangement model 
construction} \label{ssect_tracg_inc}

  We now give a more detailed description of the model construction
by De~Concini  \& Procesi via successive blowups, and then proceed with 
linking our notion of combinatorial blowups to the context of arrangement 
models.

Let $\ca$ be a complex subspace arrangement,
$\cg\,{\subseteq}\,\cl(\ca)$ a~geometric building set for~$\ca$, and 
$\{G_1,\ldots,G_t\}$ some linear extension of the partial containment 
order on associated strata in~$\dc^d$ such that $G_k\,{\supset}\,G_l$ 
implies $l<k$. The De~Concini-Procesi model $Y=Y_{\cg}$ of $\cm(\ca)$ 
is the result of blowing up the strata indexed by elements of~$\cg$ 
in the given order. Note that the linear order was chosen so that at 
each step the stratum which is to be blown up  does {\em not\/}
contain any other stratum indexed by an element of~$\cg$. At each step 
we consider intersection stratifications as described above, and we 
denote the poset of strata after blowup of $G_i$
with~$\cl^{\cg}_i(\ca)$. For the case of a~stratum $G_i$ being empty 
after previous blowups remember our agreement of considering blowups 
of $\emptyset$ as having no effect on a space. The later 
Proposition~\ref{prop_compDP_cbl} however shows that strata indexed 
by elements in~$\cg$ do not disappear during the sequence of blowups.

Let us remark that the combinatorial data of the initial
stratification, i.e., of the arrangement, prescribes much 
of the geometry of~$Y_{\cg}$: the complement $Y_{\cg}\,{\setminus}\,\cm(\ca)$
is a union of smooth irreducible divisors indexed by elements 
of~$\cg$, and these divisors intersect if and only if the set of
indices is nested in~$\cg$~\cite[Thm 3.2]{DP95}.


\begin{prop}\label{prop_compDP_cbl}
  Let $\ca$ be an arrangement of complex subspaces, $\cg$ a building
  set for $\ca$ in the sense of De~Concini~ \& Procesi, and
  $\{G_1,\ldots,G_t\}$ some linear extension of the partial
  containment order on associated strata as described above. Let
  $\bl_i^{\cg}(\ca)$ denote the geometric result of successively
  blowing up strata $G_1,\ldots,G_i$, for $1\,{\leq}\,i\,{\leq}\,t$.
  Then,
\begin{itemize}
\item[(1)] 
The poset of strata $\cl^{\cg}_i(\ca)$  of $\bl_i^{\cg}(\ca)$ can be 
described as the result of a~sequence of combinatorial blowups of the 
intersection lattice~$\cl\,{=}\,\cl(\ca)$:
\[
      \cl^{\cg}_i(\ca)\, \, = \, \, \bl_i(\cl)\, ,\qquad \mbox{for }\,
      1\,{\leq}\,i\,{\leq}\,t\, .
\]
{\rm (}Recall that $\bl_i(\cl)\,{=}\,\bl_{G_i}(\bl_{G_{i-1}}(\ldots \bl_{G_1}\cl))$
for $ 1\,{\leq}\,i\,{\leq}\,t$.{\rm )}
\item[(2)]
The union of maximal strata~$\ca_i^{\cg}$ in
$\bl_i^{\cg}(\ca)$ is a local subspace arrangement, with $\cg$ in 
$\cl^{\cg}_i(\ca)$ being a~building set for $\ca_i^{\cg}$ in the sense of \/ 
{\rm Definition~\ref{lsabsdf}}. {\rm (}Recall that $\cg$ here refers to 
the preimages 
of the original strata in $\cg\,{\subseteq}\,\cl(\ca)$ 
under the sequence of blowups.{\rm )}
\end{itemize}
\end{prop}

\noindent
\pr
We proceed by induction on the number of blowups. 
The induction start is obvious, since
the lattice of strata~$\cl_0^{\cg}(\ca)$ of the initial
stratification of $\dc^d$ coincides with the intersection lattice
$\cl(\ca)\,{=}\,\bl_0(\cl)$ of the arrangement~$\ca$. The union of maximal 
strata is the arrangement~$\ca$ itself with its given building set~$\cg$.

Assume that
$\cl_{i-1}^{\cg}(\ca)=\bl_{i-1}(\cl)$ for some $ 1\,{\leq}\,i\,{\leq}\,t$,
the union of maximal strata $\ca_{i-1}^{\cg}$ in $\bl_{i-1}^{\cg}(\ca)$ being
a local arrangement, and $\cg$ a building set for $\cl_{i-1}^{\cg}(\ca)$. 
Let $G\,{=}\,G_i$ be the next stratum to be blown up. 
First, we proceed in 4 steps to show that 
$\cl^{\cg}_i(\ca)\,{=}\,\bl_i(\cl)$. In 2 further steps we then verify
the claims in~(2).

\vskip4pt

\noindent
{\bf Step 1:} {\em Assign strata of $\bl_i^{\cg}(\ca)$ to elements 
in~$\bl_i(\cl)$.\/}

\vskip4pt

\noindent
We distinguish two types of elements in $\bl_i(\cl)$:
\[
\begin{array}{lrl}
\mbox{Type~I}:     & Y & 
                   \mbox{with }\, Y\,{\in}\,\bl_{i-1}(\cl) 
                   \mbox{ and}\, Y \not\geq G\, , \\
\mbox{Type~~I$\!$I}:\quad & [G,Y] & 
                   \mbox{with }\,Y\,{\in}\,\bl_{i-1}(\cl)\, , 
                   \, \, Y \not\geq G\, , \\
& &                \mbox{and }\, Y\,{\vee}\,G 
                   \mbox{ exists  in }\,\bl_{i-1}(\cl)\, . 
\end{array}
\]

\noindent
To $Y\,{\in}\,\bl_i(\cl)$ of type~I, assign $\bl_{G\cap Y}Y$ (recall
that blowing up an empty stratum does not change the space). Note that 
$\dim\bl_{G\cap Y}Y=\dim Y$.

\noindent
To $[G,Y]\,{\in}\,\bl_i(\cl)$ of type~I$\!$I, assign  
$(\bl_{G\cap Y}Y)\, \cap \, \widetilde G$, where $\widetilde G$ denotes
the exceptional divisor that replaces $G$ in $\bl_i^{\cg}(\ca)$. This 
description comprises $\widetilde G$ being assigned to $[G,\hat 0]$.
Note that $\dim(\bl_{G\cap Y}Y)\, \cap \, \widetilde G=\dim Y-1$.

\vskip4pt

\noindent
{\bf Step 2:} {\em Reverse inclusion order on the assigned spaces
    coincides with the partial order on~$\bl_i(\cl)$.\/}

\vskip4pt

\noindent
(1) $X,Y\,{\in}\, \bl_i(\cl)$, both of type~I: 
\[
X\, \leq_{\bl_i(\cl)} Y \, \Leftrightarrow \,
X\, \leq_{\bl_{i-1}(\cl)} Y \, \Leftrightarrow \,
X \supseteq_{\bl_{i-1}^{\cg}(\ca)} Y \, \Leftrightarrow \,
\bl_{G\cap X} X \supseteq \bl_{G\cap Y} Y ,
\]
where ``$\Leftarrow$'' in the last equivalence can be seen by first
noting that $Y \,{\setminus}\, (G\,{\cap}\,Y) \subseteq X 
\,{\setminus}\, (G\,{\cap}\,X)$, and then comparing points in the
exceptional divisors.

\vskip4pt

\noindent
(2) $X, [G,Y]\,{\in}\, \bl_i(\cl)$, $X$ of type~I,
$[G,Y]$ of type~I$\!$I: \newline
As above we conclude
\begin{eqnarray*}
X\, \leq_{\bl_i(\cl)} [G,Y] & \Leftrightarrow & 
X\, \leq_{\bl_{i-1}(\cl)} Y  \\
                            & \Leftrightarrow &
X \supseteq_{\bl_{i-1}^{\cg}(\ca)} Y \,\, \Rightarrow \,\,
\bl_{G\cap X} X \supseteq \bl_{G\cap Y} Y \cap \widetilde G\, .
\end{eqnarray*}
To prove the converse is rather subtle. Note first that 
$G\,{\cap}\,Y \subseteq G\,{\cap}\,X$. 
Assume that $G$ strictly contains~$G\,{\cap}\,X$, then 
both $G\,{\cap}\,X$ and $G\,{\cap}\,Y$ are not in the building set due 
to the linear order chosen on~$\cg$, and~$G$ is a factor of both
$G\,{\cap}\,X$ and $G\,{\cap}\,Y$. Let 
$F(G\,{\cap}\,X)\,{=}\,\{G,G_1,\ldots,G_k\}$, 
$F(G\,{\cap}\,Y)\,{=}\,\{G,H_1,\ldots,H_t\}$. $X$ written as a join
of elements in~$\bl_{i-1}(\cl)$ below the factors of~$G\,{\cap}\,X$
reads
\[
X \, = \, g_X \vee  Z_1  \vee \ldots  \vee Z_k
\] 
for some $g_X\,{\in}\,[\hat 0,G]$, $Z_i\,{\in}\,[\hat 0, G_i]$, for
$i=1,\ldots, k$. If $Z_i\,{<}\,G_i$, for some
$i\in\{1,\ldots,k\}$, we have
\begin{eqnarray*}
G\vee X 
& = & 
G \vee (g_X \vee Z_1 \vee  \ldots \vee Z_i  \vee \ldots  \vee Z_k) \\ 
& \leq & 
G \vee (g_X \vee G_1 \vee  \ldots \vee Z_i  \vee \ldots  \vee G_k) \\
& < & 
G \vee  G_1 \vee  \ldots \vee G_k \, = \,  G \vee X\, ,  
\end{eqnarray*}
by the ``necessity'' property of the 
Proposition~\ref{prop_buildingsets}(4), yielding a contradiction.
Hence,
\[
X \, = \, g_X \vee  G_1  \vee \ldots  \vee G_k\, ,
\]
and similarly, $Y\,{=}\, g_Y \vee H_1 \vee\ldots \vee H_t$
for some $g_Y\,{\in}\,[\hat 0,G]$.

For each $j\,{\in}\,\{1,\ldots,k\}$ there exists a unique 
$i_j\,{\in}\,\{1,\ldots,t\}$ such that $G_j\leq H_{i_j}$ by 
Proposition~\ref{prop_buildingsets_properties}(1). 
Thus, $\bigvee G_i\,{<}\, \bigvee H_j$, and, for showing that 
$X\leq Y$, it is enough to see that~$g_X\leq g_Y$.

We show that in an open neighborhood of any point
$y\,{\in}\,G\,{\cap}\,Y$, $g_Y\,{\subseteq}\,g_X$. This yields our
claim since strata in~$\bl_{i-1}^{\cg}(\ca)$ have pairwise transversal
intersections: if they coincide locally, they must coincide globally. 
With $\ca_{i-1}^{\cg}$ being a local arrangement, there exists
an open neighborhood of $y\,{\in}\,G\,{\cap}\,Y$  where
the stratification is biholomorphic to a stratification induced by a
subspace arrangement. We tacitly work in the arrangement setting,
using that $(\bl_{i-1}(\cl))_{\leq G\vee Y}$ is the intersection lattice 
of a product arrangement. The $\cg$-decomposition of $(G\vee Y)^{\perp}$
described in Definition~\ref{lsabsdf} yields (when transferred to the 
primal setting):
\[
g_Y\, = \,\Span(G,Y)\, .
\]
Analogously, $g_X\,{=}\,\Span(G,X)$.

In the linear setting we are concerned with, we interpret points in
the exceptional divisor of a blowup as follows:
\begin{equation}\label{eq_blowup}
\bl_{G\cap Y} Y \cap \widetilde G\,\, = \,\,
\{ (a,\Span(p,G\cap Y)) \, | \, 
a\in G\cap Y,\,  p\in Y\setminus (G\cap Y)  \}\, . 
\end{equation}
In terms of this description, the inclusion map 
$\bl_{G\cap Y} Y \cap \wti G \hookrightarrow \bl_G(\bl_{i-1}^{\cg}(\ca))$ reads
\[
  (a,\Span(p,G\cap Y)) \longmapsto  (a,\Span(p,G))\, .
\]
Therefore, $\bl_{G\cap Y} Y \cap \wti G$ being contained in
$\bl_{G\cap X} X\,{\subseteq}\, \bl_G(\bl_{i-1}^{\cg}(\ca))$ 
means that for $(a,\Span(p,G\cap Y))\in 
\bl_{G\cap Y} Y \,{\cap}\, \widetilde G$ there exists 
$q\,{\in}\, X\setminus (G\cap X)$ such that $\Span(p,G)\, = \,
\Span(q,G)$. In particular, $\Span(Y,G)\,{\subseteq}\,\Span(X,G)$, which
by our previous arguments implies that $Y\,{\subseteq}\,X$.

We assumed above that $G\,{\supset}\,G\,{\cap}\,X$. If $G\,{\cap}\,X$
coincides with~$G$, i.e., $X$ contains~$G$, then $g_X\,{=}\,X$ and a
similar reasoning applies to see that $Y\,{\subseteq}\,X$. Similarly
for $G\,{\cap}\,X\,{=}\,G\,{\cap}\,Y\,{=}\,G$.

\vskip4pt

\noindent
(3) $[G,X]$, $[G,Y]\,{\in}\, \bl_i(\cl)$, both of type~I$\!$I: 
\begin{eqnarray*}
[G,X]\leq_{\bl_i(\cl)} [G,Y] & \Leftrightarrow & 
X\, \leq_{\bl_{i-1}(\cl)} Y  \\
                            & \Leftrightarrow &
X \supseteq_{\bl_{i-1}^{\cg}(\ca)} Y \,\, \Leftrightarrow \,\,
\bl_{G\cap X} X \cap \widetilde G \supseteq \bl_{G\cap Y} Y \cap \widetilde G\, ,
\end{eqnarray*}
where ``$\Leftarrow$'' follows from~(2) and 
$\bl_{G\cap X} X\,{\supseteq}\,
\bl_{G\cap X} X \cap \widetilde G\,{\supseteq}\,
 \bl_{G\cap Y} Y \cap \widetilde G$.

\vskip4pt

\noindent
{\bf Step 3:} {\em Each of the assigned spaces is the intersection of 
maximal strata in $\bl_i^{\cl}(\ca)$.}

\vskip4pt

\noindent
It is enough to show that spaces assigned to elements of type~I
in~$\bl_i(\cl)$ are intersections of new maximal strata. Those
associated to elements of type~I$\!$I then are intersections as well
by definition.

Let $Y\,{\in}\,\bl_i(\cl)$, $Y\,{\not\geq}\,G$, 
and $Y\,{=}\,\cap_{i=1}^t\, V_i$ with $V_1,\ldots,V_t$ the maximal strata 
in~$\bl_{i-1}^{\cg}(\ca)$ containing~$Y$. We claim that
\begin{equation} \label{eq_intdiv}
\bl_{G\cap Y}Y\, \, = \, \, \bigcap_{i=1}^t\,\bl_{G\cap V_i}V_i\, . 
\end{equation}
For the inclusion~``$\subseteq$'' note that 
$\bl_{G\cap  Y}Y\,{\subseteq}\,\bl_{G\cap V_i}V_i$ is a direct
consequence of $Y\,{\subseteq}\, V_i$ as discussed in Step~2~(1).

For the reverse inclusion we need the following identity:
\begin{equation}\label{eq_join}
     \bigvee_{i=1}^t\, (G\wedge V_i)\,\, = \, \, G\wedge Y\, .
\end{equation}
This identity holds in any semilattice without referring to $G$ 
being an element of the building set.

Let $\alpha \,{\in}\, \cap_{i=1}^t\,\bl_{G\cap V_i}V_i$. In case
$\alpha \,{\in}\, \cap_{i=1}^t\, V_i\, {\setminus}\, (G\,{\cap}\,V_i)$,
we conclude that $\alpha\,{\in}\, Y\, {\setminus}\, (G\,{\cap}\,Y) $.
We thus assume that $\alpha$ is contained in the intersection of
exceptional divisors $\widetilde{G\,{\cap}\,V_i}$,  $i=1,\ldots,t$. 
We again switch to local considerations in
the neighborhood of a point~$y\,{\in}\,G\,{\cap}\,Y$, using that it 
carries a stratification biholomorphic to an arrangement
stratification.

Using the  description (\ref{eq_blowup}) of points in exceptional
divisors that are created by blowups in the arrangement setting, 
$\alpha\,{\in}\,\cap_{i=1}^t\,\widetilde{G\,{\cap}\,V_i}\,{\subseteq}\,
\cap_{i=1}^t\,\bl_{G\cap V_i}V_i$ means that there exist $a\,{\in}\,
\cap_{i=1}^t\, (G\,{\cap}\,V_i)$, and  
$p_i\,{\in}\, V_i\, {\setminus}\, (G\,{\cap}\,V_i)$ for $i=1,\ldots,t$,
with 
\[
\alpha \,\, = \, \, (a,\Span(p_i,G\,{\cap}\,V_i))\, \in \, 
\bl_{G\cap  V_i}V_i\, .
\] 
In particular, $\Span(p_i,G)\,{=}\,\Span(p_j,G)$ for
$1\,{\leq}\,i,j\,{\leq}\, t$. Thus,
\[
\Span(p_j,G) \, \subseteq \bigcap_{i=1}^t\, \Span(V_i,G)\, \, = \, \,
\Span(Y,G)
\] 
using the identity~(\ref{eq_join}). We conclude that there exists
$y\,{\in}\, Y\,{\setminus}\,(G\,{\cap}\,Y)$ such that $\Span(y,G)\,
{=}\, \Span(p_j,G)$ for all $j\,{\in}\,\{1,\ldots,k\}$, hence
\[
\alpha\, = \, (a,\Span(y,G\,{\cap}\,Y))\,\,\in\,\, \bl_{G\cap Y}Y\, .
\]

Though we are for the moment not concerned with the case of
$Y\,{\subseteq}\,G$, we note for later reference
that~(\ref{eq_intdiv})
remains true, with $\bl_Y Y\,{=}\,\emptyset$ meaning that the
intersection on the right-hand side is empty.  
Following the proof of the
inclusion~``$\supseteq$'' in~(\ref{eq_intdiv}) for
$G\,{\cap}\,Y\,{=}\,Y$,
we first find that the intersection of blowups can only contain points 
in the exceptional divisors. Assuming~$\alpha\,{\in}\,\cap_{i=1}^t\,
\widetilde{G\,{\cap}\,V_i}$ we arrive to a contradiction when
concluding that $\Span(p_j,G)\,{\subseteq}\,
\cap_{i=1}^t\,\Span(V_i,G)\,{=}$ $\Span(Y,G)\,{=}\,G$ for
$j=1,\ldots,t$.

\vskip4pt

\noindent
{\bf Step 4:} {\em Any intersection of maximal strata in 
$\bl_i^{\cg}(\ca)$ occurs as an assigned space.}

\vskip4pt

\noindent
Every intersection involving the exceptional divisor~$\widetilde G$
occurs if we can show that all other intersections occur
(intersections that additionally involve~$\widetilde G$ then are  
assigned to corresponding elements of type~I$\!$I). 

Consider
$W\, = \, \bigcap_{i=1}^t\,\bl_{G\cap V_i}V_i$,
where the $V_i$ are maximal strata in $\bl_{i-1}^{\cg}(\ca)$; recall
here that a blowup in an empty stratum does not alter the space. 
We can assume that $\cap_{i=1}^t\, V_i\,{\not =}\, \emptyset$,
otherwise the intersection~$W$ were empty. With the
identity~(\ref{eq_intdiv}) in Step~3 we conclude that either $W\,{=}\,
\emptyset$ (in case $\cap_{i=1}^t\, V_i\,{\subseteq}\,G$) or
$W\,{=}\,\bl_{G\cap \bigcap_{i=1}^t V_i} \cap_{i=1}^t\, V_i$,
in which case it is assigned to the element $\cap_{i=1}^t\, V_i$
in $\bl_i(\cl)$.

\vskip4pt

\noindent
{\bf Step 5:} {\em $\ca_{i}^{\cg}$ is a local subspace arrangement
in $\bl_i^{\cg}(\ca)$.} 

\vskip4pt

\noindent
It follows from the description~(\ref{eq_intdiv}) of strata in  
$\bl_i^{\cg}(\ca)$ that all intersections of maximal strata are
connected and smooth. It remains to show that $\ca_{i}^{\cg}$
locally looks like a subspace arrangement. Let $y\,{\in}\,\ca_{i}^{\cg}$.
We can assume that $y$ lies in the exceptional divisor~$\wti G$. 
Let $x \in G\,{\subseteq}\,\ca_{i-1}^{\cg}$ be the image of~$y$ under 
the blowdown map.

We first give a local description around $x$ in $\ca_{i-1}^{\cg}$.
By induction hypothesis, there exists a neighborhood~$N$ of~$x$, and 
an arrangement of linear subspaces~$\cb$ in~$\dc^n$
such that the pair $(N,\ca_{i-1}^{\cg}\,{\cap}\,N)$ is biholomorphic
to the pair~$(\dc^n,\cb)$. We can assume that under this biholomorphic 
map, $x$ is mapped to the origin. Let $T\,{=}\,\bigcap_{B\in \cb}\,B$ 
and note that $G\,{\cap}\,N$ is mapped to some subspace~$\Ga$ in $\cb$.

With $G$ being maximal in the building set for~$\ca_{i-1}^{\cg}$,
$\cb/T$ is a product arrangement with one of the factors being an
arrangement in~$\Ga/T$. More precisely, there exists a subspace 
$\Ga'\,{\subseteq}\,\dc^n$, and two subspace arrangements, $\cc$ in $\Ga/T$
and $\cc'$ in $\Ga'/T$, such that
\begin{itemize}
\item[(1)] $\Ga/T\, \oplus\, \Ga'/T\, \oplus\, T\, = \, \dc^n$,  
\item[(2)] $\cb\, =\ 
\{A\, \oplus\, \Ga'/T\, \oplus\, T\,| \, A\,{\in}\,\cc\} \, \cup \,
\{\Ga/T\, \oplus\, A'\, \oplus\, T\,| \, A'\,{\in}\,\cc'\}$. 
\end{itemize}

Blowing up $G$ in $\bl_{i-1}^{\cg}(\ca)$ locally corresponds to
blowing up $\Gamma$ in $\dc^n$. Let $t$ be the point on the special divisor
$\wti \Ga$ corresponding to~$y \in \wti G$, thus $t$ maps to the origin 
in~$\dc^n$ under the blowdown map. A neighborhood of~$t$ in $\bl_{\Ga}\dc^n$ 
is an $n$-dimensional open ball which can be parameterized as a direct sum
\[
M \,\oplus \, M' \, \oplus \, I \, \oplus \,T \, .
\]
Here, $M$ is an open ball around~$0$ in $\Ga/T$, $M'$ is an open ball on the
unit sphere in  $\Ga'/T$ around the point of intersection with the 
line~$\langle p \rangle$
in  $\Ga'/T$ that defines $t$ as a point in the exceptional divisor, 
$t=(0, \Span(p,\Ga)) \in \wti \Ga$ (compare (\ref{eq_blowup})), 
and $I$ an open unit ball in~$\dc$.

The maximal strata in this neighborhood are the following:
\begin{itemize}
\item[$\circ$]
      the hyperplane $M \,\oplus \, M' \, \oplus \, \{0\} \, \oplus \,T$,
      as the exceptional divisor,
\item[$\circ$]
      $(M \cap A)\,\oplus \, M' \, \oplus \, I \, \oplus \,T$, replacing
      $A\, \oplus\, \Ga'/T\, \oplus\, T$ after blowup,
\item[$\circ$]
      $M \,\oplus \, (M'\cap A') \, \oplus \, I \, \oplus \,T$, replacing
      $\Ga/T\, \oplus\, A'\, \oplus\, T$ after blowup for $A'\neq 0$.
\end{itemize} 

This proves that around $t$ in $\bl_{\Ga}\dc^n$ we have the structure of a 
local subspace arrangement, which in turn shows the local arrangement 
property around $y$ in $\ca_i^{\cg}$.

\vskip4pt

\noindent
{\bf Step 6:} {\em $\cg$ is a building set for $\ca_i^{\cg}$
 in the sense of \/{\rm Definition~\ref{lsabsdf}}.}

\vskip4pt

\noindent
$\cg$ is a~combinatorial building set by
Proposition~\ref{prop_singleblow}. Complementing this with the
dimension information about the strata, we conclude, by
Proposition~\ref{combgeomprop}(2), that $\cg$ is a~geometric building
set. 
\qed


\section{Simplicial resolutions of toric varieties}
\label{ssect_tv}

The study of toric varieties has proved to be a field of fruitful
interplay between algebraic and convex geometry: toric
varieties are determined by rational polyhedral fans, and many of their 
algebraic geometric properties are reflected by combinatorial
properties of their defining fans. 

We recall one such correspondence -- between subdivisions of fans and
special toric morphisms -- and show that so-called stellar subdivisions 
are instances of combinatorial blowups. This allows us to apply our
Main Theorem in the present context: Given a polyhedral fan, we
specify a class of {\em simplicial\/} subdivisions, and, interpreting
our notions of building sets and nested sets, we describe the 
incidence combinatorics of the subdivisions in terms of the
combinatorics of the initial fan. For background material on
toric varieties we refer to the standard sources~\cite{Da78,Od88,Ful93,Ew96}.

Let $X_{\Sigma}$ be a toric variety defined by a rational polyhedral
fan~$\Sigma$. Any subdivision of~$\Sigma$ gives rise to a proper,
birational toric morphism between the associated toric 
varieties~(cf~\cite[5.5.1]{Da78}). In particular, simplicial
subdivisions yield toric morphisms from quasi-smooth toric varieties
to the initial variety -- so-called {\em simplicial resolutions}.
Quasi-smooth toric varieties being rational homology manifolds, such
morphisms can replace smooth resolutions for (co)homological
considerations.

We define a particular, elementary, type of subdivisions:

\begin{df} \label{def_stellsd}
Let $\Sigma\,{=}\,\{\sigma\}_{\sigma \in \Sigma}\,{\subseteq}\,\dr^d$
be a polyhedral fan, i.e., a collection of closed polyhedral
cones~$\sigma$ in~$\dr^d$ such that $\sigma\,{\cap}\,\tau$ is a cone
in~$\Sigma$ for any $\sigma, \tau\,{\in}\,\Sigma$. Let ${\rm cone}(x)$ be
a ray in~$\dr^d$ generated by $x\,{\in}\,{\rm relint}\, \tau$ for some 
$\tau\,{\in}\,\Sigma$. The {\em stellar subdivision\/} ${\rm sd}(\Sigma,x)$
of $\Sigma$ in $x$ is given by the collection of cones
\[
(\,\Sigma\,\setminus\, {\rm star}(\tau, \Sigma)\, )
\,\, \cup \, \, 
\{\, {\rm cone}(x,\rho)\, | \, \rho \subseteq \sigma \,\,
                               \mbox{ for some } \,\,
                               \sigma \in {\rm star}(\tau, \Sigma)\,    
\}\, ,
\] 
where ${\rm star}(\tau, \Sigma)\,{=}\,\{\sigma\,{\in}\,\Sigma\, | \, 
\tau\,{\subseteq}\,\sigma\}$, and ${\rm cone}(x,\rho)$ the closed polyhedral 
cone
spanned by~$\rho$ and~$x$. If only concerned with the combinatorics of 
the subdivided fan, we also talk about stellar subdivision of $\Sigma$ 
in $\tau$, ${\rm sd}(\Sigma,\tau)$, meaning any stellar subdivision in~$x$
for~$x\,{\in}\,{\rm relint}\, \tau$. 
\end{df} 

\begin{prop}
Let $\cf(\Sigma)$ be the face poset of a polyhedral fan~$\Sigma$,
i.e., the set of closed cones  in~$\Sigma$ ordered by inclusion, together with the
zero cone~$\{0\}$ as a minimal element. For $\tau\,{\in}\,\Sigma$,
the face poset of the stellar subdivision of $\Sigma$ in $\tau$ can be
described as the combinatorial blowup of $\cf(\Sigma)$ at~$\tau$:
\[
  \cf({\rm sd}(\Sigma, \tau))\,\, = \,\, \bl_{\tau}\cf(\Sigma)\, . 
\] 
\end{prop}

\pr We observe that removing ${\rm star}(\tau, \Sigma)$ from $\Sigma$
corresponds to removing $\cf(\Sigma)_{\geq \tau}$ from $\cf(\Sigma)$,
while adding cones, as described in Definition~\ref{def_stellsd},
corresponds to extending $\cf(\Sigma)\,{\setminus}\,\cf(\Sigma)_{\geq
  \tau}$ by elements $[\tau,\rho]$ for $\rho\,{\in}\,\cf(\Sigma)$,
$\rho\,{\subseteq}\,\sigma$ for some $\sigma \,{\in}\,{\rm star}(\tau,
\Sigma)$. The comparison of order relations is straightforward. \qed

\vspace{0.2cm}
We apply our Main Theorem to the present context.

\begin{thm}
Let~$\Sigma$ be a polyhedral fan in~$\dr^d$ with 
face poset~$\cf(\Sigma)$. Let $\cg\,{\subseteq}\,\cf(\Sigma)$ be a building 
set 
of~$\cf(\Sigma)$ in the sense of Definition~\ref{df_buildg}, $\cn(\cg)$ the
complex of nested sets in~$\cg$ (cf.\ Definition~\ref{df_nested}).
Then, the consecutive application of stellar subdivisions in every cone 
$G\,{\in}\,\cg$ in a non-increasing order yields a simplicial subdivision
of~$\Sigma$ with face poset equal to the face poset of~$\cn(\cg)$.
\end{thm}

\noindent
As examples of building sets for face lattices of 
polyhedral fans let us mention:
\begin{itemize}
\item[(1)] the full set of faces, with the corresponding complex of nested sets 
           being the order complex of~$\cf(\Sigma)$ (stellar subdivision in all 
           cones results in the barycentric subdivision of the fan);
\item[(2)] the set of rays together with the non-simplicial faces of~$\Sigma$;
\item[(3)] the set of irreducible elements in $\cf(\Sigma)$: the set of rays 
           together with all faces of~$\Sigma$ that are not products 
           of some of their proper faces.
\end{itemize}

\begin{rem} \label {rm_smooth_tv}
  For a smooth toric variety~$X_{\Sigma}$, the union of closed
  codimension~1 torus orbits is a local subspace arrangement, in
  particular, the codimension~1 orbits form a divisor with normal
  crossings, \cite[p.~100]{Ful93}.  The intersection stratification of
  this local arrangement coincides with the torus orbit stratification
  of the toric variety. For any face $\tau$ in the defining
  fan~$\Sigma$, the torus orbit $\co_{\tau}$ together with all orbits
  corresponding to rays in~$\Sigma$ form a~geometric building set. Our
  proof in the subsection~\ref{ssect_tracg_inc} applies in this
  context with $\co_{\tau}$ playing the role of~$G$. We conclude that
  under blowup of~$X_{\Sigma}$ in the closed torus orbit $\co_{\tau}$,
  the incidence combinatorics of torus orbits changes exactly in the
  way described by a stellar subdivision of~$\Sigma$ in~$\tau$. This
  is the combinatorial part of the well-known fact that in the smooth
  case, the blowup of $X_{\Sigma}$ in a torus orbit $\co_{\tau}$
  corresponds to a regular stellar subdivision of the fan~$\Sigma$ in
  $\tau$~\cite{MO73}.
 \end{rem}


\clearemptydoublepage
\chapter*{PART II \\[2cm] Complexes of Trees \\[0.1cm] 
and Quotient Constructions}
\addcontentsline{toc}{chapter}{PART II. Complexes of Trees and 
Quotient Constructions}
\vspace{9cm}
\mbox{ }\hfill
\begin{minipage}{7.6cm}
  { In the reproof of chance}\\
{ Lies the true proof of men}\\[0.2cm]
    -William Shakespeare, {\it Troilus and Cressida}
\end{minipage}

\clearemptydoublepage
\chapter[Spaces of Complex Monic Polynomials and Complexes of Forests]
{Rational Homology of Spaces of Complex Monic Polynomials with
  Multiple Roots and Complexes of Marked Forests }

\section{The stratification by root multiplicities of the space of 
complex monic polynomials}
              
Let $n$ be an~integer, $n\geq 2$. We view $n$-dimensional complex
space $\dc^n$ as the space of all monic polynomials with complex
coefficients of degree $n$ by identifying
$a=(a_0,\dots,a_{n-1})\in\dc^n$ with
$f_a(z)=z^n+a_{n-1}z^{n-1}+\dots+a_0$.  To each
$\lam=(\lam_1,\dots,\lam_t)\vdash n$ one can associate a~topological
space as follows (we refer the reader to the Appendix~A for
a~description of our conventions on the terminology of number and set
partitions).
\begin{df} \label{df1.1}
  $\wsl$ is the set of all $a\in\dc^n$, for which the roots of
  $f_a(z)$ can be partitioned into sets of sizes
  $\lam_1,\dots,\lam_t$, so that within each set the roots are equal.
  Clearly, $\wsl$ is a~closed subset of $\,\dc^n$. Let $\sla$ be the
  one-point compactification of~$\wsl$.
\end{df}

It is easy to see that the Definition~\ref{df1.1} is equivalent to the
description of $\sla$ given in the Chapter~1, so our notation is
consistent in that point.
In this chapter we shall focus on the (reduced) rational Betti numbers
of the spaces $\sla$. In~\cite{Ar70a}, V.I.\ Arnol'd has computed
$\tb_*(\sla,\dq)$ for $\lam=(k^m,1^{n-km})$.
\begin{thm}\label{arn} {\rm (Arnol'd, \cite{Ar70a}).}
  Let $\lam=(k^m,1^{n-km})$ for some natural numbers $k\geq 2$, $m$,
  and $n\geq km$.  Then
\begin{equation}\label{are}
 \tb_i(\sla,\dq)=\begin{cases} 1,& \text{for } i=2l(\lam);\\
                               0,& \text{otherwise.}
  \end{cases}
\end{equation}
\end{thm}

  In \cite{SuW97} Sundaram and Welker conjectured that 
   
\begin{conj}\label{c3}
  For any number partition $\lam$, $\tb_i(\sla,\dq)=0$ unless $i=2l(\lam)$.
\end{conj}

In this chapter we shall give a new, combinatorial proof of the
Theorem~\ref{arn} and disprove Conjecture~\ref{c3}. To do that, we
shall introduce a family of topological spaces $\xlm$, indexed by
pairs of number partitions $(\lam,\mu)$, satisfying $\lam\vdash\mu$.
$\xlm$ will be defined so that the following equality is satisfied

\begin{equation}\label{rswf}
  \tb_i(\sla,\dq)=\sum_{\lam\vdash\mu\vdash n}\tb_{i-2l(\mu)-1}(\xlm,\dq).
\end{equation}
  
\nin
  Here is the summary of the chapter.

  
\vskip4pt 
\nin 
{\bf Section~\ref{s2}.} We define the topological spaces $\xlm$ and 
derive \eqref{rswf}.

\vskip4pt \nin {\bf Section~\ref{s3}.} We give a~combinatorial
description of the cell structure of the triangulated spaces~$\xlm$ in
terms of marked forests, see Theorem~\ref{main1}. This description is
the backbone of the chapter, it serves as both language and intuition
for the material in the subsequent sections. One consequence of
Theorem~\ref{main1} is that homology groups of $\xlm$ may be computed
from a~chain complex, whose components are freely generated by marked
forests, and the boundary operator is described in terms of
a~combinatorial operation on such forests (deletion of level sets).

\vskip4pt 
\noindent 
{\bf Section \ref{s4}.} We prove a~general theorem about
collapsibility of certain triangulated spaces. The direct
combinatorial argument is heavily relying on the combinatorial cell
description, from Section~\ref{s3}, of~$\xlm$. We would like to
mention that for the special cases $\lam=(k,1^{n-k})$ and $\lam=(k^m)$
the new proof of the Theorem~\ref{arn} was also found in~\cite{SuW97},
the argument there also made use of the Theorem~\ref{swt} (as
Proposition~\ref{tgen} shows, the case $\lam=(k^m)$ is especially
simple). However, the Theorem~\ref{main2} is the first combinatorial
(modulo Theorem~\ref{swt}) proof of the result of Arnol'd in the
general case.

\vskip4pt \nin {\bf Section~\ref{s5}.} We disprove the conjecture of
Sundaram and Wel\-ker. Besides giving a~counterexample we prove this
conjecture for a~class of number partitions, which we call generic
partitions.

 \section{Orbit arrangements and spaces $\xlm$} \label{s2}
  
\subsection{Reformulation in the language of orbit arrangements}

\noindent
Following \cite{SuW97}, we give a~different interpretation of
the numbers $\tb_i(\sla,\dq)$, for general $\lam$. First, let us
observe that the symmetric group $\cs_n$ acts on $\dc^n$ by permuting
the coordinates, so we can consider the space $\dc^n/\cs_n$ endowed
with the quotient topology. It is a~classical fact that the map
$\phi:\dc^n\ra\dc^n/\cs_n$, mapping a~polynomial to the (unordered)
set of its roots, is a~homeomorphism, which extends to the one-point
compactifications. Therefore $\dc^n\cong\dc^n/\cs_n$ and
$\sla\cong\phi(\sla)=\phi(\wsl)\cup\{\infty\}$.

   $\phi(\wsl)$ can be viewed as the configuration space of $n$ unmarked points 
on~$\dc$ such that the number partition given by the coincidences among the points 
is refined by~$\lam$. For example, $\phi(\wti\Sigma_{(2,1^{n-2})})$ is the 
configuration space of $n$ unmarked points on $\dc$ such that at least 2 points 
coincide. Using this point of view, $\phi(\sla)$ can be described in the language 
of orbit subspace arrangements.

\begin{df}
  For $\pi\vdash[n]$, $\pi=(S_1,\dots,S_t)$,
  $S_j=\{i_1^j,\dots,i_{|S_j|}^j\}$, $1\leq j\leq t$, $K_\pi$ is
  the~sub\-space given by the equations
  $x_{i_1^1}=\dots=x_{i_{|S_1|}^1},\dots,x_{i_1^t}=\dots=x_{i_{|S_t|}^t}$.
  For $\lam\vdash n$, set $I_\lam=\{\pi\vdash[n]\,|$ $\text{\rm
    type}\,(\pi)=\lam\}$ and define $\ca_\lam=\{K_\pi\,|\,\pi\in
  I_\lam\}$. $\ca_\lam$'s are called {\bf orbit arrangements}.
\end{df}
The orbit arrangements were introduced in \cite{Bj94} and studied in
further detail in~\cite{Ko97}. They provide the appropriate language
to describe $\phi(\sla)$, indeed
\begin{equation}\label{fisla}
  \phi(\sla)=\Gamma_{\ca_\lam}^{\cs_n}.
\end{equation} 

An~important special case is that of the braid arrangement
$\ca_{n-1}=\ca_{(2,1^{n-2})}$, which corresponds under $\phi$ to
${\wti\Sigma}_{(2,1^{n-2})}$, the~space of all monic complex
polynomials of degree~$n$ with at least one multiple root. The name
``braid arrangement'' stems from the fact that $\dc^n\sm
V_{\ca_{n-1}}$ is a~classifying space of the colored braid group,
see~\cite{Ar69}. The intersection lattice $\cl_{\ca_{n-1}}$ is usually
denoted $\Pi_n$. It is the~poset consisting of all set partitions
of~$[n]$, where the partial order relation is refinement. Furthermore,
for $\lam\vdash n$, the intersection lattice of $\ca_\lam$ is denoted
$\Pi_\lam$.  It is the~subposet of $\Pi_n$ consisting of all elements
which are joins of elements of type~$\lam$, with the minimal element
$\hat 0$ attached.
          
\subsection{Applying Sundaram-Welker formula}
  
\noindent   
 The following 
formula of S.~Sun\-da\-ram and V.~Welker, \cite{SuW97}, is vital for our approach. 

\begin{thm}\label{swt}{\rm (\cite[Theorem 2.4(ii) and Lemma 2.7(ii)]{SuW97}).}\newline 
  Let $\ca$ be an~arbitrary subspace arrangement in $\dc^n$ with
  an~action of a~finite group $G\subset\U_n(\dc)$. Let $\cd_\ca$ be
  the intersection of $V_\ca$ with the $(2n-1)$-sphere (often called
  the link of $\ca$). Then there is the following isomorphism of
  $G$-modules.
\begin{equation}\label{swf}
  \wti H_i(\cd_\ca,\dq)\cong_G\bigoplus_{x\in\cl_\ca^{>\hat 0}/G} 
  \text{\rm Ind}_{\st(x)}^G (\wti H_{i-\dim x}(\da(\cl_\ca(\hat 0,x)),\dq)),
\end{equation}
where the sum is taken over representatives of the orbits of $G$ in
$\cl_\ca\sm\{\hat 0\}$, under the action of $G$, one representative
for each orbit.
\end{thm} 
   
Clearly $\Gamma_\ca^G\cong\susp\,(\cd_\ca/G)$. Recall that if a~finite
group~$G$ acts on a~finite cell complex $K$ then $\tb_i(K/G,\dq)$ is
equal to the multiplicity of the trivial representation in the induced
representation of $G$ on the $\dq$-vector space $\wti H_i(K,\dq)$, see
for example \cite[Theorem 1]{Co56}, \cite{Br72}. Hence, it follows
from \eqref{swf}, and the Frobenius reciprocity law, that
\begin{equation}\label{star}
  \tb_i(\Gamma_\ca^G,\dq)=\sum_{x\in\cl_\ca^{>\hat 0}/G} 
  \tb_{i-\dim x-1}(\da(\cl_\ca(\hat 0,x))/{\st(x)},\dq),
\end{equation}
where $\st(x)$ denotes the stabilizer of $x$.
   
\subsection{Spaces $\xlm$ and their properties}

\noindent
Let us now restate this identity in the special case of orbit
arrangements.  As mentioned above, the intersection lattice of
$\ca_\lam$ is $\Pi_\lam$. It has an~action of the symmetric group
$\cs_n$, which, for any $\pi\in\Pi_\lam$ induces an~action of
$\st(\pi)$ on $\da(\Pi_\lam(\hat 0,\pi))$.

\vskip4pt

{\bf Notation. }{\it Let $\xlm$ denote the topological space
  $\da(\pl(\hat 0,\pi))/\st(\pi)$, where the set partition~$\pi$ has
  type $\mu$.  If there is no set partition $\pi\in\pl$ of type $\mu$,
  i.e., if $\mu$ cannot be obtained as a join of $\lambda$'s, then let
  $\xlm$ be a~point.}

\vskip4pt

For fixed $\mu$, the space $\xlm$ does not depend on the choice of
$\pi$. Observe that $\xlm$ is in general not a~simplicial complex,
however it is a~triangulated space, (a~regular CW complex with each
cell being a~simplex, see~\cite[Chapter~I, Section~1]{GeM96}), with
its cell structure inherited from the simplicial complex
$\da(\Pi_\lam(\hat 0,\pi))$. In general, whenever $G$ is a~finite
group which acts on a~poset $P$ in an~order-preserving way, $\da(P)/G$
is a~triangulated space whose cells are orbits of simplices of
$\da(P)$ under the action of $G$; this is obviously not true in
general for an~action of a~finite group on a~finite simplicial
complex.
        
  Clearly, \eqref{fisla} together with \eqref{star}, and the fact that
$\phi$ is a~homeomorphism, implies \eqref{rswf}.
Let us quickly analyze \eqref{rswf}. $X_{\lam,\lam}=\emptyset$ makes
a~contribution~1 in dimension $2l(\lam)$. Assume $\mu\neq\lam$, then $1\leq l(\mu)
\leq l(\lam)-1$ and $\xlm\neq\emptyset$. $\dim\xlm=l(\lam)-l(\mu)-1$, hence 
$\tb_{i-2l(\mu)-1}(\xlm,\dq)=0$ unless $0\leq i-2l(\mu)-1\leq l(\lam)-l(\mu)-1$,
that is $2l(\mu)+1\leq i\leq l(\lam)+l(\mu)$. It follows from \eqref{rswf} that
$$\tb_{2l(\lam)}(\sla,\dq)=1,\text{  and }\tb_i(\sla,\dq)=0\text{ unless }
3\leq i\leq 2l(\lam).$$

 The purpose of this chapter is to investigate the values $\tb_i(\sla,\dq)$ for
$3\leq i\leq 2l(\lam)-1$, by studying $\tb_i(\xlm,\dq)$. We shall prove that 
the latter are equal to~0 for a~certain set of pairs $(\lam,\mu)$, 
$\lam\vdash\mu$, of partitions, including the case in Theorem~\ref{arn}, 
($\lam=(k^m,1^{n-km})$, $\mu$ is arbitrary such that $\lam\vdash\mu$),  
and we shall give an~example that this is not the case in general.

              \section{The cell structure of $\xlm$ and marked forests}
          \label{s3}
                 
 \subsection{The terminology of marked forests}

\noindent  
  In order to index the simplices of $\xlm$ we need to introduce some terminology 
for certain types of trees with additional data. For an~arbitrary forest of rooted
trees $T$ (we only consider finite graphs), let $V(T)$ denote the set of the 
vertices of~$T$, $R(T)\subseteq V(T)$ denote the set of the roots of~$T$ and 
$L(T)\subseteq V(T)$ denote the set of the leaves of~$T$. For any integer $i\geq 0$,
let $l_i(T)$ be the number of $v\in V(T)$ such that, $v$ has distance $i$ to 
the root in its connected component. 

\begin{df}
  A forest of rooted trees $T$ is called a~{\bf graded forest of rank $r$} if 
$l_{r+2}(T)=0$, $l_{r+1}(T)=|L(T)|$, and the sequence $l_0(T),\dots,l_{r+1}(T)$ 
is strictly increasing.  
\end{df}

For $v,w\in V(T)$, $w$ is called {\it a~child} of $v$ if there is an~edge between $w$ 
and $v$ and the unique path from $w$ to the corresponding root passes through $v$.
For $v\in V(T)$, we call the distance from $v$ to the closest leaf the {\it height}
of $v$. For example, in a~graded forest of rank $r$, leaves have height 0 and roots
have height $r+1$.

\begin{df}\label{df2.2}
  A {\bf marked forest of rank $r$} is a~pair $(T,\eta)$, where $T$ is a~graded
forest of rank $r$ and $\eta$ is a~function from $V(T)$ to the set of natural 
numbers such that for any vertex $v\in V(T)\sm L(T)$ we have 
\begin{equation}\label{mc1}
\eta(v)=\sum_{w\in\text{\rm{children}}(v)}\eta(w).
\end{equation}
\end{df}

We remark that the set of the marked forests of rank $r$, such that
not all leaves have label~1, is equal to the set of graded forests of
rank $r+1$.  Indeed, instead of labeling the vertices with natural
numbers so that~\eqref{mc1} is satisfied, one can as well attach a~new
level of leaves so that each ``old leaf'' $v$ has $\eta(v)$ children.
Then the old labels will correspond to the numbers of the new leaves
below each vertex. For our context it is more convenient to use labels
rather than auxiliary leaves, i.e.,~it is more handy to label all
vertices rather than just the leaves, so we stick to the terminology
of Definition~\ref{df2.2}.

   For a~marked forest $(T,\eta)$ of rank $r$ and $0\leq i\leq r+1$, we have 
a~number partition $\lam_i(T,\eta)=\{\eta(v)\,|\,v\text{ has height }i\}$. Clearly 
$\lam_0(T,\eta)\vdash\dots\vdash\lam_r(T,\eta)\vdash\lam_{r+1}(T,\eta)$.

\begin{df}\label{lmd}
  Let $\lam\vdash\mu\vdash n$, $\lam\neq\mu$. A {\bf $\lmu$-forest of
    rank $r$} is a~marked forest of rank $r$, $(T,\eta)$ such that
  $\mu=\lam_{r+1}(T,\eta)$ and $\lam\vdash\lam_0(T,\eta)$.
\end{df}

To simplify the language, we call $((2,1^{n-2}),\mu)$-forests simply
{\it $\mu$-forests}, and $((2,1^{n-2}),(n))$-forests simply {\it
  $n$-trees}.
 
\vskip-10pt
\hskip-100pt $$\epsffile{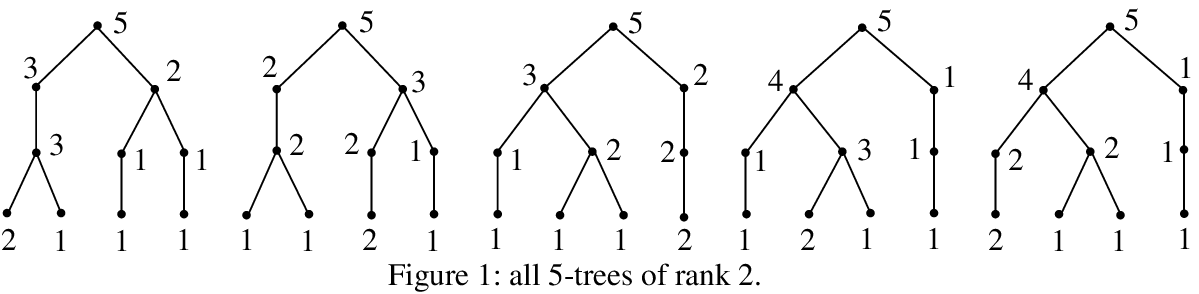}$$
 \begin{center}\vskip-2pt
Figure 3.1. All 5-trees of rank 2.
\end{center}

Whenever $(T,\eta)$ is a~$\lmu$-forest of rank $r$ and $0\leq i\leq
r$, we can obtain a~$\lmu$-forest $(T^i,\eta^i)$ of rank $r-1$ by
deleting from $T$ all the vertices of height~$i$ and connecting the
vertices of height~$i+1$ to their grandchildren (unless $i=0$);
$\eta^i$ is the restriction of $\eta$ to $V(T^i)$. In other words,
$(T^i,\eta^i)$ is obtained from $(T,\eta)$ by removing the entire
$i$th level, counting from the leaves, and filling in the gap in
an~obvious way. This allows us to define a~boundary operator by
\begin{equation}\label{eqb}
  \partial(T,\eta)=\sum_{i=0}^r(-1)^i(T^i,\eta^i).
\end{equation}
This paves the way to explicit combinatorial computations by means of
marked forests.

\vskip-10pt
\hskip-100pt $$\epsffile{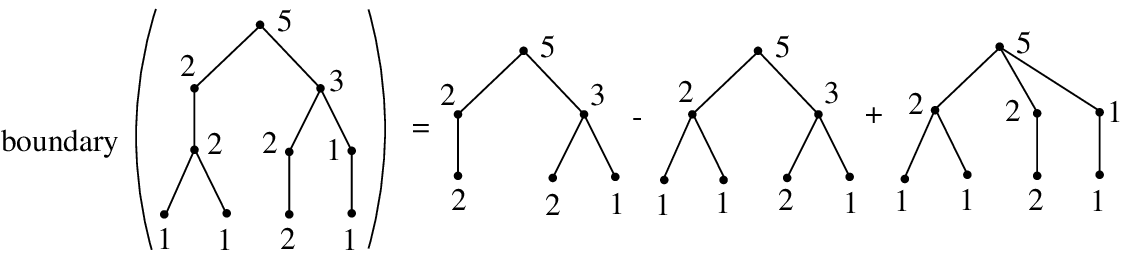}$$
\begin{center}\vskip-2pt
Figure 3.2. An example of a boundary computation.
\end{center}

  For a~given set partition $\pi$ one can define the notion of 
a~{\bf $\pi$-forest $(T,\zeta)$ of rank $r$} almost identically to 
the case of number partitions described above. The difference is that 
$\zeta$ maps $V(T)$ to the set of finite sets, rather than the set of 
natural numbers. The condition~\eqref{mc1} is replaced~by 
\begin{equation}\label{mc2}
\zeta(v)=\bigcup_{w\in\text{\rm{children}}(v)}\zeta(w),
\end{equation}
and $\pi=\{\zeta(v)\,|\,v\in R(T)\}$. For $0\leq i\leq r+1$, analogously to 
$\lam_i(T,\eta)$, we define $\pi_i(T,\zeta)$ to be the set partition which is read
off from the vertices of $T$ having height~$i$. 

  Let $\mu$ be the type of $\pi$, then there exists a~canonical $\mu$-forest 
$(T,|\zeta|)$ associated to each $\pi$-forest $(T,\zeta)$, where $|\zeta|$ 
is obtained as the composition of $\zeta$ with the map which maps finite sets to 
their sizes.

\subsection{The main theorem}

\noindent
Let us describe how to associate a $\lmu$-forest, $\psi(\sigma)$, of
rank $r$ to an $r$-simplex $\sigma$ of $\xlm$. The simplex $\sigma$ is
an~$\st(\pi)$-orbit of $r$-simplices of $\da(\Pi_\lam(\hat 0,\pi))$,
where $\pi$ is a set partition of type~$\mu$. Take a representative of
this orbit, a~chain $c=(x_r>\dots>x_0)$.  Now we define
$\psi(\sigma)=(T,\eta)$. Each element $x_i$ corresponds to the $i$th
level in $T$, counting from the leaves. Each block $b$ of $x_i$
corresponds to a node in the tree; on this node we define the value
of~$\eta$ to be $|b|$. We define the edges of the tree $T$ by
connecting each node corresponding to a block $b$ of $x_i$ to all
nodes corresponding to the blocks of $x_{i-1}$ contained in $b$, we do
that for all $b$ and $i$. The top $(r+1)$th level is added
artificially, its nodes correspond to the blocks of $\pi$, and the
edges from the top level to the $r$th level connect each block of
$\pi$ to the blocks of $x_r$ contained in it.  For example, the value
of $\psi$ on the $\cs_5$-orbit of the chain
$(123)(45)>(123)(4)(5)>(13)(2)(4)(5)$ is the first 5-tree on the
Figure~3.1.

  We are now ready to state and prove the main result of this section.

\begin{thm}\label{main1}
  Assume $\lam\vdash\mu\vdash n$, $\lam\neq\mu$. The correspondence
  $\psi$ of the $r$-simplices of $\xlm$ and $\lmu$-forests $(T,\eta)$
  of rank $r$ is a~bijection. Under this bijection, the boundary
  operator of the triangulated space $\xlm$ corresponds to the
  boundary ope\-ra\-tor described in~\eqref{eqb}.
\end{thm}
  
In particular, the simplices of $\da(\Pi_n)/\cs_n$ along with the cell
inclusion structure are described by the $n$-trees. Indeed,
$\da(\Pi_5)/\cs_5$ is shown in the Figure~3.3.
  
\vskip-5pt
\hskip-100pt $$\epsffile{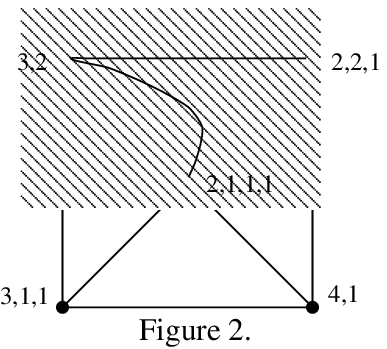}$$
\begin{center}\vskip-5pt
Figure 3.3 
\end{center}

\noindent
The five triangles may be labeled by the five 5-trees of rank 2 in
Figure~3.1.

\vskip4pt

\noindent
{\bf Proof of the Theorem \ref{main1}.} By the definitions of $\Pi_n$
and of $\da$, the $r$-simplices of $\da(\Pi_n)$ can be indexed by
$([n])$-trees of rank~$r$ (we write $([n])$ to emphasize that the set
$[n]$ is viewed here as a~set partition consisting of only one set).
Furthermore, the cell inclusions in $\da(\Pi_n)$ correspond to level
deletion in $([n])$-trees as is described above for the case of number
partitions, because the levels in the $([n])$-trees correspond to the
elements of $\Pi_n$, and the edges in the $([n])$-trees correspond to
block inclusions of two consecutive elements in the chain.

More generally, the $r$-simplices $\sig$ of $\da(\pl(\hat 0,\pi))$ can
be indexed by $\pi$-forests $(T(\sig),\zeta(\sig))$ such that $\lam$
refines the type of $\pi_{0}(T(\sig),\zeta(\sig))$. The definition of
$\psi$ can now be rephrased as {\it associating to $\sigma$ the
  $\lmu$-forest $(T(\sig),|\zeta(\sig)|)$,} where $\mu=\,$type$\,\pi$.

The group action of $\st(\pi)$ on $\da(\pl(\hat 0,\pi))$ corresponds
to relabeling elements within the sets of $\zeta(\sig)$. This shows
that for $g\in\st(\pi)$ we have $T(g\sig)=T(\sig)$ and
$|\zeta(g\sig)|=|\zeta(\sig)|$. Therefore $\psi(\sig)$ is
well-defined, it does not depend on the choice of the representative
of the, corresponding to $\sig$, $\st(\pi)$-orbit of chains.

$\psi$ is surjective, we shall now show that it is also injective.  If
$\sig_1,\sig_2$ are two different $r$-simplices of $\da(\pl(\hat
0,\pi))$ such that $T(\sig_1)=T(\sig_2)$ and $|\zeta(\sig_1)|=
|\zeta(\sig_2)|$, then there exists $g\in\st(\pi)$ such that
$\zeta(g\sig_2)= \zeta(\sig_1)$. Indeed, let $T=T(\sig_1)=T(\sig_2)$
and let $\alpha_1$, resp.~$\alpha_2$, be the string concatenated from
the values of $\zeta(\sig_1)$, resp.~$\zeta(\sig_2)$, on the leaves of
$T$; the order of leaves of $T$ is arbitrary, but the same for
$T(\sig_1)$ and $T(\sig_2)$, the order of elements within each
$\zeta(\sig_1)(v)$, resp.~$\zeta(\sig_2)(v)$, for $v\in L(T)$ is also
chosen arbitrarily. Then $g\in\cs_n$ which maps $\alpha_2$ to
$\alpha_1$ satisfies the necessary conditions:
$\zeta(g\sig_2)=\zeta(\sig_1)$ on the leaves of~$T$, and hence
by~\eqref{mc2} on all vertices of~$T$. Furthermore, since
$g\pi=g\pi_{r+1}(T,\zeta(\sig_2))=\pi_{r+1}(T,\zeta(g\sig_2))=
\pi_{r+1}(T,\zeta(\sig_1))=\pi$, we have $g\in\st(\pi)$.

  This shows that $\psi$ is a~bijection. Since the levels of 
the $(\lam,\mu)$-forests correspond to the $\st(\pi)$-orbits of 
the vertices of $\da(\Pi_\lam(\hat 0,\pi))$ (hence to the 
vertices of $\xlm$), the boundary operator of $\xlm$ corresponds
under $\psi$ to the level deletion in $\lmu$-forests, i.e.,~the
boundary operator described in \eqref{eqb}.
\qed
 
\subsection{Remarks} 

\nin
1. While the presence of the root in an~$([n])$-tree is just
a~formality (two marked $([n])$-trees are equal iff the deletion of
the root gives equal marked forests), the presence of the roots in
a~$\pi$-forest is vital. In fact, if roots were not taken into account
(as seems natural, since the partition read off from the roots does
not correspond to any vertex in $\da(\pl(\hat 0,\pi))$) the argument
above would be false already for vertices: if $\tau_1,\tau_2\in
\pl(\hat 0,\pi)$, such that type$\,(\tau_1)=\,$type$\,(\tau_2)$
(i.e.,~the corresponding marked forests of rank~0 are equal once the
roots are removed), there may not exist $g\in\st(\pi)$, such that
$g\tau_2=\tau_1$ (although such $g\in\cs_n$ certainly exists).
   
\vskip4pt
  
\nin
2. Marked forests equipped with an~order on the children of each
vertex were used by Vassiliev, \cite{Vas94}, to label cells in
a~certain CW-complex structure on the space $\wti\dr^n(m)$, the
one-point compactification of the configuration space of~$m$ unmarked
distinct points in $\dr^n$. Vassiliev's cell decomposition of
$\wti\dr^n(m)$ is a~generalization of the earlier Fuchs' cell
decomposition of $\wti\dr^2(m)=\wti\dc(m)$, \cite{Fuc70}, which
allowed Fuchs to compute the ring $H^*(Br(m),\dz_2)$, where $Br(m)$ is
Artin's braid group on $m$ strings, see also~\cite{Vai78}. Beyond
a~certain similarity of the combinatorial objects used for labeling
the cells, cf.~\cite[Lemma 3.3.1]{Vas94} and Theorem~\ref{main1}, the
connection between the results of this chapter and the results of
Vassiliev and Fuchs seems unclear.

As yet another instance of a~similar situation, we would like to
mention the labeling of the components in the stratification of
$\overline{M}_{0,n}$ (the Deligne-Knudsen-Mumford compactification of
the moduli space of stable projective complex cur\-ves of genus $0$ with
$n$ punctures) with trees with $n$ labeled leaves, see, e.g.,
\cite{FM94,Kn83}.

           \section{A new proof of a~theorem of Arnol'd}   \label{s4}
           
\subsection{Formulation of the main theorem and its corollaries}

\noindent
   In this section we take a~look at a~rather general question of which
$\dq$-acyclicity of the spaces $\xlm$ is a~special case:
\begin{quote}
   Let $\pi\in\Pi_n$ and let $Q$ be an~$\st(\pi)$-invariant subposet of
$\Pi_n(\hat 0,\pi)$. When is the multiplicity of the trivial representation
in the induced representation of $\st(\pi)$ on $\wti H_i(\da(Q),\dq)$ equal 
to~0 for all $i$, in other words, when is $\da(Q)/\st(\pi)$ $\dq$-acyclic?
\end{quote}
 
\begin{df}
  Let $\Lambda$ be a~subset of the set of all number partitions of $n$
  such that $(1^n),(n)\notin\Lambda$. Define $\Pi_\Lambda$ to be the
  subposet of $\Pi_n$ consisting of all set partitions $\pi$ such that
  $(\text{\rm type }\pi)\in\Lambda$.
\end{df}     

Clearly, $\Pi_\Lambda$ is $\cs_n$-invariant and, more generally,
$\Pi_\Lambda(\hat 0,\pi)$ is $\st(\pi)$-invariant. Vice versa, any
$\cs_n$-invariant subposet of $\Pi_n\sm\{(\{1\},\dots,\{n\}),([n])\}$
is of the form $\Pi_\Lambda$ for some~$\Lambda$.

The following theorem is the main result of this section.  The proof
is a~combination of the language of marked forests from
Section~\ref{s3} and the ideas used in the proof
of~\cite[Theorem~4.1]{Ko00a}.

\begin{thm}\label{main2}
  Let $2\leq k<n$. Assume $\Lambda$ is a~subset of the set of all
  number par\-ti\-tions of~$n$ such that $(1^n),(n)\notin\Lambda$ and
  $\Lambda$ satisfies the following condition:
\begin{quote}
  {\rm Condition $C_k$.} If $\mu\in\Lambda$, such that
  $\mu=(\mu_1,\dots,\mu_t)$, where $\mu_i=kq_i+r_i$, $0\leq r_i<k$ for
  $i\in[t]$, then
  $\gamma_k(\mu)=(k^{q_1+\dots+q_t},1^{r_1+\dots+r_t})\in\Lambda$.
\end{quote}
Then for any $\mu\in\Lambda\cup\{(n)\}$ the triangulated space
$\xLm=\da(\Pi_\Lambda(\hat 0,\pi))/\st(\pi)$, where $\mu=\,${\rm
  type}$\,\pi$, is collapsible, in particular the multiplicity of the
trivial representation in the induced $\st(\pi)$-representation on
$\wti H_i(\da(\Pi_\Lambda(\hat 0,\pi)),\dq)$ is equal to~0 for
all~$i$.
\end{thm}

\begin{rem}
  Consider the following special case: $k=2$ and $\Lambda$ is equal to
  the set of all number partitions of~$n$ except for $(1^n)$ and
  $(n)$.  Then the Condition $C_k$ is obviously satisfied. Since in
  this case $\xLm=\da(\Pi_n)/\cs_n$, we conclude that the complex
  $\da(\Pi_n)/\cs_n$ is collapsible.
\end{rem}

Slightly more generally, we have the following result.

\begin{crl}
  Assume $\lam=(k^m,1^{n-km})$, $\lam\vdash\mu$, then $\xlm$ is
  collapsible.  In particular, $\xlm$ is $\dq$-acyclic, therefore
  Theorem~\ref{arn} follows.
\end{crl}
\pr Clearly $\xlm=\xLm$ for
$\Lam=\{\tau\,|\,\lam\vdash\tau\vdash\mu,\tau\neq(n)\}$. It is easy to
check that Condition $C_k$ is satisfied for the case
$\lam=(k^m,1^{n-km})$, therefore Theorem~\ref{arn} follows from
Theorem~\ref{main2} via~\eqref{rswf}. \qed

\vskip4pt

Another consequence of Theorem~\ref{main2} concerns the multiplicity
of the trivial character in certain representations of the symmetric
group. 

Let $\bf k$ be a field such that either char$\,{\bf k}=0$ or
char$\,{\bf k}>n$.  Following a~conjecture of R.\ Stanley, \cite[page
151]{St82}, P.\ Hanlon has proved in~\cite[Theorem 3.1]{Ha83} that if
$\Pi_n^t$ denotes the $\{1,\dots,t\}$ rank selection of the partition
lattice, then the multiplicity of the trivial character in the natural
representation $\cs_n\ra GL(H_i(\da(\Pi_n^t),\bf k))$ induced by the
standard permutation $\cs_n$-representation on the set $[n]$, is equal
to $0$ for all $i$ and $t$.  The following corollary generalizes his
result.

\begin{crl}
  Assume $\Lambda$ is as in the Theorem~\ref{main2}, then the
  multiplicity of the trivial character in the representation
  $\cs_n\ra GL(H_i(\da(\Pi_\Lambda),\bf k))$ is $0$ for all $i$.
\end{crl}
\pr We know that the complex $X_\Lambda=\da(\Pi_\Lambda)/\cs_n$ is
collapsible. By \cite[Theorem 1]{Co56}, $\beta_i(X_\Lambda,{\bf k})$
is equal to the multiplicity of the trivial character in the
representation $\cs_n\ra GL(H_i(\da(\Pi_\Lambda),\bf k))$. Thus the
statement of the corollary is equivalent to saying that $X_\Lambda$ is
${\bf k}$-acyclic, which in turn is immediate from Theorem~\ref{main2}. 
\qed

\vskip4pt


Theorem~\ref{main2} can be viewed as an~attempt to provide a~common
framework for these results in the spirit of the question stated in
the beginning of this section.
    
  \subsection{Auxiliary propositions}

\noindent   
 First we need some terminology. For an arbitrary cell complex $\da$
we denote by $V(\da)$ the set of vertices of $\da$. Assume $\da$ is 
a regular CW complex and $\da'$ is its subcomplex. We denote the set 
of the simplices of $\da$ which are not simplices of $\da'$ by 
$\da\sm\da'$. We use the sign $\succ$ to denote the cover relation 
in the cell structure of~$\da$.

  Assume that, in addition, $\da$ is a triangulated space with some
linear order $\ll$ on the set of vertices. For $\sig\in\da\sm\da'$
we may write $\sig=(x_1,\dots,x_t)$, this notation is slightly
inaccurate since the set of vertices does not determine the simplex 
uniquely, all we mean is that $\sig$ has vertices $x_1\ll\dots\ll x_t$. 
In that case, we let $\xi(\sig)=i$ if $x_1,\dots,x_{i-1}\in V(\da')$ 
and $x_i\notin V(\da')$.

\begin{prop}\label{pr1}
  Let $\da$ be a~regular CW complex and $\da'$ a~subcomplex of $\da$,
  then the following are equivalent:
   
 a) there is a~sequence of collapses leading from $\da$ to $\da'$;
 
 b) there is a~matching of cells of $\da\sm\da'$: $\sig\leftrightarrow
 \phi(\sig)$, such that $\phi(\sig)\succ\sig$ and there is no sequence
 $\sig_1,\dots\sig_t\in\da\sm\da'$ such that $\phi(\sig_1)\succ\sig_2,
 \phi(\sig_2)\succ\sig_3,\dots,\phi(\sig_t)\succ\sig_1$ (such matching
 is called acyclic).
\end{prop}
\pr 

\vskip4pt

\noindent
\underline{a) $\Rightarrow$ b)}: Let the elementary collapses define
the matching $\phi$. Assume there is a~sequence
$\sig_1,\dots,\sig_t\in\da\sm\da'$ such that
$\phi(\sig_1)\succ\sig_2,\phi(\sig_2)\succ\sig_3,\dots,\phi(\sig_t)\succ\sig_1$.
Without loss of generality we can assume that the collapse
$(\sig_1,\phi(\sig_1))$ precedes collapses $(\sig_i,\phi(\sig_i))$ for
$2\leq i\leq t$. Then $\phi(\sig_t)\succ\sig_1$ yields
a~contradiction.

\vskip4pt

\noindent
\underline{b) $\Rightarrow$ a)}: The proof is again very easy, various
versions of it were given in \cite[Corollary 3.5]{Fo98},
\cite[Proposition 3.7]{BBLSW}, and \cite[Theorem 3.2]{Ko00a}. \qed

\begin{prop}\label{pr2}
  Let $\da$ be a~triangulated space with some linear order $\ll$ on
  its set of vertices $V(\da)$. Let $V'\subseteq V(\da)$ and $\da'$ be
  the subcomplex of $\da$ induced by~$V'$. Assume we have a~partition
  $V(\da)=\cup_{z\in V'}V_z$ such that $z=\min_{\ll}V_z$. For
  $\sig\in\da\sm\da'$, let $\chi(\sig)\in V'$ be defined by
  $x_{\xi(\sig)}\in V_{\chi(\sig)}$.  Finally assume that the
  following condition is satisfied:
\begin{quote}
  Condition $\aleph$. If $\sig\in\da\sm\da'$, $\sig=(x_1,\dots,x_t)$, is such that 
either $\xi(\sig)=1$ or $x_{\xi(\sig)-1}\neq\chi(\sig)$, then there exists a~unique 
simplex $\sig'=(x_1,\dots,x_{\xi(\sig)-1},\chi(\sig),x_{\xi(\sig)},\dots,x_t)$ such 
that $\sig'\sm\chi(\sig)=\sig$.
\end{quote}
  Then there is a~sequence of collapses leading from $\da$ to $\da'$.
\end{prop}
\pr Let $U$ denote the set of all $\sig\in\da\sm\da'$, $\sig=(x_1,\dots,x_t)$, 
such that $x_{\xi(\sig)-1}\neq\chi(\sig)$ or $\xi(\sig)=1$. The matching $\phi$ is 
defined by Condition $\aleph$: for $\sig\in U$ we set $\phi(\sig)=\sig'$. By 
Proposition~\ref{pr1} it is enough to check that this matching is acyclic.

For $\sig\in U$ we have $\xi(\phi(\sig))=\xi(\sig)+1$. Moreover, if
$\phi(\sig)\succ\sig'$ and $\sig'\in U$, then $\sig'=\phi(\sig)\sm
x_{\xi(\sig)}$, hence $\xi(\sig')\geq\xi(\phi(\sig))$. Therefore, if
there is a~sequence $\sig_1,\dots,\sig_t\in\da\sm\da'$ such that
$\phi(\sig_1)\succ\sig_2,\phi(\sig_2)
\succ\sig_3,\dots,\phi(\sig_t)\succ\sig_1$, then we have $\xi(\sig_1)<
\xi(\phi(\sig_1))\leq\xi(\sig_2)<\xi(\phi(\sig_2))\leq\dots<
\xi(\phi(\sig_t))\leq \xi(\sig_1)$ which yields a~contradiction. \qed
   
\subsection{Proof of Theorem~\ref{main2}}

\noindent
We define a~$\Lmu$-forest of rank $r$ to be a~marked forest $(T,\eta)$
of rank $r$ such that $\lam_{r+1}(T,\eta)=\mu$ and
$\lam_i(T,\eta)\in\Lam$, for $0\leq i\leq r$.  It follows from the
discussion in Section~\ref{s3} and in particular from
Theorem~\ref{main1} that the $r$-simplices of $\xLm$ can be indexed by
$\Lmu$-forests of rank $r$ so that the boundary relation of $\xLm$
corresponds to level deletion in the marked forests.

We call number partitions of the form $(k^m,1^{n-km})$, for some~$m$,
{\it special}. Let $K$ be the subcomplex of $\xLm$ induced by the set
of all special partitions. We adopt the notations $\xi(\sig)$ and
$\chi(\sig)$ used in Proposition~\ref{pr2} to the context of $\xLm$
and its subcomplex $K$.  The linear order on $V(\xLm)$ can be taken to
be any linear extension of the partial order on $V(\xLm)$ given by the
negative of the length function. The partition of $V(\xLm)$ is given
by: for $v\in V(\xLm)\sm V(K)$, $z\in V(K)$, we have $v\in V_z$ iff
$z=\gamma_k(v)$.

Let us show that the subcomplex $K$ satisfies Condition $\aleph$. Let
$\sig\in\xLm\sm K$, $\sig=(x_1,\dots,x_t)$, and assume $\xi(\sig)=1$
or $\chi(\sig) \neq x_{\xi(\sig)-1}$. In the language of marked
forests this can be reformulated as: $\sig$ is a~$\Lmu$-forest
$(T,\eta)$ of rank $t$ such that $\lam_{\xi(\sig)-1} (T,\eta)$ is not
special and if $\xi(\sig)>1$ then $\lam_0(T,\eta),\dots,$
$\lam_{\xi(\sig)-2}(T,\eta)$ are special, and
$\lam_{\xi(\sig)-2}(T,\eta)\neq \gamma_k(\lam_{\xi(\sig)-1}(T,\eta))$.
In other words, on all vertices of height 0 to $\xi(\sig)-2$ the
function~$\eta$ takes only values 1 or $k$ and for the vertices of
height $\xi(\sig)-1$ it is no longer true. Moreover, there exists
a~vertex of height $\xi(\sig)-1$ which has at least~$k$ children on
which $\eta$ is equal to~1. It is now clear that there exists a~{\it
  unique} $\Lmu$-forest $(\wti T,\ti\eta)$ of rank $r+1$ such that
\begin{itemize}
\item $(\wti T^{\xi(\sig)-1},\ti\eta^{\xi(\sig)-1})=(T,\eta)$;
\item $\lam_{\xi(\sig)-1}(\wti T,\ti\eta)=\gamma_k(\lam_{\xi(\sig)}(\wti T,\ti\eta))$,
i.e.,~$\ti\eta$ takes only values~1 or~$k$ on the vertices of height $\xi(\sig)-1$ 
and each vertex of height $\xi(\sig)$ in $(\wti T,\ti\eta)$ has no more than~$k-1$
children labeled~1.
\end{itemize}
   To construct $(\wti T,\ti\eta)$, extend $(T,\eta)$ by splitting each vertex
of height $\xi(\sig)-1$ into vertices marked~$k$ and~1 so that the number of $k$'s
is maximized. The uniqueness of $(\wti T,\ti\eta)$ follows from the definition
of the notion of isomorphism of marked forests.
 
   We have precisely checked Condition $\aleph$ and therefore by 
Proposition~\ref{pr1} we conclude that there is a~sequence of collapses leading
from $\xLm$ to $K$.

  It remains to see that $K$ is collapsible. If $\mu=(n)$, then $K$ is
a~simplex, so we can assume that $\mu\in\Lam$. If $\mu=\gamma_k(\mu)$,
then $K$ is again a simplex. Otherwise it is easy to see that there is 
a~unique vertex in $\xLm$ labeled $\gamma_k(\mu)$ and that $K$ is 
a~cone with an~apex in this vertex.
\qed

              \section{On a conjecture of Sundaram and Welker}
           \label{s5}    

\subsection{A counterexample to the general conjecture} 
 
\noindent
 The original formulation of Conjecture~\ref{c3} in \cite{SuW97} was
  \begin{conj}\label{c1}\cite[Conjectures 4.12 and 4.13]{SuW97}.
    Let $\lam$ and $\mu$ be different set par\-ti\-tions, such that
    $\lam\vdash\mu$.  Let $\pi\in\pl$ be a~par\-tition of type $\mu$.
    Then the multiplicity of the trivial representation in the
    $\st(\pi)$-module $\wti H_*(\da(\pl(\hat 0,\pi)),\dq)$ is~$0$.
\end{conj}

 In our terms Conjecture \ref{c1} is equivalent to
\begin{conj}\label{c2}
  For $\lam\vdash\mu$, $\lam\neq\mu$, the space $\xlm$ is $\dq$-acyclic.
\end{conj}
               
  We shall give an~example
when $\xlm$ is not even connected. It turns out that if one is only interested
in counting the number of connected components of $\xlm$, then there is a~simpler
poset model which we now proceed to describe. 

\begin{df}\label{pm}
  Assume $\lam\vdash\mu\vdash n$, $\lam\neq\mu$. The $\lmu$-forests of rank $0$ can 
be partially ordered as follows: $(T_1,\eta_1)\prec(T_2,\eta_2)$ if there exists
a~$\lmu$-forest $(T,\eta)$ of rank $1$ such that $(T_1,\eta_1)=(T^1,\eta^1)$ and 
$(T_2,\eta_2)=(T^0,\eta^0)$. We call the obtained poset $P_{\lam,\mu}$.
\end{df}      
In other words, elements of $P_{\lam,\mu}$ are number partitions
$\tau\neq\mu$ such that $\lam\vdash\tau\vdash\mu$, together with
a~bracketing which shows how to form $\mu$ out of $\tau$, the order of
the brackets and of the terms within the brackets is neglected. For
example $(1,1,1)(3,1)(2,2)$ and $(3)(2,1,1)(2,1,1)$ are two different
elements of $P_{(2,1^9),(4^2,3)}$, while $(1,1,1)(2,2)(3,1)$ is equal
to the first mentioned element. These bracketed partitions are ordered
by refinement, preserving the bracket structure.

\begin{prop}\label{pr}
  $X_{\lam,\mu}$ and $\da(P_{\lam,\mu})$ have the same number of connected 
components, i.e., $\beta_0(X_{\lam,\mu})=\beta_0(\da(P_{\lam,\mu}))$.
\end{prop}

\pr We know that $\da(P_{\lam,\mu})$ and $X_{\lam,\mu}$ have the same
set of vertices and that there is an~edge between two vertices $a$ and
$b$ of $\da(P_{\lam,\mu})$ iff $a\prec b$ or $b\prec a$, which is, by
the Definition \ref{pm}, the case iff there is an~edge between the
corresponding vertices of $X_{\lam,\mu}$. This shows that
$\da(P_{\lam,\mu})$ and $\xlm$ have the same number of connected
components. \qed

\vskip4pt

Note that 1-skeleta of $\xlm$ and $\da(P_{\lam,\mu})$ need not be
equal.  $\da(P_{\lam,\mu})$ can intuitively be thought of as
a~simplicial complex obtained by forgetting the multiplicities of
simplices in the triangulated space $\xlm$.

\vskip4pt

\noindent
{\bf Counterexample.} For $n=23$, $\lam=(7,6,4,3,2,1)$,
$\mu=(10,8,5)$, $\xlm$ is disconnected. $P_{\lam,\mu}$ is shown on the
Figure~3.4. Clearly $\da(P_{\lam,\mu})$ is not connected, hence, by
the~Proposition~\ref{pr}, neither is $\xlm$, which disproves
Conjecture~\ref{c1}.
 
\vskip-5pt
\hskip-100pt $$\epsffile{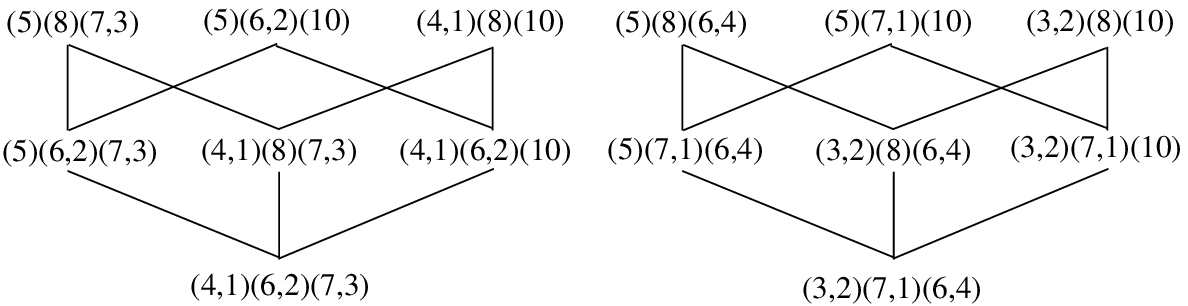}$$
\begin{center}
\vskip-5pt Figure 3.4
\end{center}

\begin{rem}
 In the counterexample above, one can actually verify that 
$\xlm=\da(P_{\lam,\mu})$. However, we choose to use posets $P_{\lam,\mu}$ for
two reasons: 
\begin{enumerate}
\item it is easier to produce series of counterexamples to Sundaram-Welker
conjecture using $\da(P_{\lam,\mu})$ rather than $\xlm$; 
\item we feel that posets $P_{\lam,\mu}$ are of independent interest,
  since they are in a~certain sense the ``naive'' quotient of
  $\Pi_\lam(\hat 0,\pi)$ by $\st(\pi)$.
\end{enumerate}
\end{rem}

 We believe that, in general, connected components of $\xlm$ may be
 not acyclic.

\subsection{Verification of the conjecture in a special case}
  
\begin{df}\label{dgen}
  We say that a~number partition $\lam=(\lam_1,\dots,\lam_t)$ is
  generic (also called free of resonances in~\cite{ShW98}, having no
  equal sub-sums in~\cite{Ko97}) if whenever $\sum_{i\in
    I}\lam_i=\sum_{j\in J}\lam_j$, for some $I,J\subseteq[t]$, we have
  $\{\lam_i\}_{i\in I}=\{\lam_j\}_{j\in J}$ as multisets.
\end{df}

  For example $\lam=(k^m)$ is generic.
  
\begin{prop}\label{tgen}  
  If $\lam$ is generic, then the stratum $\sla$ is homeomorphic to the
  $2l(\lam)$-dimensional sphere.    
\end{prop}
\pr
For a~generic partition
$\lambda=(\underbrace{\lambda_1,\dots,\lambda_1}_{k_1},
\underbrace{\lambda_2,\dots,\lambda_2}_{k_2},\dots,
\underbrace{\lambda_t,\dots,\lambda_t}_{k_t})$ 
we have 
$\sla\cong(S^2)^{(k_1)}\wedge(S^2)^{(k_2)}\wedge\dots\wedge(S^2)^{(k_t)}
\cong S^{2k_1+\dots+2k_t}=S^{2l(\lambda)}$.
\qed


\clearemptydoublepage
\chapter{Complexes of Directed Trees and Their Quotients}

\section{The objects of study and the main questions}

 
To any directed graph one can associate an~abstract simplicial complex
in the~following way.

\begin{df}\label{def}
  Let $G$ be an~arbitrary directed graph. $\dg$ is the~simplicial
  complex $\dg$ constructed as follows: the~vertices of $\dg$ are
  given by the~edges of $G$ and faces are all directed forests which
  are subgraphs of $G$.
\end{df}

  In~\cite{St97} R.\ Stanley asked the~following two questions.

\vskip4pt

{\bf Question 1.} {\it Let $G_n$ be the~complete directed graph on $n$ vertices, 
i.e.,~a~graph having exactly one edge in each direction between any pair of 
vertices, all together $n(n-1)$ edges. The~complex $\dgn$ is obviously pure, but 
is it shellable? }

\vskip4pt

There has been a~recent upsurge of activity in studying the homotopy
type of simplicial complexes constructed from monotone properties of
graphs: the vertices of such a~complex are all possible edges of the
graph and the simplices are given by graphs which satisfy given
monotone property, see~\cite{BBLSW,BWe98}. The question above can be
reformulated in this language: what is the homotopy type of the
simplicial complex corresponding to the monotone property of a~graph
being a~directed forest?

\vskip4pt
{\bf Question 2.} {\it Is the~complex $\dg$ shellable in general? }
\vskip4pt
 
In general, one may ask what are the~homology groups $H_*(\dg)$ and
whether they can be linked to the~combinatorial invariants of
the~graph $G$ in a~simple way.

In this chapter we answer affirmatively the Question~1 in
Theorem~\ref{thm1} and negatively the Question~2 in
Example~~\ref{exam5.2.2}, Section~\ref{s4.2}. Furthermore, we start
the general investigation by computing the~homology groups of $\dg$
for the cases when $G$ is {\it essentially a~tree} (see
Definition~\ref{df4.1}) and when $G$ is a~{\it double directed cycle}.

\vskip4pt 

There is a~natural action of $\cs_n$ on $\da(G_n)$ induced by the
standard permutation action of $\cs_n$ on $[n]$, thus one can form
the~topological quotient $X_n=\da(G_n)/\cs_n$; for example, the case
$n=3$ is shown on the Figure~4.3. The quotient complexes of
combinatorially defined spaces still tend to be rather complicated in
general, however sometimes (as is the case for $\da(\Pi_n)/\cs_n$)
their homology groups are torsion free.

\vskip4pt 
{\bf Question 3.} {\it Are the homology groups $H_i(X_n,{\Bbb Z})$ 
torsion free in general? } 
\vskip4pt

In Section~\ref{s6} we show that the
groups $H_*(X_n,{\Bbb Z})$ are, in general, not free, and also give
a~formula for $\beta_{n-2}(X_n,{\Bbb Q})$.

 \section{First examples and properties} \label{s4.2}

\subsection{Conventions} 

\nin
For brevity, we write $H_*(G)$ instead of $H_*(\dg,{\Bbb Z})$.
We will also use the following standard fact: if $\da$ is a~simplicial
complex and $F_1,\dots,F_t$ is the set of maximal simplices such that
$\da\sm\cup_{i=1}^t F_i$ is contractible, then $\da\simeq \vee_{i=1}^t
S^{\dim F_i}$. In this case, we call
$\{F_1,\dots,F_t\}$ a~{\it generating set} for $\da$. Clearly, for
a~fixed $\da$, there may be several generating sets, but the multiset
$\{\dim F_1,\dots, \dim F_t\}$ is defined uniquely by~$\da$.

Let $\da$ be a simplicial complex, $F$ be one of its maximal simplices
and $\ti F$ be a subsimplex of $F$, such that $\dim F=\dim\ti F+1$ and
$F$ is the only maximal simplex which contains $\ti F$. Then, removing
$F$ and $\ti F$ from $\da$ is called an {\it elementary collapse}.
Obviously, $\da\sm\{\ti F,F\}$ is a~strong deformation retract of
$\da$. For further references on the topological concepts used here we
refer the~reader to the Appendix~D, the textbook by J.\ Munkres,
\cite{Mu84}, and the thorough survey article by A.\ 
Bj\"orner,~\cite{Bj95}.

\subsection{First examples} 

\noindent
Let us give a~couple of examples to illustrate how irregular $\da(G)$
can be.

\begin{exam} 
 {\rm A graph $G$ for which $\dg$ is not pure.}
\end{exam}
  $$\begin{array}{c}
  \epsffile{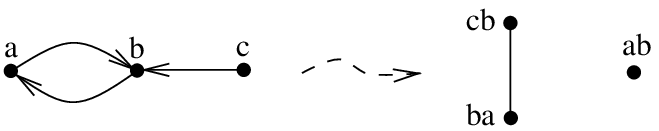}\\
  \text{Figure 4.1}
   \end{array}$$

\begin{exam} \label{exam5.2.2}
{\rm The complex $\dg$ does not have to be shellable.}
  Let $G=C_5$, a~double directed cycle with 5 vertices, see the
  subsection~4.3 for the definition. It is easy to see that $\da(C_5)$
  is a~pure simplicial complex of dimension~3. On the other hand, by
  Proposition~\ref{cyc}, $\da(C_5)\simeq S^2\vee S^3\vee S^3$. This
  implies that $H_2(\da(C_5),{\Bbb Z})={\Bbb Z}$, in particular
  $\da(C_5)$ is not shellable.
\end{exam}

\subsection{Elementary properties} 

\noindent
It is not difficult to derive the~simplest properties of our
construction.  For example $\da(G_1)*\da(G_2)=\da(G_1\uplus G_2)$.
Also it is easy to characterize those graphs $G$ for which $\dg$ is
pure of full dimension.  Namely, for $x\in V(G)$, let
$S(x)=\{y\,|\,(y\ra x)\in E(G)\}$. Then we have
\begin{prop}
  Assume $V(G)=n$. The~following are equivalent
\begin{enumerate}
\item[a)] $\da(G)$ has a~maximal simplex of dimension less than $n-1$;
\item[b)] there exist two disjoint subtrees $T_1$ and $T_2$ of $G$
  with roots $x_1$, resp.~$x_2$ such that $S(x_1)\subseteq V(T_1)$,
  $S(x_2)\subseteq V(T_2)$.  In particular, the sets $S(x_1)$ and
  $S(x_2)$ are disjoint.
\end{enumerate}
  In other words, $\dg$ is pure of dimension $n-1$ iff b) is not true.
\end{prop}
\pr 

\vskip4pt

\noindent
\underline{$a)\Rightarrow b)$}: Let $F$ be a~maximal simplex with
fewer than $n-1$ edges. Then $F$ defines a~forest consisting of two or
more trees. Let $T_1$ and $T_2$ be two different maximal subtrees of
this forest with roots $x_1$, resp. $x_2$. If there exists $y\in
S(x_1)\sm V(T_1)$, then $(y\ra x_1)\cup F$ is a~simplex of $\da(G)$.
This contradicts to the fact that $F$ is a~maximal simplex. Thus
$S(x_1)\subseteq V(T_1)$ and, symmetrically, $S(x_2)\subseteq V(T_2)$.

\vskip4pt

\noindent
\underline{$b)\Rightarrow a)$}: Let $F$ be an arbitrary maximal
simplex of $\da(G)$ such that $E(T_1)\cup E(T_2)\subseteq F$. Since
$S(x_1)\subseteq V(T_1)$ and $S(x_2)\subseteq V(T_2)$, there are no
edges in $F$ which point to either $x_1$ or $x_2$. This proves that
$|F|\leq n-2$. \qed

  \section{Graphs with complete source} \label{sec3}

\subsection{Shellability of complexes of directed trees} 

\nin
In the~next theorem we answer the~first question 
of R.\ Stanley.

\begin{thm} \label{thm1}
  If the~graph $G$ has a~complete source, then the~complex $\dg$ is
  shellable.  In particular $\dgn$ is shellable for any $n\geq 1$.
  More precisely, $\dgn$ is homotopy equivalent to a~wedge of
  $(n-1)^{n-1}$ spheres of dimension $(n-2)$ and the~representatives
  of the~cohomology classes are labeled by the~spanning directed trees
  having vertex $1$ as a~leaf.
\end{thm}

\pr Let $G$ be a~directed graph on $n$ vertices labeled by the~set
$[n]$ and $1$ be the~complete source of $G$. Clearly every partial
subforest of $G$ can be completed to a~tree: just add edges pointing
from the~complete source to the~roots of the~trees of this subforest
(except for the~one which contains the~complete source). Thus $\dg$ is
pure. We shall now describe a~labeling of the~maximal faces of $G$.
 
 For an~edge $(x\ra y)$ we define $\lambda(x\ra y)=x$. For a~graph $H$
 we define $\lambda(H)=(\lambda_1,\dots,\lambda_m)$, where $m=|E(H)|$,
 $\{\lambda_1,\dots,\lambda_m\}=\{\lambda(e)\,|\,e\in E(H)\}$ (as
 multisets) and $\lambda_1\leq\dots\leq\lambda_m$. The~function
 $\lambda$ describes a~labeling of maximal faces of $\dg$ by sequences
 of $n-1$ numbers. Let us order these sequences (and hence the~maximal
 faces) in a~lexicographic order.  We shall next verify that this
 ordering is a~shelling order.

  Let $A$ and $B$ be directed trees on $n$ vertices, such that 
$$\lambda(A)=(\alpha_1,\dots,\alpha_{n-1})\preceq(\beta_1,\dots,\beta_{n-1})=\lambda(B).$$
  Let $C$ be a~graph defined by $V(C)=[n]$,
  $E(C)=E(B)\cap E(A)$. Clearly, $\lambda(C)$ is a~substring of
  $\lambda(A)$ and $\lambda(B)$. $C$ is a~forest, denote its trees by
  $T_1,\dots,T_s$ such that $1\in V(T_1)$ and let $r_i$ be the~root of
  $T_i$. It is clear that edges from $E(A)\sm E(B)$ as well as from
  $E(B)\sm E(A)$ have the~form $(a\ra r_i)$ for some $a\in[n]$, $i\in
  [s]$.

Choose $(x\ra y)\in E(B)\sm E(A)$ such that $x\neq 1$. It must exist,
since otherwise we have $\alpha_i=\beta_i$, $i\in[n-1]$, and all of
the~edges from $E(A)\sm E(B)$ and $E(B)\sm E(A)$ are $(1\ra r_i)$, for
$i=2,\dots,s$, which would imply $A=B$.

Let $B'$ be defined by $E(B')=E(B)\sm\{(x\ra y)\}$. $B'$ consists of
two trees $T'$ and $T''$ with roots $r'$ and $r''$. Assume that $1\in
V(T')$. Define a~new tree $\tilde B$ by $$E(\tilde
B)=E(B')\uplus\{(1\ra r'')\}.$$
It is clear that $\lambda(\tilde
B)=(\lambda(B)\sm\{x\})\uplus\{1\}$ (as multisets), hence
$\lambda(\tilde B)\preceq\lambda(B)$. Furthermore, by construction,
$C$ is a~subgraph of $\tilde B$. Hence we have verified
Con\-di\-ti\-on~(S), thus $\dg$ is shellable.

It is now easy to choose representatives of cohomology classes of
$\dg$. They are labeled by maximal faces $A$, such that
$1\notin\lambda(A)$ (the representatives themselves are functions
which evaluate to 1 on such $A$ and 0 on all other maximal faces).
Furthermore, when $G=G_n$, it is easy to enumerate all such maximal
faces. Denote this number by $f(n)$. The~condition $1\notin\lambda(A)$
simply means that 1 is a~leaf in $A$, so we can instead consider trees
on $n-1$ vertices with some marked vertex. The~number of such trees is
clearly $(n-1)g(n-1)$, where $g(n)$ is the~number of all rooted
labeled trees on $n$ vertices. Since it is well known that
$g(n)=n^{n-1}$ we conclude that $f(n)=(n-1)g(n-1)=(n-1)^{n-1}$. \qed
\vskip4pt

In the~last part of the~argument above we have have found a~new proof
of the combinatorial formula $\tilde\chi(\dgn)=(n-1)^{n-1}$, cf.\ 
\cite{Pi96,St99}.

\subsection{Algebraic consequences of the shelling} 

\noindent
The~result of Theorem \ref{thm1} can be reformulated in the~algebraic
language as follows.
\begin{crl}
  Pick an~index set $\Gamma\subseteq[n]\times[n]$ such that
  $(1,i)\in\Gamma$ for $2\leq i\leq n$, $(i,i)\notin\Gamma$ for $1\leq
  i\leq n$. For an~arbitrary field~$k$, let
  $k[\Gamma]=k[x_{ij}]_{(i,j)\in\Gamma}$. Let $I$ be the~ideal of
  $k[\Gamma]$ generated by the~following monomials:
\begin{itemize}
\item $x_{ij} x_{kj}$, for $(i,j),(k,j)\in\Gamma$;
\item $x_{i_1 i_2}x_{i_2 i_3}\dots x_{i_t i_1}$, for $t\geq 2$, $(i_t,i_1),
(i_j,i_{j+1})\in\Gamma$, for $j\in [t-1]$.
\end{itemize}
  Then the~ring $k[\Gamma]/I$ is Cohen-Macaulay (i.e.,~the~dimension of 
$k[\Gamma]/I$ is equal to its depth, see~\cite{BH93} for further details).
\end{crl}
\pr The~ring $k[\Gamma]/I$ is the~Stanley-Reisner ring of $\da(G)$,
where $G$ is a~directed graph with a~complete source such that
$E(G)=\Gamma$. See~\cite{Re76,St96} for the definition and properties
of Stanley-Reisner rings. \qed 

\vskip4pt

\subsection{A few words on the related $\cs_n$-representations} 

\noindent
The~natural action of the~symmetric group $S_n$ on $V(G_n)$ induces
an~action on $\da(G_n)$ in an~obvious way. This action determines
a~linear representation of $S_n$ on $H_{n-2}(G_n)$. It would be
interesting to understand the~structure of this representation better,
but it seems to be hard. However, R.~Stanley was able to compute its
character,~\cite{St97}. We study the topological quotient 
$\da(\Pi_n)/\cs_n$ in Section~\ref{s6}.

Consider instead a~slightly different, but better behaving action. Let
us act with $S_{n-1}$ on $\da(G_n)$ by permuting the~vertices
$\{2,\dots,n\}$. It follows from the~description of the~cohomology
classes of $\da(G_n)$ that $S_{n-1}$ permutes them and thus
the~representation of $S_{n-1}$ is a~permutation representation:
$S_{n-1}$ permutes double-rooted trees on $n-1$ vertices. Again,
the~structure of this representation, such as decomposition into
irreducibles, is unclear.

\section{Computations for other classes of graphs}

\subsection{Graphs which are essentially trees} 

\begin{df} \label{df4.1}
  A~graph $G$ is called {\bf essentially a~tree} if it turns into an
  undirected tree when one replaces all directed edges/pairs of
  directed edges going in opposite direction by an~edge.
\end{df}

The~following 3 propositions will provide us with a~procedure to
compute homology groups of $\dg$ when $G$ is essentially a~tree.

\begin{prop}\label{t1}
  Let $G$ be a~directed graph and let $x\in V(G)$ with $S(x)=\{y\}$,
  for some $y\in V(G)$. If no edge of $G$ has $x$ as a~source, then
  $\da(G)$ is contractible.
\end{prop}
\pr $\dg$ is a~cone with apex $(y\ra x)$.
\qed

\begin{prop}\label{t2}
  Let $G$ be a~graph, $x\in V(G)$, such that $(x\ra y),(y\ra x)\in
  E(G)$ and there are no other edges where $x$ is a~source or a~sink.
  Then $\dg\simeq\susp\, \da(\tilde G)$, where $\tilde G$ is defined
  by $V(\tilde G)=V(G)\sm\{x\}$, $E(\tilde G)=E(G)\sm(\{(x\la
  y)\}\uplus\{(y\la z)\,|$ $z\in V(G)\})$.
\end{prop}
\pr Let $A={\text{st}\,}(x\ra y)$, $B={\text{st}\,}(y\ra x)$. Then
$\dg=A\cup B$ since every forest not containing edge $(x\ra y)$ can be
extended with the~edge $(y\ra x)$. $A\cap B$ contains those forests
which can be extended with both $(x\ra y)$ and $(y\ra x)$. This means
we have to delete all edges having $x$ or $y$ as a~sink, which gives
$A\cap B=\da(\tilde G)$. Both $A$ and $B$ are contractible, hence, by
\cite[Lemma 10.4(ii)]{Bj95}, $\dg\simeq\susp\,\da(\tilde G)$.  \qed

\begin{prop}\label{t3}
  Let $G$ be a~graph and $x_1,\dots,x_k,y\in V(G)$ such that $(x_1\ra
  y),\dots, (x_k\ra y)\in E(G)$ and there are no further edges which
  have $x_i$ as a~source or a~sink for $i\in[k]$.

   Assume furthermore, that for some $z\in V(G)\sm\{x_1,\dots,x_k,y\}$ 
one of the~following is true:
\be
\item [a)] $(y\ra z)\in E(G)$;
\item [b)] $(z\ra y)\in E(G)$;
\item [c)] $(y\ra z),(z\ra y)\in E(G)$;
\ee
  and there are no other edges having $y$ as a~sink or a~source.

  Then, corresponding to these cases, we have
\be
\item [a)] $\dg\simeq\susp_k\da(\tilde G)$, where $V(\tilde G)=V(G)\sm
  \{x_1,\dots,x_k\}$, $E(\tilde G)=E(G)\sm\{(x_1\ra y),\dots,(x_k\ra
  y),(z\ra y)\}$;
\item [b)] $\dg\simeq\susp_{k+1}\da(G')$, where $G'$ is the subgraph
  of $G$ induced by $V(G)\sm\{x_1,\dots,x_k,y\}$;
\item [c)] If $k=1$, then $\dg\simeq\susp\,\da(G')$; if $k\geq 2$,
  then, for all $t\in\Bbb Z$, we have
  $H_t(G)=H_{t-1}(G')\oplus(H_{t-1}(\tilde G))^{\oplus(k-1)}$.  \ee
\end{prop}
\pr Set $A_i={\text{st}\,}(x_i\ra y)$ for $i\in[k]$. Clearly, all $A_i$ are
contractible and the~intersection of two or more of them is equal to
$\da(\tilde G)$. Furthermore, in case a) we have $\dg=\cup_{i=1}^k
A_i$ and hence the~conclusion of a) follows by~\cite[Lemma
10.4(ii)]{Bj95}.

Assume now that b) or c) holds. In both cases $(z\ra y)\in E(G)$. Set
$B={\text{st}\,}(z\ra y)$, then $\dg=(\cup_{i=1}^k A_i)\cup B$. In the~case b)
the~intersection of any two or more of the~complexes $A_1,\dots,A_k,B$
is equal to $\da(G')$, hence the~conclusion of b) again follows
by~\cite[Lemma 10.4(ii)]{Bj95}.

Let us show c) by induction on $k$. Assume $k=1$, then $\da(G)=A_1\cup
B$.  Both $A_1$ and $B$ are contractible, hence by~\cite[Lemma
10.4(ii)]{Bj95} we get $\da(G)\simeq\susp\,(A_1\cap
B)=\susp\,\da(G')$.

Assume now that $k\geq 2$. Let $A=(\cup_{i=1}^{k-1}A_i)\cup B$ and
$A'=A_{k-1}\cup A_k$. Then $A\cup A'=\da(G)$ and $A\cap A'=A_{k-1}$,
which is a~cone. Thus
$$H_t(G)=H_t(A)\oplus H_t(A')=H_{t-1}(G')\oplus (H_{t-1}(\tilde
G))^{\oplus(k-2)}\oplus H_{t-1}(\tilde G),$$
where the first equality
follows by a~Mayer-Vietoris argument and the second equality follows
from the induction assumption. \qed

\vskip4pt

So, given any graph $G$, which is essentially a~tree, we have
recursive procedure to compute homology groups $H_*(G)$. If some class
of trees which behaves well under recursion is specified, then closed
formulae can be derived if so desired. Observe, for example, that if
$G$ is a~double directed tree with two leaves or more attached to
the~same vertex, then $\dg$ is contractible (just apply Proposition
\ref{t2} and then Proposition \ref{t1}).

\subsection{Double directed strings} 

\noindent
Another interesting specific example is a~double directed string on
$n+1$ vertices $L_n$, which is defined by $V(L_n)=[n+1]$,
$E(L_n)=\{(i\ra i+1), (i+1\ra i)\,|\,i\in[n]\}$. Let $\ti L_n$ be the
directed graph defined by $V(\ti L_n)=[n+2]$, $E(\ti
L_n)=E(L_n)\uplus\{(n+1\la n+2)\}$. The complexes $\da(L_n)$ and
$\da(\ti L_n)$ have an~alternative description.

\begin{df}
  Complex $\cl_n$ has $n$ vertices indexed by the~set $[n]$ and $F\in
  2^{[n]}$ is a~face of $\cl_n$ iff it does not contain $\{i,i+1\}$
  for $i\in [n-1]$. 
\end{df}

   It is easy to see that 
\begin{equation} \label{a}
\cl_{2n}\simeq\da(L_n),\quad\cl_{2n+1}\simeq\da(\ti L_n).
\end{equation}

\newpage

\begin{prop}
$\,$
\begin{enumerate}
\item[(1)] $\cl_{n+3}\simeq\susp\,\cl_n$. 
\item[(2)] The generating simplices for $\cl_n$ are
  $(2,5,\dots,3k+2)$, if $n=3k+2$, and $(2,5,\dots,3k-1)$, if $n=3k$.
\item[(3)] $\cl_1$ is contractible, $\cl_2=S^0$, and $\cl_3\simeq
  S^0$. Hence
\begin{equation} \label{b}
\cl_n\simeq
\begin{cases} S^{k-1},        &\text{ if } n=3k;\\
              \text{a point}, &\text{ if } n=3k+1;\\
              S^k,            &\text{ if } n=3k+2.
\end{cases}
\end{equation}
\end{enumerate}
\end{prop}

\pr Let $C=\cl_{n+3}\sm\{3\}$. Since every maximal simplex of $C$
contains exactly one of the vertices of 1 and 2, we have
$C\simeq\susp\,(\cl_{n+3}\sm\{1,2,3\}) \simeq\susp\,(\cl_n)$. Let us
show that $C$ is a~deformation retract of $\cl_{n+3}$. Order the
simplices of $\cl_{n+3}$ which have vertex~3, but do not have vertex~1
in any order $S_1,\dots,S_k$ respecting inclusions, i.e., if $S_j$ is
a~subcomplex of $S_i$, then $i<j$. Remove pairs of simplices
$(S_i,S_i\uplus \{1\})$ in the increasing order of $i$. These removals
are elementary collapses, they correspond to deformation retracts.
Since $\cl_{n+3}=C\uplus\{S_i\,|\,i\in
[k]\}\uplus\{S_i\uplus\{1\}\,|\,i\in [k]\}$, we conclude that $C$ is
a~deformation retract of $\cl_{n+3}$. This proves (1).

(2) follows from (1) and the fact that if $\da$ is a~simplicial
complex which is homotopy equivalent to a~sphere and $F$ is its
generating simplex, then $F\cup\{a\}$ is a~generating simplex of
$\da*\{a,b\}=\susp\,\da$. To see this, it is enough to show that
$\da*\{a,b\}\sm\{F\cup\{a\}\}$ is contractible. Indeed, removing $F$
and $F\cup\{b\}$ from $\da*\{a,b\}\sm\{F\cup\{a\}\}$ is an~elementary
collapse, and $\da*\{a,b\}\sm\{F,F\cup\{a\},F\cup\{b\}\}=(\da\sm
F)*\{a,b\}$, which is contractible. The verification of (3) is left to
the reader. \qed

\vskip4pt

 Using \eqref{a} and \eqref{b} we immediately derive that
$$\da(L_n)\simeq\begin{cases}
  S^{2k-1}, & \text{ if } n=3k;\\
  S^{2k}, & \text{ if } n=3k+1;\\
  \text{a point}, & \text{ if } n=3k+2.
\end{cases}$$

\subsection{Cycles} \label{sect5} 

\nin
Let $n\geq 3$ and denote by $C_n$ a~double directed cycle, i.e.,
a~directed graph defined by $V(C_n)={\Bbb Z}_n$, $E(C_n)=\{(i\ra
i+1),(i+1\ra i)\,|\, i\in{\Bbb Z}_n\}$. In this section we determine
the homotopy type of $\da(C_n)$.

Again we would like to formulate our complexes in a~slightly different
language.
\begin{df} \label{df2}
  Complex $\cc_n$ has $n$ vertices indexed by the~set ${\Bbb Z}_n$ and
  $F\in 2^{{\Bbb Z}_n}$ is a~face of $\cc_n$ iff it does not contain
  $\{i,i+1\}$ for $i\in{\Bbb Z}_n$.
\end{df} 

Similar to before, we have a~relation $\da(C_n)=\tilde\cc_{2n}$, where
$\ti\cc_{2n}$ is obtained from $\cc_{2n}$ by deleting the simplices
$(1,3,\dots,2n-1)$ and $(2,4,\dots,2n)$.


\begin{prop} \label{cyc}
   The~homotopy type of $\cc_n$ is given by
\begin{equation}\label{hccn} 
\cc_n\simeq\begin{cases}
  S^{k-1}\vee S^{k-1}, & \text{ if } n=3k;\\
  S^{k-1}, & \text{ if } n=3k+1;\\
  S^k, & \text{ if } n=3k+2.
\end{cases}
\end{equation}
Therefore
\begin{equation}\label{hcn}
\da(C_n)\simeq\begin{cases}
  S^{2k-1}\vee S^{2k-1}\vee S^{3k-2}\vee S^{3k-2}, & \text{ if } n=3k;\\
  S^{2k}\vee S^{3k-1}\vee S^{3k-1}, & \text{ if } n=3k+1;\\
  S^{2k}\vee S^{3k}\vee S^{3k}, & \text{ if } n=3k+2.
\end{cases}
\end{equation}
\end{prop}
\pr Let $A={\text{st}\,}_{\cc_n}(1)$ and $B=\cc_n\sm\{1\}$. Then
$A\cup B=\cc_n$, $A\cap B=\cc_n\sm\{1,2,n\}=\cl_{n-3}$, $A$ is
contractible and $B=\cl_{n-1}$.

\vskip3pt

\nin If $n=3k+2$, then $B$ is contractible, hence
$\cc_n\simeq\susp\,(A\cap B)= \susp\,(\cl_{n-3})\simeq S^k$ and
$(1,4,\dots,n-1)$ is the generating simplex.

\vskip3pt

\nin If $n=3k+1$, then $A\cap B$ is contractible, hence $\cc_n\simeq
B=\cl_{n-1}\simeq S^{k-1}$.

\vskip3pt

\nin If $n=3k$, let $F=(3,6,\dots,n)$ be a~generating simplex of $B$.
Since $F\in B\sm(A\cap B)$, one can shrink $B\sm F$ to a~point inside
$\cc_n$ and obtain
$$\cc_n\simeq B\vee\susp\,(A\cap B)\simeq\cl_{3k-1}\vee\cl_{3k}\simeq
S^{k-1}\vee S^{k-1}.$$
Furthermore the generating simplices are
$(3,6,\dots,n)$ and $(1,4,\dots,n-2)$.  This shows \eqref{hccn}.

Let us now see \eqref{hcn}. Let us single out the following simplices
of $\cc_{2n}$: $F_1=(1,3,\dots,2n-1)$, $F_2=(2,4,\dots,2n)$, $\ti
F_1=(1,3,\dots, 2n-3)$, $\ti F_2=(2,4,\dots,2n-2)$. Recall that
$\da(C_n)=\cc_{2n}\sm\{F_1,F_2\}$. $F_1$, resp. $F_2$, is the only
maximal simplex which contains $\ti F_1$, resp.  $\ti F_2$, hence to
remove $F_1$, $F_2$, $\ti F_1$, and $\ti F_2$, from $\cc_{2n}$ means
to perform two elementary collapses. Let
$\hat\cc_{2n}=\cc_{2n}\sm\{F_1,F_2, \ti F_1,\ti F_2\}$. Let us show
that the boundaries of $\ti F_1$ and $\ti F_2$ can be shrunk to a
point within $\hat\cc_{2n}$.

Simplices $F_1$ and $\ti F_1$ lie in the subcomplex
$\cc_{2n}\sm\{2n\}\simeq \cl_{2n-1}$. We can choose
$(2,5,\dots,2n-2)$, if $2n=3k+1$, and $(2,5,\dots, 2n-1)$, if $2n=3k$,
as a~generating simplex of $(\cc_{2n}\sm\{2n\})\sm\{F_1, \ti F_1\}$.
$\ti F_1$ does not contain this generating simplex, hence its boundary
can be shrunk to a point within $(\cc_{2n}\sm\{2n\})\sm\{F_1,\ti
F_1\}$, and hence within $\hat\cc_{2n}$. Analogously the boundary of
$\ti F_2$ can be shrunk to a point within~$\hat\cc_{2n}$.
   
  Thus we obtain $\da(C_n)\simeq\cc_{2n}\sm\{F_1,F_2\}=\hat\cc_{2n}\cup
\{\ti F_1,\ti F_2\}\simeq\hat\cc_{2n}\vee S^{n-2}\vee S^{n-2}\simeq\cc_{2n}
\vee S^{n-2}\vee S^{n-2}$.
\qed

 \section{$\cs_n$-quotients of complexes of directed forests}   
\label{s6}

%
%

\subsection{A~combinatorial description for the cell structure of
  $X_n$} 

\noindent
As mentioned in the introduction, let $X_n$ be the topological
quotient $\da(G_n)/\cs_n$. Clearly, the action of~$\cs_n$ on
$\da(G_n)$ is not free. What is worse, the elements of~$\cs_n$ may fix
the simplices of $\da(G_n)$ without fixing them pointwise: for example
for $n=3$ the permutation $(23)$ ``flips'' the 1-simplex given by the
directed tree $2\longleftarrow 1 \longrightarrow 3$. Therefore, one
does not have a~bijection between the orbits of simplices of
$\da(G_n)$ and simplices of $X_n$. 

To rectify the situation, let us consider the~barycentric subdivision
$B_n=\text{Bsd}\,(\da(G_n))$. We have a~simplicial $\cs_n$-action on
$B_n$ induced from the $\cs_n$-action on $\da(G_n)$ and, clearly,
$B_n/\cs_n$ is homeomorphic to~$X_n$. Furthermore, if an~element of
$\cs_n$ fixes a~simplex of $B_n$ then it fixes it pointwise. In this
situation, it is well-known, e.g.,~see~\cite{Br72}, that the quotient
projection $B_n\ra X_n$ induces a~simplicial structure on $X_n$, in
which simplices of $X_n$ correspond to $\cs_n$-orbits of the simplices
of $B_n$ with appropriate boundary relation.

Let us now give a~combinatorial description of the $\cs_n$-orbits of
the simplices of $B_n$. Let $\sigma$ be a~simplex of $B_n$, then
$\sigma$ is a~chain $(T_1,T_2,\dots,T_{\dim(\sigma)+1})$ of forests on
$n$ labeled vertices, such that $T_i$ is a~subgraph of $T_{i+1}$, for
$i=1,\dots, \dim(\sigma)$. One can view this data in a~slightly
different way: it is a~forest with $\dim(\sigma)+1$ integer labels on
edges (labels on different edges may coincide). Indeed, given a~chain
of forests as above, take the forest $T_{\dim(\sigma)+1}$ and put
label 1 an~all edges of the forest~$T_1$, label 2 on all edges of
$T_2$, which are not labeled yet, etc. Vice versa, given a~forest $T$
with a~labeling, let $T_1$ be the forest consisting of all edges of
$T$ with the smallest label, let $T_2$ be the forest consisting of all
edges of $T$ with one of the two smallest labels, etc. To make the
described correspondence a~bijection, one should identify all labeled
forests on which labelings produce the same order on edges.

\vskip4pt

Formally: {\it the $p$-simplices of $B_n$ are in bijection with the
  set of all pairs $(T,\phi^T)$, where $T$ is a~directed forest on
  $n$~labeled vertices and $\phi^T:E(T)\ra\bbz$, such that
  $|\text{Im}\,\phi^T|=p+1$, modulo the following equivalence
  relation: $(T_1,\phi^{T_1})\sim (T_2,\phi^{T_2})$ if $T_1=T_2$ and
  there exists an~order-preserving injection $\psi:\bbz\ra\bbz$, such
  that $\phi^{T_1}\circ\psi=\phi^{T_2}$.}

\vskip4pt

The boundary operator can be described as follows: for a~$p$-simplex
$(T,\phi^T)$, $p\geq 1$, we have
$$\partial(T,\phi^T)=\sum_{i=1}^{p+1}(-1)^{p+i+1} (T_i,\phi^{T_i}),$$
where, for $i=1,\dots,p$, we have $T_i=T$ and $\phi^{T_i}$ takes the
same values as $\phi^T$ except for the edges on which $\phi^T$ takes
$i$th and $(i+1)$st largest values (say $a$ and $b$), on these edges
$\phi^{T_i}$ takes value $a$. Furthermore, $T_{p+1}$ is obtained from
$T$ be removing the edges with the highest value of $\phi^T$,
$\phi^{T_{p+1}}$ is the restriction of $\phi^T$. Of course, this
description of the boundary map is just a~rephrasing of the deletion
of the $i$th forest from the chain of forests in the original
description. However, we will find it more convenient to work with the
labeled forests rather than the chains of forests.

The orbits of the action of $\cs_n$ can be obtained by forgetting the
numbering of the vertices. Thus, using the fact that simplices
of~$X_n$ and $\cs_n$-orbits of simplices of $B_n$ are the same thing,
we get the following description.

\vskip4pt

{\it The $p$-simplices of~$X_n$ are in bijection with pairs $(T,\phi^T)$,
where $T$ is a~directed forest on $n$~unlabeled vertices and $\phi^T$ is
an~edge labeling of $T$ with $p+1$ labels, modulo a~certain equivalence 
relation. This equivalence relation and the boundary operator are exactly 
as in the description of simplices of $B_n$.}

$$\epsffile{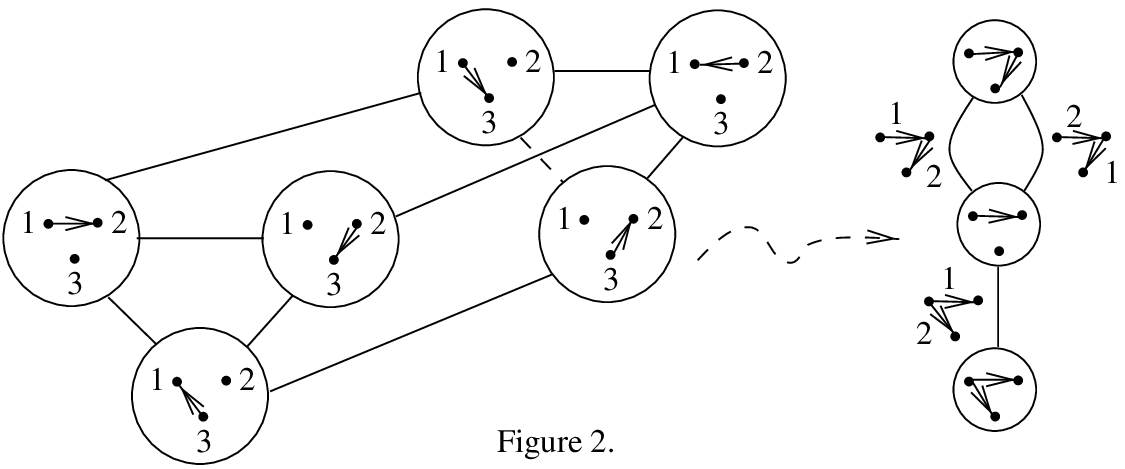}$$
$$\text{Figure 4.2}$$
  
On Figure~4.2 we show the case $n=3$: on the left hand side we have
$\da(G_3)$, on the right hand side is $X_3=\da(G_3)/\cs_3$.  The
labeled forests next to the edges indicate the bijection described
above, labeling on the forests corresponding to the vertices in $X_3$
is omitted. $\cs_3$ acts on $\da(G_3)$ as follows: 3-cycles act as
rotations around the line which goes through the middles of the
triangles, each transposition acts as a~central symmetry on one of the
quadrangles, and as a~``flip'' on the edge which is parallel to that
quadrangle.

\subsection{Filtration and description of the $E^1$ tableau} 

\nin
There is a~natural filtration on the chain complex
associated to the simplicial structure on $X_n$ described above. Let
$F_i$ be the union of all simplices $(T,\phi^T)$ where $T$ has at most
$i$ edges. Clearly, $\emptyset=F_0\subset F_1\subset\dots\subset
F_{n-1}=X_n$.
 
Recall that $E_{p,k}^1=H_p (F_k,F_{k-1})$, here we use the indexing
from the Appendix~D. In other words, the homology is computed with
``truncated'' boundary operator: the last term, where some edges are
deleted from the forest, is omitted. Clearly,
\begin{equation}\label{eq6.1}
   E_{p,k}^1=\bigoplus_{T}H_p(E_T),
\end{equation} 
where the sum is over all forests with $k$ edges and $E_T$ is a~chain 
complex generated by the simplices $(T,\phi^T)$, for various labelings 
$\phi^T$, with the truncated boundary operator as above.

\vskip4pt

Let us now describe a~simplicial complex whose reduced homology
groups, after a~shift in the index by 1, are equal to the nonreduced
homology groups of $E_T$. The arrangement of $k(k-1)/2$ hyperplanes
$x_i=x_j$ in $\bbr^k$ cuts the space $S^{k-1}\cap H$ into simplices,
where $H$ is the~hyperplane given by the equation
$x_1+x_2+\dots+x_k=0$. Denote this simplicial complex $A_k$. The
permutation action of $\cs_k$ on $[k]$ induces an~$\cs_k$-action on
$A_k$. It is easy to see that if an~element of $\cs_k$ fixes a~simplex
of $A_k$, then it fixes it pointwise. Hence, for any subgroup
$\Gamma\subseteq\cs_k$, the $\Gamma$-orbits of the simplices of $A_k$
are in a~natural bijection with the simplices of $A_k/\Gamma$.

\vskip4pt

  Let $T$ be an~arbitrary forest with $n$ vertices and $k$ edges.
Assume that vertices, resp.~edges, are labeled with numbers $1,\dots,n$,
resp.~$1,\dots,k$. $\cs_n$ acts on $[n]$ by permutation, let $\stab(T)$
be stabilizer of $T$ under this action, that is the maximal subgroup of 
$\cs_n$ which fixes $T$. Then $\stab(T)$ acts on $E(T)$, i.e.,~we have
a~homomorphism $\chi:\stab(T)\ra\cs_k$. Let $\cs(T)=\im\chi$. Clearly
$\cs(T)$ does not depend on the choice of the labeling of vertices.
However, relabeling the edges changes $\cs(T)$ to a~conjugate subgroup.
Therefore, for a~forest $T$ without labeling on vertices and edges,
$\cs(T)$ can be defined, but only up to a~conjugation.

\begin{prop}\label{pr6.1} 
  The chain complex of $A_k/\cs(T)$ and $E_T$ (with a~shift by 1 in
  the indexing) are isomorphic. In particular,
  $\rh_p(A_k/\cs(T))=H_{p+1}(E_T)$.
\end{prop} 
\pr Label the $k$ edges of $T$ with numbers $1,\dots,k$. As mentioned
above, the $p$-simplices of $E_T$ are in bijection with labelings of
the edges of~$T$ with numbers $1,\dots,p+1$ (using each number at
least once). Taking in account the chosen labeling of the edges, this
is the same as to divide the set $[k]$ into an~ordered tuple of $p+1$
non-empty sets, modulo the symmetries of $[k]$ induced by the
symmetries of~$T$. Clearly, these symmetries of $[k]$ are precisely
the elements of~$\cs(T)$.

The $(p-1)$-simplices of $A_k$ are in bijection with dividing $[k]$
into an ordered tuple of $p+1$ non-empty sets: by the values of the
coordinates.  Therefore we conclude that the $p$-simplices of $E_T$
are in a~natural bijection with the $(p-1)$-simplices of $A_k/\cs(T)$.
Here the unique 0-simplex of $E_T$, $(T,{\bf 1})$, ($\bf 1$ is the
constant function taking value~1), corresponds in $A_k/\cs(T)$ to the
empty set, which is a~$(-1)$-simplex. One verifies immediately that
the boundary operators of $E_T$ and $A_k/\cs(T)$ commute with the
described bijection. Therefore $E_T$ and $A_k/\cs(T)$ are isomorphic
as chain complexes (after a~shift in the indexing). In particular,
$\rh_p(A_k/\cs(T))=H_{p+1}(E_T)$. \qed
  
\subsection{$\bbq$ coefficients} 

\nin Proposition~\ref{pr6.1} allows us to give a~description of
$E_{*,*}^1$-entries in the case when the homology groups are computed
with rational coefficients. 

Indeed, it is well known that, when a~finite group $\Gamma$ acts on
a~finite simplicial complex $X$, one has
$\rh_i(X/\Gamma,\bbq)=\rh_i^\Gamma(X,\bbq)$, where
$\rh_i^\Gamma(X,\bbq)$ is the maximal vector subspace of
$\rh_i(X,\bbq)$ on which $\Gamma$ acts trivially (more generally
$\bbq$ can be replaced with a~field whose characteristic does not
divide $|\Gamma|$). Since $A_k$ is homeomorphic to $S^{k-2}$ we have
$\rh_{k-2}(A_k,\bbq)=\bbq$ and $\rh_i(A_k,\bbq)=0$ for $i\neq k-2$.

It is easy to compute $\rh_{k-2}^{\cs(T)}(A_k,\bbq)$. In fact, for
$\pi\in\cs_k$, $\alpha\in\rh_{k-2}(A_k,\bbq)$, one has
$\pi(\alpha)=(-1)^{\sgn\pi}\alpha$, where $\sgn$ denotes the sign
homomorphism $\sgn:\cs_k\ra\{-1,1\}$. Therefore
$$\rh_{k-2}(A_k/\cs(T),\bbq)=\rh_{k-2}^{\cs(T)}(A_k,\bbq)=
\begin{cases}
  \bbq,&\text{ if } \cs(T)\subseteq\ca_k,\\
  0, & \text{ otherwise,}
\end{cases}$$
where $\ca_k$ is the alternating group, $\ca_k=\sgn^{-1}(1)$.

Combined with the Proposition~\ref{pr6.1} this gives
$H_i(E_T,\bbq)=\bbq$, if $i=|E(T)|-1$ and
$\cs(T)\subseteq\ca_{|E(T)|}$, and $H_i(E_T,\bbq)=0$ in all other
cases. Therefore it follows from~\eqref{eq6.1} that $\rk
E_{k-1,k}^1=f_{k,n}$, where $f_{k,n}$ is equal to the number of
forests $T$ with $k$ edges and $n$ vertices, such that
$\cs(T)\subseteq\ca_k$.  $\rk E_{p,k}^1=0$ for $p\neq k-1$. Note that
$\beta_i(X_n,\bbq)=0$, for $i\neq n-2$, because
$\beta_i(\da(G_n),\bbq)=0$, for $i\neq n-2$ (by the
Theorem~\ref{thm1}), and
$\beta_i(X_n,\bbq)=\beta_i^{\cs_n}(\da(G_n),\bbq)$.  In particular, by
computing the Euler characteristic of $X_n$ in two different ways, we
obtain

\begin{thm}
  For $n\geq 3$, $\beta_{n-2}(X_n,\bbq)=\sum_{k=2}^{n-1}
  (-1)^{n+k+1}f_{k,n}$.
\end{thm}

  The first values of $f_{k,n}$ are given in the Table~4.3.
Note that there are zeroes on and below the main diagonal and that
the rows stabilize at the entry $(k,2k-1)$ (for $k\geq 2$).

$$\begin{tabular}{ c c c | c c c | c c c | c c c | c c c | c c c | c c c |}
& $k\backslash n$ &&& 1 &&& 2 &&& 3 &&& 4 &&& 5 &&& 6 &\\ \hline 
& 1 &&&   0  &&& {\bf 1} &&& 1 &&&  1 &&&  1 &&&  1 &\\ \hline
& 2 &&&   0  &&& 0 &&& {\bf 1} &&&  1 &&&  1 &&&  1 &\\ \hline
& 3 &&&   0  &&& 0 &&& 0 &&&  2 &&&  {\bf 3} &&&  3 &\\ \hline
& 4 &&&   0  &&& 0 &&& 0 &&&  0 &&&  4 &&&  7 &\\ \hline 
& 5 &&&   0  &&& 0 &&& 0 &&&  0 &&&  0 &&&  8 &\\ \hline
\end{tabular} $$
$$\text{Table 4.3}$$

\subsection{$\bbz$ coefficients} 

\nin The case of integer coefficients is more complicated. In general,
we do not even know the entries of the first tableau. However, we do
know that it is different from the rational case, i.e.,~torsion may
occur.

For example, let $T$ be the~forest with 8 vertices and 6 edges
depicted on Figure~4.4. Clearly, $\cs(T)=\{\text{id},(12)(34)(56)\}$. It
is easy to see that $A_6/\cs(T)$ is a~double suspension (by which we
mean suspension of suspension) of ${\Bbb R}{\Bbb P}^2$, thus the only
nonzero homology group is $\rh_3(A_6/\cs(T),\bbz)=\bbz_2$. In
particular, $E_{4,6}^1$ is not free.

$$\epsffile{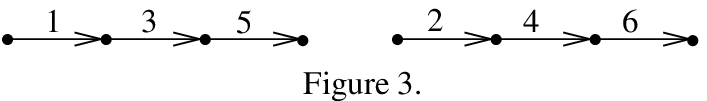}$$
$$\text{Figure 4.4}$$

On the positive side, we can describe the values which $d^1$ takes on
the ``rational'' generators of $E_{*,*}^1$. Let us call a~forest {\it
  admissible} if $\cs(T)\subseteq\ca_{|E(T)|}$. For every admissible
forest $T$ with $k$ edges we fix some order on the edges,
i.e.,~a~bijection $\psi_T:E(T)\ra[k]$. This determines uniquely
an~integer generator $e_T$ of $H_{k-1}(E_T,\bbz)$ by
\begin{equation}\label{eqet}
e_T=\sum_{\cs(T)g}\sgn(g)(T,g\circ\psi_T),
\end{equation}
where we sum over all right cosets of $\cs(T)$, (we choose one
representative for each coset). Observe that the sign of $g$,
resp.~the simplex $(T,g\circ\psi_T)$, are the same for different
representatives of the same right coset class, because
$\cs(T)\subseteq\ca_k$, resp.~by the definition of~$\cs(T)$.

\begin{prop}\label{pr6.2}
  For an~admissible forest $T$, we have
\begin{equation}\label{d1eq}
  d^1(e_T)=\sum_{\alpha}\sgn(\tilde\psi_{T,\alpha}\circ\psi_T^{-1})
  \lambda_{T,\alpha}e_{T\sm\alpha},
\end{equation}
where the sum is over $\cs(T)$-orbits of $E(T)$, for which there
exists a~representative $\alpha$, such that $T\sm\alpha$ is
admissible, we choose one representative for each orbit; note that the
admissibility of $T\sm\alpha$ depends only on the $\cs(T)$-orbit of
$\alpha$, not on the choice of the representative. Notation in the
formula: $T\sm\alpha$ denotes the forest obtained from $T$ by removing
the edge $\alpha$; $\tilde\psi_{T,\alpha}:E(T)\ra[k]$ is defined by
$\tilde\psi_{T,\alpha}|_{T\sm\alpha}=\psi_{T\sm\alpha}$ and
$\tilde\psi_{T,\alpha}(\alpha)=k$;
$\lambda_{T,\alpha}=[\cs(T\sm\alpha): \wti\cs(T)]$, where $\wti\cs(T)$
consists of those permutations of edges of $T\sm\alpha$ which can be
extended to $T$ by fixing the additional edge.
\end{prop}
\pr For an~admissible forest~$T$ with $k$ edges and a~bijection
$\phi:E(T)\ra[k]$, let $(\wti T,\tilde\phi)$ denote a~face simplex of
$(T,\phi)$, where $\wti T$ is obtained from $T$ by removing the edge
with the highest label, $\tilde\phi$ is the restriction of $\phi$ to
$\wti T$. In our notations $(\wti
T,\tilde\phi)=(T\sm\phi^{-1}(k),\phi|_{E(T\sm \phi^{-1}(k))})$.
However, for convenience, we use the notation ``tilde'' in the rest of
the proof.

  According to the general theory for spectral sequences,  
$d^1(e_T)=\partial(e_T)$, where $\partial$ denotes the usual
boundary operator, and we view $\partial(e_T)$ as embedded into the
relative homology group $H_{k-2}(F_{k-1},F_{k-2})$. $\partial(e_T)$ is
a~linear combination of simplices which are obtained from the simplices
$(T,g\circ\psi_T)$ by either merging two labels, or omitting the edge with 
the top label. $e_T\in H_{k-1}(F_k,F_{k-1})$ means that the application
of the ``truncated'' boundary operator to $e_T$ gives 0, therefore all
the simplices obtained by merging two labels will cancel out. 
Furthermore, since $\partial(e_T)\in H_{k-2}(F_{k-1},F_{k-2})$, 
$\dim F_{k-1}=k-2$, and the group $H_{k-2}(F_{k-1},F_{k-2})$ is freely 
generated by $e_U$, where $U$ is an~admissible forest with $k-1$ edges,
we can conclude that also the contributions $(\wti T,\tilde\phi)$, where
$\wti T$ is not admissible, will cancel out. Combining these arguments
with~\eqref{eqet} we obtain:
\begin{equation}\label{raz}
  d^1(e_T)=\sum_{\cs(T)g}\sgn(g)(\wti T,\wti{g\circ\psi_T}),
\end{equation}
where we have only those terms left in the sum, for which $\wti T$ is 
admissible. After regrouping we get
\begin{equation}\label{dva}
\sum_{\cs(T)g}\sgn(g)(\wti T,\wti{g\circ\psi_T})=\sum_{\alpha}
\sum_{\cs(T)g}\sgn(g)(\wti T,\wti{g\circ\psi_T}),
\end{equation}
where in the second term the first sum is taken over all
$\cs(T)$-orbits of~$[k]$, for which $\wti T$ is admissible, while the
second sum is taken over all right cosets $\cs(T)g$ which have
a~representative $g$ such that $g\circ\psi_T(\alpha)=k$, we take one
representative per coset. To verify~\eqref{dva} we just need to
observe that the $\cs(T)$-orbit of $(g\circ\psi_T)^{-1}(k)$ does not
depend on the choice of the representative of $\cs(T)g$; this follows
from the definition of $\cs(T)$.

  Finally, one can see that, for $\alpha$ being an~edge of~$T$, such that
$T\sm\alpha$ is admissible, 
\begin{equation}\label{tri}
\sum_{\cs(T)g}\sgn(g)
(\wti T,\wti{g\circ\psi_T})=\sgn(\tilde\psi_{T,\alpha}
\circ\psi_T^{-1})\lambda_{T,\alpha}\sum_{\cs(T\sm\alpha)h}
\sgn(h)(T\sm\alpha,h\circ\psi_{T\sm\alpha}),
\end{equation}
where the sum in the first term is again taken over all right cosets
$\cs(T)g$ which have a~representative $g$ such that
$g\circ\psi_T(\alpha)=k$, and the sum in the second term is simply
over all right cosets of $\cs(T\sm\alpha)$.

Indeed, on the left hand side we have a~sum over all labelings of
$E(T)$ with numbers $1,\dots,k$, such that $\alpha$ gets a~label~$k$,
and we consider these labelings up to a~symmetry of $T$; each labeling
comes in with a~sign of the permutation $g$, which is obtained by
reading off this labeling in the order prescribed by $\psi_T$. On the
right hand side the same sum is regrouped, using the observation that
to label $E(T)$ with $[k]$, so that $\alpha$ gets a~label~$k$, is the
same as to label $E(T\sm\alpha)$ with $[k-1]$. The only details which
need attention are the multiplicity and the sign.

Every $\cs(T)$-orbit of labelings of $E(T)$ with $[k]$ so that
$\alpha$ gets a~label~$k$ corresponds to
$[\cs(T\sm\alpha):\wti\cs(T)]$ of $\cs(T\sm\alpha)$-orbits of
labelings of $E(T\sm\alpha)$ with $[k-1]$, since we identify labelings
by the actions of different groups:
$\cs(T\sm\alpha)\supseteq\wti\cs(T)$. Each of this
$\cs(T\sm\alpha)$-orbits comes with the same sign, because
$\cs(T\sm\alpha)\subseteq\ca_{k-1}$. The sign
$\sgn(\tilde\psi_{T,\alpha}\circ\psi_T^{-1})$ corresponds to the
change of the order in which we read off the edges: instead of reading
them off according to $\psi_T$, we first read off along
$\psi_{T\sm\alpha}$ and then read off the edge $\alpha$ last.
Formally: $g\circ\psi_T= \tilde h\circ\tilde\psi_{T,\alpha}$, and
$\sgn\tilde h=\sgn h$, hence $\sgn g=\sgn h
\,\,\sgn(\tilde\psi_{T,\alpha}\circ\psi_T^{-1})$, where $\tilde h$ is
defined by $\tilde h|_{[k-1]}=h$, $\tilde h(k)=k$.

Combining \eqref{raz}, \eqref{dva} and \eqref{tri} we obtain
\eqref{d1eq}. \qed

\subsection{Homology groups of $X_n$ for $n=2,3,4,5,6$} 

$X_2$ is just a~point. As shown in Figure~4.3, $X_3\simeq S^1$, where
$\simeq$ denotes homotopy equivalence. With a~bit of labor, one can
manually verify that $X_4\simeq S^2$. Furthermore, one can see that
$H_3(X_5,\bbz)=\bbz^2$ and $\rh_i(X_5,\bbz)=0$ for $i\neq 3$. We leave
this to the reader, while confining ourselves to the case $n=6$. On
Figure~4.5 we have all forests on 6 vertices. We denote some of the
forests by two digits. The numbers over the edges denote the order in
which we read the labels, i.e.,~the bijection $\psi_T$.

 It is easy to see that $A_k/\cs(T)$ is homeomorphic to $S^{k-2}$ for 
all admissible $T$, and is contractible otherwise. The only nontrivial
cases are 41, 47, 48, 51, 55, and 59, all of which can be verified
directly. Therefore, the only nontrivial entries of $E_{*,*}^1$
($\bbz$ coefficients) will lie on the $(k-1,k)$-diagonal. Thus
$H_*(X_6,\bbz)$ can be computed from the chain complex 
$0\lar E_{0,1}^1\stackrel{d^1}{\lla}E_{1,2}^1\stackrel{d^1}{\lla} 
E_{2,3}^1\stackrel{d^1}{\lla} E_{3,4}^1\stackrel{d^1}{\lla} E_{4,5}^1\lar 0$.

  By Proposition~\ref{pr6.2} we have the following relations:
$$\begin{tabular}{ c c  c c c   }
 $d^1(11)=0$, && $d^1(21)=0$, &&  \\ 
 $d^1(31)=2\cdot 21$, && $d^1(32)=21$, && $d^1(33)=21$,  \\  
 $d^1(41)=32-33$, && $d^1(42)=31-32-33$, && $d^1(43)=31-2\cdot 33$,\\  
 $d^1(44)=0$, && $d^1(45)=31-32-33$, && $d^1(46)=32-33$, \\  
  && $d^1(47)=0$, &&  \\ 
\end{tabular} $$
$$\begin{tabular}{ c   c c } 
 $d^1(51)=41-46$, && $d^1(52)=-42+43+44-46$, \\
 $d^1(53)=-42+45$, && $d^1(54)=2\cdot 41+42-43-46+2\cdot 47$, \\
 $d^1(55)=41-43+45$, && $d^1(56)=42+44-45-2\cdot 47$, \\ 
 $d^1(57)=-43+44+45+46$, && $d^1(58)=2\cdot 44+2\cdot 47$, \\  
\end{tabular} $$
here the two-digit strings denote the corresponding forests on
Figure~4.5. Thus $\rh_3(X_6,\bbz)=\bbz_2$, $\rh_4(X_6,\bbz)=\bbz^3$
and $\rh_i(X_6,\bbz)=0$ for $i\neq 3,4$.

\vskip5pt {\it Therefore we conclude that 6 is the smallest value of
  $n$, for which the homology groups $H_*(X_n,\bbz)$ are not free.}

\newpage
  
$$\epsffile{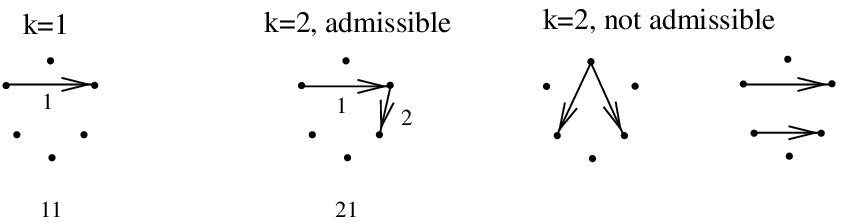}$$
$$\epsffile{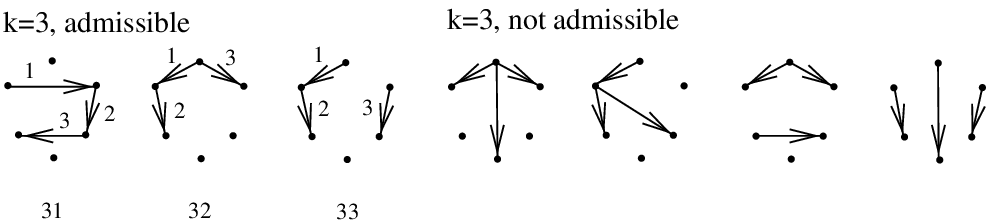}$$
$$\epsffile{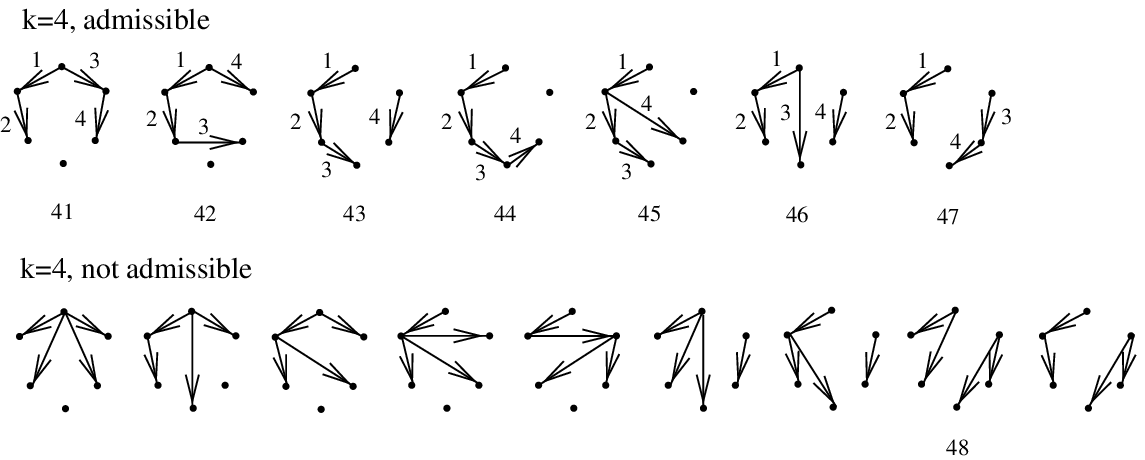}$$
$$\epsffile{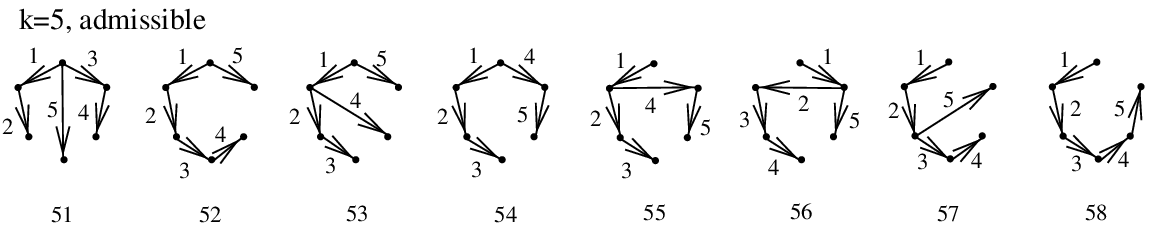}$$
$$\epsffile{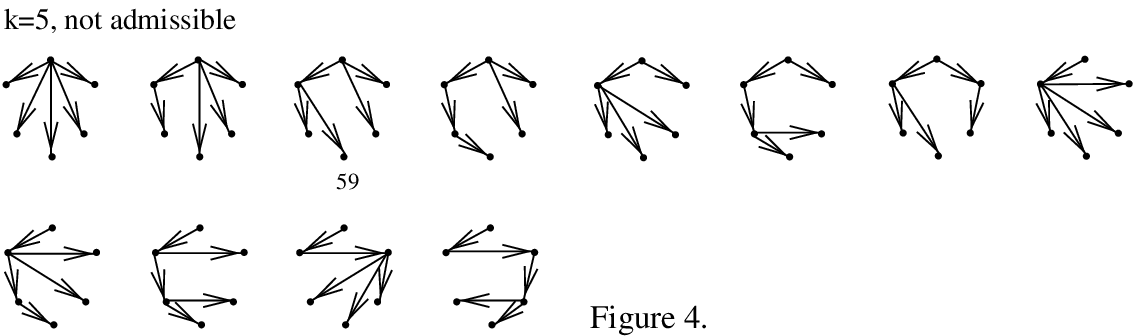}$$
$$\text{Figure 4.5}$$


\clearemptydoublepage
\chapter{Group Actions on Posets}
 
\section{Preamble}

Assume that we have a finite group $G$ acting on a poset $P$ in an
order-preserving way. The~purpose of this chapter is to compare the
various constructions of the quotient, associated with this action.
Our~basic suggestion is to view $P$ as a category and the group action
as a functor from $G$ to $\cat$. Then, it is natural to define $P/G$
to be the colimit of this functor.  As a~result $P/G$ is in general 
a~category, not a~poset.
                 
After getting a hand on the formal setting in Section~\ref{s5.2} we
proceed in Section~\ref{s5.3} with imposing different conditions on
the group action. We give conditions for each of the following
properties to be satisfied:
\begin{enumerate}
\item[(1)] the morphisms of $P/G$ are exactly the orbits of the
  morphisms of~$P$, we call it {\it regularity}; 
\item[(2)] the quotient construction commutes with Quillen's
  nerve functor;
\item[(3)]  $P/G$ is again a poset.
\end{enumerate}
  
Furthermore, we study the class of categories which can be seen as the
``quotient closure'' of the set of all finite posets: loopfree
categories.

\section{Formalization of group actions and the main question}
\label{s5.2}
 
\subsection{Preliminaries} 

\noindent
For a~small category $K$ denote the~set of its objects by $\co(K)$ and
the~set of its morphisms by $\cm(K)$. For every $a\in\co(K)$ there is
exactly one identity morphism which we denote $\id_a$, this allows us
to identify $\co(K)$ with a~subset of $\cm(K)$. If $m$ is a~morphism
of $K$ from $a$ to $b$, we write $m\in\cm_K(a,b)$, $\bo^\bu m=a$ and
$\bo_\bu m=b$. The~morphism $m$ has an~inverse $m^{-1}\in\cm_K(b,a)$,
if $m\circ m^{-1}=\id_a$ and $m^{-1}\circ m =\id_b$. If only
the~identity mor\-phisms have inverses in $K$ then $K$ is said to be
a~category without inverses.
  
We denote the~category of all small categories by $\cat$. If
$K_1,K_2\in\co(\cat)$ we denote by $\cf(K_1,K_2)$ the set of functors
from $K_1$ to $K_2$. We need three full subcategories of $\cat$: ${\bf
  P}$ the~category of posets, (which are categories with at most one
morphism, denoted $(x\ra y)$, between any two objects $x,y$), ${\bf
  L}$ the~category of loopfree categories (see Definition~\ref{5dc}),
and $\grp$ the~category of groups, (which are categories with a single
element, morphisms given by the group elements and the law of
composition given by group multiplication).  Finally, {\bf 1} is
the~terminal object of $\cat$, that is, the~category with one element,
and one (identity) morphism. The~other two categories we use are
$\tp$, the~category of topological spaces, and $\ssc$, the~category of
simplicial sets.

We are also interested in the~functors $\tda:\cat\ra\ssc$ and ${\cal
  R}:\ssc\ra\tp$.  The~composition is denoted $\ti\da:\cat\ra\tp$.
Here, $\da$ is the~nerve functor, see Appendix~B,
or~\cite{Qu73,Qu78,Se68}. In particular, the~simplices of $\tda(K)$
are chains of morphisms in $K$, with dege\-ne\-rate simplices
corresponding to chains that include identity morphisms,
see~\cite{GeM96,We94}. ${\cal R}$ is the~topological realization
functor, see~\cite{Mil57}.


We recall here the definition of a colimit (see~\cite{ML98,Mit65}).  
\begin{df} \label{5col}
  Let $K_1$ and $K_2 $ be categories and
  $X\in\cf(K_1,K_2)$. A~{\bf sink} of $X$ is a~pair consisting of
  $L\in\co(K_2)$, and a~collection of morphisms $\{\lam_s\in\cm_{K_2}
  (X(s),L)\}_{s\in\co(K_1)}$, such that if
  $\alpha\in\cm_{K_1}(s_1,s_2)$ then $\lam_{s_2}\circ
  X(\alpha)=\lam_{s_1}$. (One way to think of this collection of
  morphisms is as a natural transformation between the functors $X$
  and $X'=X_1\circ X_2$, where $X_2$ is the terminal functor
  $X_2:K_1\ra\bf 1$ and $X_1:{\bf 1}\ra K_2$ takes the object of
  $\,\bf 1$ to $L$). When $(L,\{\lam_s\})$ is universal with respect
  to this property we call it the {\bf colimit} of $X$ and write
  $L=\il X$.
\end{df}

\subsection{Definition of the quotient and formulation of 
  the main problem}

\noindent
  Our main object of study is described in the~following definition.

\begin{df}\label{5adf}
  We say that a group $G$ {\bf acts on} a category $K$ if there is a
  functor $\ca_K:G\ra\cat$ which takes the~unique object of G to $K$.
  The~colimit of $\ca_K$ is called the {\bf quotient} of $K$ by the
  action of $G$ and is denoted by $K/G$.
\end{df}

To simplify notations, we identify $\ca_K g$ with $g$ itself.
Furthermore, in Definition~\ref{5adf} the category $\cat$ can be
replaced with any category $C$, then $K,K/G\in\co(C)$. Important
special case is $C=\ssc$. It arises when $K\in\co(\cat)$ and we
consider $\il\tda\circ\ca_K=\tda(K)/G$.

\vskip8pt
\nin
{\bf Main Problem.} $\,$ {\it Understand the relation between the
  topological and the categorical quotients, that is, between
  $\tda(K/G)$ and $\tda(K)/G$.}
\vskip8pt

To start with, by the~universal property of colimits there exists
a~canonical sur\-jec\-tion $\lam:\tda(K)/G\ra\tda(K/G)$. In the~next
section we give com\-bi\-na\-torial con\-di\-tions under which this
map is an isomorphism.

The general theory tells us that if $G$ acts on the~category $K$, then
the~colimit $K/G$ exists, since $\cat$ is cocomplete. We shall now
give an explicit description.

\vskip4pt

\noindent
{\bf An explicit description of the category $K/G$.}

\vskip4pt

\noindent
When $x$ is a morphism of~$K$, denote by $Gx$ the~orbit of $x$ under
the~action of $G$. We have $\co(K/G)=\{Ga\,|\,a\in\co(K)\}$. The
situation with morphisms is more complicated.  Define a relation $\lr$
on the~set $\cm(K)$ by setting $x\lr y$, iff there are decompositions
$x=x_1\circ\dots\circ x_t$ and $y=y_1\circ\dots\circ y_t$ with
$Gy_i=Gx_i$ for all $i\in [t]$. The~relation $\lr$ is reflexive and
symmetric since $G$ has identity and inverses, however it is not in
general transitive.  Let $\sim$ be the~transitive closure of $\lr$, it
is clearly an equivalence relation. Denote the~$\sim$ equivalence
class of $x$ by $[x]$. Note that $\sim$ is the~minimal equivalence
relation on $\cm(K)$ closed under the $G$ action and under
composition; that is, with $a\sim ga$ for any $g\in G$, and if $x\sim
x'$ and $y\sim y'$ and $x\circ x'$ and $y\circ y'$ are defined then
$x\circ x'\sim y \circ y'$. It is not difficult to check that the~set
$\{[x]\,|\,x\in\cm(K)\}$ with the~relations $\bo_\bu[x]=[\bo_\bu x]$,
$\bo^\bu[x]=[\bo^\bu x]$ and $[x]\circ[y]=[x\circ y]$ (whenever
the~composition $x\circ y$ is defined), are the morphisms of
the~category $K/G$.

\vskip4pt

Note that if $P$ is a poset with a $G$ action, the~quotient taken in
$\cat$ need not be a poset, and hence may differ from the~poset
quotient.
\begin{exam} 
  Let $P$ be the center poset in the figure below. Let $\cs_2$ act on
  $P$ by simultaneously permuting $a$ with $b$ and $c$ with $d$. (I)
  shows $P/\cs_2$ in ${\bf P}$ and (II) shows $P/\cs_2$ in $\cat$.
  Note that in this case the quotient in $\cat$ commutes with the
  functor $\da$ (the canonical surjection $\lam$ is an isomorphism),
  whereas the quotient in ${\bf P}$ does not.
\end{exam}

$$\epsffile{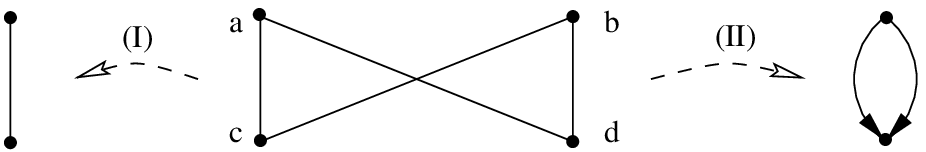}$$
$$\text{Figure 5.1}$$

  \section{Conditions on group actions} \label{s5.3}

\subsection{Outline of the results and surjectiveness of
  the canonical map}
  
\noindent
In this section we consider combinatorial conditions for a~group $G$
acting on a~category $K$ which ensure that the quotient by the group
action commutes with the nerve functor. If $\ca_K:G\ra\cat$ is a~group
action on a~category $K$ then $\tda\circ\ca_K:G \ra \ssc$ is the
associated group action on the nerve of $K$.  It is clear that
$\tda(K/G)$ is a~sink for $\tda\circ\ca_K$, and hence, as previously
mentioned, the universal property of colimits gives a~canonical map 
$\lam:\tda(K)/G\ra\tda(K/G)$. We wish to find conditions under which 
$\lam$ is an isomorphism.
  
First we prove in Proposition~\ref{5sur} that $\lam$ is always
surjective. Furthermore, $Ga=[a]$ for $a\in\co(K)$, which means that,
restricted to $0$-skeleta, $\lam$ is an~iso\-morphism. If the two
simplicial spaces were simplicial complexes (only one face for any
fixed vertex set), this would suffice to show isomorphism. Neither one
is a~simplicial complex in general, but while the quotient of
a~complex $\tda(K)/G$ can have simplices with fairly arbitrary face
sets in common, $\tda(K/G)$ has only one face for any fixed edge set,
since it is a~nerve of a~category. Thus for $\lam$ to be
an~iso\-morphism it is necessary and sufficient to find conditions
under which


\be
\item[{1)}] $\lam$ is an isomorphism restricted to $1$-skeleta;
\item[{2)}] $\tda(K)/G$ has only one face with any given set of edges.
  \ee

We will give conditions equivalent to $\lambda$ being an~isomorphism,
and then give some stronger conditions that are often easier to check,
the strongest of which is also inherited by the action of any subgroup
$H$ of $G$ acting on $K$.

First note that a~simplex of $\tda(K/G)$ is a~sequence
$([m_1],\dots,[m_t])$, $m_i\in\cm(K)$, with $\bo_\bu [m_{i-1}]=\bo^\bu
[m_i]$, which we will call a~{\it chain}. On the other hand a simplex
of $\tda(K)/G$ is an~orbit of a~sequence $(n_1,\dots,n_t)$,
$n_i\in\cm(K)$, with $\bo_\bu n_{i-1}=\bo^\bu n_i$, which we denote
$G(n_1,\dots,n_t)$. The~canonical map $\lam$ is given by
$\lam(G(n_1,\dots,n_t))=([n_1],\dots,[n_t])$.

\begin{prop} \label{5sur}
  Let $K$ be a category and $G$ a group acting on $K$. The canonical map
$\lam:\tda(K)/G\ra\tda(K/G)$ is surjective.  
\end{prop}
\pr By the above description of $\lam$ it suffices to fix a chain
$([m_1],\dots,[m_t])$ and find a chain $(n_1,\dots, n_t)$ with
$[n_i]=[m_i]$. The proof is by induction on $t$. The case $t=1$
is obvious, just take $n_1=m_1$.

Assume now that we have found $n_1,\dots,n_{t-1}$, so that
$[n_i]=[m_i]$, for $i=1,\dots,t-1$, and $n_1,\dots,n_{t-1}$ compose,
i.e., $\bo^\bu n_i=\bo_\bu n_{i+1}$, for $i=1,\dots,t-2$. Since
$[\bo_\bu n_{t-1}]=[\bo_\bu m_{t-1}]=[\bo^\bu m_t]$, we can find $g\in
G$, such that $g\bo^\bu m_t=\bo_\bu n_{t-1}$. If we now take $n_t=g
m_t$, we see that $n_{t-1}$ and $n_t$ compose, and $[n_t]=[m_t]$,
which provides a proof for the induction step. 
\qed


\subsection{Conditions for injectiveness of the canonical projection} 

\begin{df} \label{5df3.2}
  Let $K$ be a category and $G$ a group acting on $K$. We say that
  this action satisfies {\bf Condition (R)} if the following is true:
  If $x,y_a,y_b\in\cm(K)$, $\bo_\bu x=\bo^\bu y_a= \bo^\bu y_b$ and
  $Gy_a=Gy_b$, then $G(x\circ y_a)=G(x\circ y_b)$.  
\end{df}

$$\epsffile{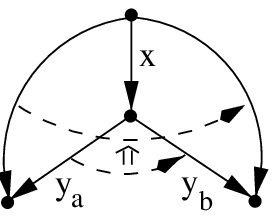}$$
$$\text{Figure 5.2}$$

We say in such case that $G$ acts {\it regularly} on~$K$.

\begin{prop} \label{5reg}
  Let $K$ be a category and $G$ a group acting on $K$. This action
  satisfies Condition (R) iff the~canonical surjection $\lam:\tda(K)/G
  \ra \tda(K/G)$ is injective on $1$-skeleta.
\end{prop}
\pr The injectiveness of $\lam$ on $1$-skeleta is equivalent to
requiring that $Gm=[m]$, for all $m\in\cm(K)$, while Condition~(R) is
equivalent to requiring that $G(m \circ Gn)=G(m \circ n)$, for all
$m,n\in\cm(K)$ with $\bo_\bu m=\bo^\bu n$; here $m\circ Gn$ means the
set of all $m\circ gn$ for which the composition is defined.

Assume that $\lam$ is injective on $1$-skeleta. The we have the
following computation:
$$G(m \circ Gn)=Gm \circ Gn=[m]\circ [n]=[m \circ n]=G(m \circ n),$$
hence the Condition (R) is satisfied.

Reversely, assume that the Condition~(R) is satisfied, that is $G(m
\circ Gn)=G(m \circ n)$.  Since the equivalence class $[m]$ is
generated by $G$ and composition, it suffices to show that orbits are
preserved by com\-po\-sition, which is precisely $G(m \circ Gn)=G(m
\circ n)$. \qed

\vskip4pt

The following theorem is the main result of this chapter. It provides
us with combinatorial conditions which are equivalent to $\lambda$
being an isomorphism.

\begin{thm} \label{5tc}
   Let $K$ be a category and $G$ a group acting on $K$. 
The~following two assertions are equivalent for any $t\geq 2$:
\begin{enumerate}
\item[(1$_t$)] {\bf Condition (C$_{\text{\bf t}}$).}  If
  $m_1,\dots,m_{t-1},m_{a},m_{b}\in\cm(K)$ with $\bo^\bu m_i=\bo_\bu
  m_{i-1}$ for all $2\leq i\leq t-1$, $\bo^\bu m_{a}=\bo^\bu
  m_{b}=\bo_\bu m_{t-1}$, and $Gm_{a}=Gm_{b}$, then there is some
  $g\in G$ such that $gm_{a}=m_{b}$ and $g m_i=m_i$ for $1 \leq i \leq
  t-1$.
\item[(2$_t$)] The~canonical surjection $\lam : \tda(K)/G \ra \tda(K/G)$ is 
injective on $t$-skeleta.  
\end{enumerate}
  In particular, $\lambda$ is an isomorphism iff {\bf (C$_{\text{\bf t}}$)} 
is satisfied for all $t\geq 2$. If this is the case, we say that Condition 
{\bf (C)} is satisfied
\end{thm}

$$\epsffile{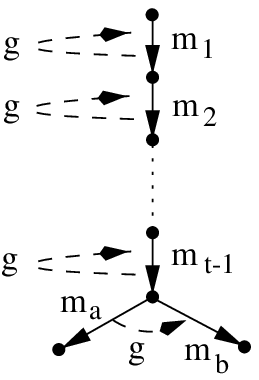}$$
$$\text{Figure 5.3}$$

\pr $(1_t)$ is equivalent to
$G(m_1,\dots,m_t)=G(m_1,\dots,m_{t-1},Gm_t)$; this notation is used,
as before, for all sequences $(m_1,\dots,m_{t-1},gm_t)$ which are
chains, that is for which $m_1\circ\dots\circ m_{t-1}\circ gm_t$ is
defined. $(2_t)$ implies Condition (R) above, and so can be restated
as $G(m_1,\dots,m_t)=(Gm_1,\dots,Gm_t)$.
\begin{multline*}
\underline{(2_t)\Ra (1_t):} \,\,G(m_1,\dots,m_t)=(Gm_1,\dots,Gm_t)=
G(Gm_1,\dots,Gm_t) \\
\supseteq G(m_1,\dots,m_{t-1},Gm_t)\supseteq G(m_1,\dots,m_t).  
\end{multline*}
\vskip-30pt
\begin{multline*}
\underline{(1_2)\Ra (2_2):} \,\,G(m_1,m_2)=G(m_1,Gm_2) \\
=\{g_1(m_1,g_2 m_2)\,|\,\bo_\bu m_1=\bo^\bu g_2 m_2\}\, \\ 
=\{(g_1 m_1,g_2 m_2)\,|\,\bo_\bu g_1 m_1=\bo^\bu g_2 m_2\}=(Gm_1,Gm_2).
\end{multline*}

\hskip-2pt $\underline{(1_t)\Ra (2_t),\,t\geq 3:}$ We use induction on $t$.
\begin{equation*}
\begin{split}
 G(m_1, \dots, m_t)&=G(m_1, \dots, m_{t-1}, Gm_t) \\
&=\{(gm_1,\dots,gm_{t-1},\ti gm_t\,|\,\bo_\bu gm_{t-1}=\bo^\bu \ti gm_t)\} \\
&=\{(g_1 m_1,\dots,g_t m_t)\,|\,\bo_\bu g_i m_i=\bo^\bu g_{i+1}m_{i+1},
i\in[t-1]\} \\
&=(Gm_1,\dots,Gm_t).\,\,\,\,\mqed
\end{split}
\end{equation*}

\begin{exam}
  {\rm A~group action which satisfies Condition (C$_t$), but does not
    satisfy Condition~(C$_{t+1}$).}  Let $P_{t+1}$ be the order sum of
  $t+1$ copies of the 2-element antichain.  The~automorphism group of
  $P_{t+1}$ is the direct product of $t+1$ copies of ${\Bbb Z}_2$.
  Take $G$ to be the index $2$ subgroup consisting of elements with an
  even number of nonidentity terms in the product.
\end{exam}

\nin  The~following condition implies Condition (C), and is often 
easier to check.

\vskip4pt

\nin
{\bf Condition (S).} There exists a set $\{S_m\}_{m\in\cm(K)}$, 
$S_m\subseteq\,\,$Stab$\,(m)$, such that
\begin{enumerate}
\item[(1)] $S_m\subseteq S_{\bo^\bu m}\subseteq S_{m'}$, for any
  $m'\in\cm(K)$, such that $\bo_\bu m'=\bo^\bu m$;
\item[(2)]\vskip-5pt $S_{\bo^\bu m}$ acts transitively on 
$\{g m\,|\,g\in\,$Stab$\,(\bo^\bu m)\}$, for any $m\in\cm(K)$.
\end{enumerate}

\begin{prop} 
 Condition (S) implies Condition (C). 
\end{prop}
\vskip-4pt
\pr Let $m_1,\dots,m_{t-1},m_a,m_b$ and $g$ be as in Condition (C),
then, since $g\in\,$Stab$\,(\bo^\bu m_a)$, there must exist $\ti g\in
S_{\bo^\bu m_a}$ such that $\ti g(m_a)=m_b$. From (1) above one can
conclude that $\ti g(m_i)=m_i$, for $i\in [t-1]$. \qed\vskip4pt

We say that the {\bf strong} Condition (S) is satisfied if Condition
(S) is satisfied with $S_a=\,$Stab$\,(a)$. Clearly, in such a~case
part (2) of the Condition~(S) is obsolete. 

\begin{exam}
  {\rm A~group action satisfying Condition~(S), but not the strong
    Condition~(S).}  Let $K=\cb_n$, lattice of all subsets of $[n]$
  ordered by inclusion, and let $G=\cs_n$ act on $\cb_n$ by permuting
  the ground set $[n]$. Clearly, for $A\subseteq[n]$, we have
  Stab$\,(A)=\cs_A\times\cs_{[n]\sm A}$, where, for $X\subseteq[n]$,
  $\cs_X$ denotes the subgroup of $\cs_n$ which fixes elements of
  $[n]\sm X$ and acts as a permutation group on the set $X$. Since
  $A>B$ means $A\supset B$, condition (1) of (S) is not satisfied for
  $S_A=\,$Stab$\,(A)$: $\cs_A\times\cs_{[n]\sm
    A}\not\supseteq\cs_B\times\cs_{[n]\sm B}$. However, we can set
  $S_A=\cs_A$. It is easy to check that for this choice of
  $\{S_A\}_{A\in\cb_n}$ Condition~(S) is satisfied.
\end{exam}

We close the discussion of the conditions stated above by the
following proposition.

\begin{prop}\label{5int} $\,$
  
\nin  1) The~sets of group actions which satisfy Condition (C) or
  Condition (S) are closed under taking the restriction of the group
  action to a subcategory.

\nin 2) Assume a finite group $G$ acts on a poset $P$, so that Condition
(S) is satisfied. Let $x\in P$ and $S_x\subseteq
H\subseteq\,$Stab$\,(x)$, then Condition (S) is satisfied for the
action of $H$ on $P_{\leq x}$.

\nin 3) Assume a finite group $G$ acts on a category $K$, so that
Condition (S) is satisfied with $S_a=$Stab$(a)$ (strong version), and
$H$ is a subgroup of $G$.  Then the strong version of Condition $(S)$
is again satisfied for the action of $H$ on $K$.
\end{prop}
\vskip-4pt
\pr 1) and 3) are obvious. To show 2) observe that for $a\leq x$ we
have $S_a\subseteq S_x\subseteq H$, hence $S_a\subseteq
H\cap\,$Stab$\,(a)$.  Thus condition (1) remains true. Condition~(2)
is true since
$\{g(b)\,|\,g\in\,$Stab$\,(a)\}\supseteq\{g(b)\,|\,g\in\,$Stab$\,(a)\cap
H\}$. \qed

\subsection{Conditions for the categories to be closed under taking
quotients} 

Next, we are concerned with finding out what categories one may get as
a~quotient of a poset by a group action. In particular, we ask: {\it
  in which cases is the quotient again a poset?} To answer that
question, it is convenient to use the following class of
cate\-go\-ries.

\begin{df}\label{5dc}
  A~category is called {\bf loopfree} if it has no inverses and no
  nonidentity automorphisms.
\end{df}

Intuitively, one may think of loopfree categories as those which can
be drawn so that all nontrivial morphisms point down. To familiarize
us with the notion of a~loopfree category we make the following
observations:
\begin{itemize}
\item $K$ is loopfree iff for any $x,y\in\co(K)$, $x\neq y$, only one
  of the sets $\cm_K(x,y)$ and $\cm_k(y,x)$ is non-empty and
  $\cm_K(x,x)=\{\id_x\}$;
\item a poset is a loopfree category;
\item a barycentric subdivision of an arbitrary category is a loopfree
  category;
\item a barycentric subdivision of a loopfree category is a poset;
\item if $K$ is a loopfree category, then there exists a~partial order
  $\geq$ on the set $\co(K)$ such that $\cm_K(x,y)\neq\emptyset$
  implies $x\geq y$.
\end{itemize}

\begin{df}\label{5hor}
  Suppose $K$ is a~small category, and $T\in\cf(K,K)$. We say that $T$
  is {\bf hori\-zontal} if for any $x\in\co(K)$, if $T(x)\neq x$, then
  $\cm_K(x,T(x))=\cm_K(T(x),x)=\emptyset$. When a~group $G$ acts on
  $K$, we say that the action is hori\-zontal if each $g\in G$ is
  a~horizontal functor.
\end{df}

When $K$ is a~finite loopfree category, the action is always
horizontal. Another example of horizontal actions is given by rank
preserving action on a~(not necessarily finite) poset.  We have the
following useful property:

\begin{prop}\label{5star}
  Let $P$ be a finite loopfree category and $T\in\cf(P,P)$ be a
  hori\-zontal functor.  Let $\ti T\in\cf(\da(P),\da(P))$ be the
  induced functor, i.e., $\ti T=\da(T)$. Then $\da(P_T)=\da(P)_{\ti
    T}$, where $P_T$ denotes the subcategory of $P$ fixed by $T$ and
  $\da(P)_{\ti T}$ denotes the subcomplex of $\da(P)$ fixed by $\ti
  T$.
\end{prop}  

\pr Obviously, $\da(P_T)\subseteq\da(P)_{\ti T}$. On the other hand,
if for some $x\in\da(P)$ we have $\ti T(x)=x$, then the minimal
simplex $\sigma$, which contains $x$, is fixed as a set and, since the
order of simplices is preserved by $T$, $\sigma$ is fixed by $T$
pointwise, thus $x\in\da(P_T)$. \qed\vskip4pt

The~class of loopfree categories can be seen as the closure of the
class of posets under the operation of taking the quotient by a
horizontal group action.  More precisely, we have:
\begin{prop}
  The~quotient of a loopfree category by a~horizontal action is again
  a loopfree category. In particular, the quotient of a poset by
  a~horizontal action is a~loopfree category.
\end{prop}
\pr Let $K$ be a loopfree category and assume $G$ acts on $K$
horizontally. First observe that $\cm_{K/G}([x])=\{\id_{[x]}\}$.
Because if $m\in\cm_{K/G}([x])$, then there exist $x_1,x_2\in\co(K)$,
$\ti m\in\cm_K(x_1,x_2)$, such that $[x_1]=[x_2]$, $[\ti m]=m$. Then
$g x_1=x_2$ for some $g\in G$, hence, since $g$ is a horizontal
functor, $x_1=x_2$ and since $K$ is loopfree we get $\ti m=\id_{x_1}$.

Let us show that for $[x]\neq[y]$ at most one of the sets
$M_{K/G}([x],[y])$ and $M_{K/G}([y],[x])$ is nonempty. Assume the
contrary and pick $m_1\in M_{K/G}([x],[y])$, $m_2\in
M_{K/G}([y],[x])$. Then there exist $x_1,x_2,y_1,y_2\in\co(K)$, $\ti
m_1\in\cm_K(x_1,y_1)$, $\ti m_2\in\cm_K(y_2,x_2)$ such that
$[x_1]=[x_2]=[x]$, $[y_1]=[y_2]=[y]$, $[\ti m_1]=[m_1]$, $[\ti
m_2]=[m_2]$.  Choose $g\in G$ such that $g y_1=y_2$. Then $[g
x_1]=[x_2]=[x]$ and we have $g \ti m_1\in\cm_K(g x_1,y_2)$, so $\ti
m_2\circ g \ti m_1\in \cm_K(g x_1,x_2)$. Since $K$ is loopfree we
conclude that $g x_1=x_2$, but then both $\cm_K(x_2,y_2)$ and
$\cm_K(y_2,x_2)$ are nonempty, which contradicts to the fact that $K$
is loopfree. \qed\vskip4pt

Next, we shall state a~condition under which the quotient of
a~loopfree category is a~poset.

\begin{prop} \label{5sreg}
  Let $K$ be a loopfree category and let $G$ act on $K$. The~following
  two assertions are equivalent:
\begin{enumerate}
\item[(1)] {\bf Condition (SR).} If $x,y\in\cm(K)$, $\bo^\bu x=\bo^\bu
  y$ and $G\bo_\bu x=G\bo_\bu y$, then $Gx=Gy$.
\item[(2)] $G$ acts regularly on $K$ and $K/G$ is a poset.
\end{enumerate}
\end{prop}
\pr $(2)\Ra (1)$. Follows immediately from the regularity of the
action of $G$ and the fact that there must be only one morphism
between $[\bo^\bu x](=[\bo^\bu y])$ and $[\bo_\bu x](=[\bo_\bu y])$.

$$\epsffile{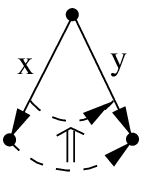}$$
$$\text{Figure 5.4}$$

$(1)\Ra (2)$. Obviously (SR) $\Ra$ (R), hence the action of $G$ is
regular. Furthermore, if $x,y\in\cm(K)$ and there exist $g_1,g_2\in G$
such that $g_1 \bo^\bu x=\bo^\bu y$ and $g_2 \bo_\bu x=\bo_\bu y$,
then we can replace $x$ by $g_1 x$ and reduce the situation to the one
described in Condition (SR), namely that $\bo^\bu x=\bo^\bu y$.
Applying Condition (SR) and acting with $g_1^{-1}$ yields the result. 
\qed

\vskip4pt

 When $K$ is a poset, Condition (SR) can be stated in simpler terms.

\vskip4pt

\nin
{\bf Condition (SRP).} If $a,b,c\in K$, such that $a\geq b$, $a\geq c$
and there exists $g\in G$ such that $g(b)=c$, then there exists $\ti
g\in G$ such that $\ti g(a)=a$ and $\ti g(b)=c$. 

$$\epsffile{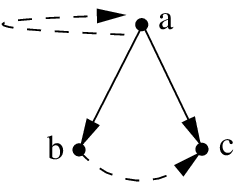}$$
$$\text{Figure 5.5}$$


That is, for any $a,b\in P$, such that $a\geq b$, we require that
the stabilizor of $a$ acts transitively on $Gb$.

\begin{prop}
  Let $P$ be a poset and assume $G$ acts on $P$. The~action of $G$ on
  $P$ induces an action on the barycentric subdivision $\bd P$ (the
  poset of all chains of $P$ ordered by inclusion). This action
  satisfies Condition (S), hence it is regular and $\da(\bd
  P)/G\cong\da((\bd P)/G)$. Moreover, if the action of $G$ on $P$ is
  horizontal, then $(\bd P)/G$ is a poset.
\end{prop}

\pr Let us choose chains $b$, $c$ and $a=(a_1>\dots>a_t)$, such that
$a\geq b$ and $a\geq c$. Then $b=(a_{i_1}>\dots>a_{i_l})$,
$c=(a_{j_1}>\dots>a_{j_l})$. Assume also that there exists $g\in G$
such that $g(a_{i_s})=a_{j_s}$ for $s\in[l]$. If $g$ fixes $a$ then it
fixes every $a_i$, $i\in[t]$, hence $b=c$ and Condition (S) follows.

If, moreover, the action of $G$ is horizontal, then again
$a_{i_s}=a_{j_s}$, for $s\in[l]$, hence $b=c$ and Condition (SRP)
follows. \qed\vskip4pt


\appendix

\clearemptydoublepage
\chapter{Combinatorial Tools}



\section{Number and set partitions}
 
\noindent 
Let $n$ be a~natural number. We denote the set $\{1,\dots,n\}$ by
$[n]$.  

\begin{df}
  A~{\bf number partition} of~$n$ is a~set $\{\lam_1,\dots,\lam_t\}$
  of natural numbers, such that $\lam_1+\dots+\lam_t=n$.
\end{df}

The usual con\-ven\-tion is to write
$\lam=(\lam_1,\dots,\lam_t)$, where $\lam_1\geq\dots\geq\lam_t$, and
$\lam\vdash n$. The {\it length} of $\lam$, denoted $l(\lam)$, is the
number of components of $\lam$, say, in the previous sentence
$l(\lam)=t$. We also use the power notation:
$(n^{\alpha_n},\dots,1^{\alpha_1})=(\underbrace
{n,\dots,n}_{\alpha_n},\dots,\underbrace{1,\dots,1}_{\alpha_1})$.

\begin{df}
We say that $\pi$ is an {\bf ordered set partition} of $[n]$ with $m$
parts (sometimes called {\it blocks}) when $\pi=(\pi_1,\dots,\pi_m)$,
$\pi_i\neq\emptyset$, $[n]=\cup_{i=1}^m\pi_i$, and
$\pi_i\cap\pi_j=\emptyset$, for $i\neq j$. If the order of the parts
is not specified, then $\pi$ is just called a~{\bf set partition}. 
\end{df}

We denote the set of all set partitions, resp.\ ordered set
partitions, of a~set $A$ by $P(A)$, resp.\ $OP(A)$. For $P([n])$,
resp.\ $OP([n])$, we use the shorthand notations $P(n)$, resp.\ 
$OP(n)$. Furthermore, for every set $A$, we let $\un:OP(A)\ra P(A)$ be
the map which takes the ordered partition to the associated unordered
partition.

Whenever we write $\pi\vdash[n]$, it implicitly implies that $\pi$ is
a~set partition, as opposed to a~number partition. A~set partition
$\pi\vdash[n]$, $\pi=(S_1,\dots,S_t)$, is said to have {\it type}
$\lam$, where $\lam\vdash n$ is the~number partition
$\lam=\{|S_1|,\dots,|S_t|\}$.

\begin{df} \mbox{ }

\noindent
(1) For two set partitions $\pi,\ti\pi\vdash S$,
$\pi=(S_1,\dots,S_t)$, $\ti\pi= (\ti S_1,\dots,\ti S_q)$ we write
$\pi\vdash\ti\pi$, and say that $\pi$ {\bf refines} $\tilde\pi$, if
there exists $\iota\vdash[t]$, $\iota=\{I_1,\dots,I_q\}$, such that
$\ti S_i= \cup_{j\in I_i} S_j$, for $i\in[q]$.

\noindent
(2) Analogously, for two number partitions
$\lam=(\lam_1,\dots,\lam_t)$, $\mu=(\mu_1,\dots,\mu_q)$ we write
$\lam\vdash\mu$, and say that $\lam$ {\bf refines} $\mu$, if there
exists $\iota\vdash[t]$, $\iota=\{I_1,\dots,I_q\}$, such that $\mu_i=
\sum_{j\in I_i}\lam_j$, for $i\in[q]$.
\end{df}

\nin Clearly $\pi\vdash[n]$ and $\lam\vdash n$ are special cases of
these notations. Finally observe that if $\pi,\ti\pi$ are two set
partitions, such that $\pi\vdash\ti\pi$, then $(\text{type }\pi)\vdash
(\text{type }\ti\pi)$.

\section{Graphs} 

\begin{df}
  A~{\bf directed graph} $G$ is a~pair of sets $(V(G),E(G))$ such that
  $E(G)\subseteq V(G)\times V(G)\sm\{(x,x)\,|\,x\in V(G)\}$. $V(G)$ is
  called the~set of vertices of $G$, $E(G)$, the~set of edges of $G$.
\end{df}

To support the~intuition, we sometimes write $(x\ra y)$ instead of
$(x,y)$ and call this an {\it edge from} $x$ {\it to} $y$.  A~directed
graph $H=(V(H),E(H))$ is called {\it a~subgraph} of $G$ if
$V(H)\subseteq V(G)$ and $E(H)\subseteq E(G)$. Furthermore, $H$ is
called a {\it subgraph induced by $V(H)$} if $E(H)=E(G)\cap(V(H)\times
V(H))$. For $x,y\in V(G)$, a~{\it directed path} from $x$ to $y$ is an
ordered tuple $((x_1,x_2),(x_2,x_3),\dots,(x_{k-1},x_k))$, such that
$(x_i,x_{i+1})\in E(G)$, for $i\in[k-1]$, and $x_1=x$, $x_k=y$.

\begin{df} $\,$
  
  a) A~directed graph $G$ is called a~{\bf directed tree} with root
  $x\in V(G)$ if for every $y\in V(G)$ there is a~unique directed path
  from $x$ to $y$.
  
  b) A~directed graph $G$ is called a~{\bf directed forest} if there
  exists a~decomposition $V(G)=\uplus_{i\in I}A_i$, (where $\uplus$
  means disjoint union), such that each subgraph induced by $A_i$, for
  $i\in I$, is a~directed tree, and there are no edges between $A_i$
  and $A_j$ for $i\neq j$.
\end{df}
   
If $G$ is a~directed forest, $x\in V(G)$, and there are no edges
$(y\ra x)$, for $y\in V(G)$, we say that $x$ is a {\it root}.

\begin{df}
We say that $x\in V(G)$ is a~{\bf complete source} of $G$ if $(x\ra
y)\in E(G)$ for all $y\in V(G)\sm\{x\}$.
\end{df}


\clearemptydoublepage
\chapter{Posets  and Related Topological Constructions}


\section{Basic notions} 

\noindent 
All posets discussed in this thesis are finite.

\begin{df}
  A poset~$\cl$ is called a {\bf meet-semilattice} if any two elements \linebreak
  $x,y\,{\in}\,\cl$ have a greatest lower bound, i.e., the set
  $\{z\,{\in}\,\cl\,|\,z\,{\leq}\,x, z\,{\leq}\,y\}$ has a~ma\-xi\-mal
  element, called the {\em meet}, $x\,{\wedge}\,y$, of~$x$ and~$y$.
\end{df}

For a~subset $A\,{=}\,\{a_1,\ldots,a_t\}\subseteq\cl$ we let
$\bigwedge A\,{=}\, a_1\,{\wedge}\,\ldots \,{\wedge}\,a_t$ denote the
unique greatest lower bound of $A$, called the~{\em meet}. In
particular, meet-semilattices have a unique minimal element
denoted~$\hat 0$. Minimal elements in $\cl\,{\setminus}\,\{\hat 0\}$
are called the {\em atoms} in~$\cl$.

Symmetrically, meet-semilattices share the following property: for any
subset $A\,{=}\,\{a_1,\ldots,a_t\}\,{\subseteq}\,\cl$ the set
$\{x\,{\in}\,\cl \,|\, x\,{\geq}\, a \text{ for all } a\,{\in}\, A\}$
is either empty or it has a unique minimal element, denoted $\bigvee
A\,{=}\, a_1\,{\vee}\,\ldots \,{\vee}\, a_t$, called the~{\em join}
of~$A$. If the meet-semilattice needs to be specified, we write
$(\bigvee A)_{\cl}\,{=}\, (a_1\,{\vee}\,\ldots \,{\vee}\, a_t)_{\cl}$
for the join of~$A$ in~$\cl$. For brevity, we talk about semilattices
throughout the second chapter, meaning meet-semilattices.

\begin{df}
  For arbitrary posets $P$ and $Q$, the poset $P\times Q$, called
the {\bf direct product}, consists of all pairs $(p,q)$, $p\in P$, $q\in Q$,
ordered by the rule: $(p,q)\leq(\tilde p,\tilde q)$ iff ($p\leq\tilde p$
and $q\leq\tilde q$). 
\end{df}

Let $P$ be an arbitrary poset. For $x\,{\in}\,P$ set: $P_{\leq
  x}\,{=}\,\{y\,{\in}\,P\,{|}\, y\,{\leq}\,x\}$; $P_{<x}$, and
$P_{\geq x}$, $P_{>x}$ are defined analogously. For subsets
$\cg\,{\subseteq}\,P$ with the induced order, and $x\in P$, we define
$\cg_{\leq x}\,{=}\,\{y\,{\in}\,\cg\,{|}\,y\,{\leq}\,x\}$, and
$\cg_{<x}$ again analogously. For intervals in~$P$ we use the
following standard notations: $[x,y]\,{=}\,\{z\,{\in}\, P\,|\,
x\,{\leq}\,z \,{\leq}\,y\}$, $[x,y)\,{=}\,\{z\,{\in}\, P\,|\,
x\,{\leq}\,z \,{<}\,y\}$, etc. We refer to~\cite[Ch.\,3]{St86} for
further details.

\section{Order complexes of posets}

\begin{df}
  For a~poset $P$, let $\da(P)$ denote the {\bf nerve} of~$P$ viewed as
  a~category in the usual way: it is a~simplicial complex with
  $i$-dimensional simplices corresponding to chains of $i+1$ elements
  of $P$ (chains are totally ordered sets of elements of~$P$). In
  particular, vertices of~$\da(P)$ correspond to the elements of~$P$.
  We call $\da(P)$ the {\bf order complex} of~$P$.
\end{df}

The concept of the nerve of a~category goes back at least to D.\ 
Quillen, \cite{Qu73}, and probably even further back to G.\ Segal,
\cite{Se68}. In its combinatorial guise of the order complex, it
appears in the Goresky-MacPherson formula and serves as one of the
main bridges between combinatorics and topology.

\section{Shellability} 

\begin{df}
 A~simplicial complex $\da$ is called {\bf shellable} if there exists
  an~ordering $F_1,\dots,F_t$ of the~ma\-xi\-mal faces of $\da$, such
  that $F_{i+1}\cap (\cup_{j=1}^{i}F_j)$ is a~pure simplicial complex
  of dimension $\dim F_{i+1}-1$ for all $i\in [t-1]$. 
\end{df}

Such an ordering is said to satisfy Condition~(S). Sometimes it is
useful to replace Condition~(S) with an~equivalent Condition~(S$'$):
for $1\leq i<k\leq t$ there exist $j\leq k$ and $x\in F_k$ such that
$F_i\cap F_k\subseteq F_j\cap F_k=F_k\sm\{x\}$.
  
If a~simplicial complex is shellable, then $\da$ is homotopy
equivalent to a~wedge of spheres, indexed by those simplices $F_i$,
for which $F_i\cap (\cup_{j=1}^{i-1}F_j)$ is equal to the full
boundary of $F_i$, and each sphere has the~dimension $\dim F_i$, for
the corresponding~$i$. In particular, the representatives of
cohomology classes are given by the cochains dual to these simplices.
See~\cite{Bj80,Bj95} for more information on shellability.


\clearemptydoublepage

\chapter{Subspace Arrangements}


\section{Definition and related constructions}

\begin{df}
  A set $\ca=\{\ck_1,\dots,\ck_t\}$ of affine linear subspaces in
  a~vector space $\cv$, such that $\ck_i\not\subseteq\ck_j$ for $i\neq
  j$, is called a~{\bf subspace arrangement} in~$\cv$. If all the subspaces
$\ck_1,\dots,\ck_t$ are also required to contain the origin, then 
the subspace arrangement is called {\bf central}.
\end{df}

The topological spaces which one customarily associates to 
a~subspace arrangement are:
\begin{itemize}
\item $V_\ca=\bigcup_{i=1}^t\ck_i$, the union of subspaces;
\item $\cm_\ca=\cv\sm V_\ca$, the {\it complement} of the arrangement.
\end{itemize}

\begin{df}
  Let $\ca$ be a~subspace arrangement in $\bc^n$, and let $G$ be
  a~subgroup of $\GL_n(\dc)$, such that $V_\ca$ is invariant under the
  action of $G$. We say that $G$ {\bf acts on} $\ca$. In that case,
  $\Gamma_\ca^G$ denotes the one-point compactification of~$V_\ca/G$.
\end{df}

\section{Goresky-MacPherson theorem}

\nin
The intersection data of a~subspace arrangement may be represented by
a~poset.

\begin{df}
  To a~subspace arrangement $\ca$ in $\cv$ one can associate
  a~partially ordered set $\cl_\ca$, called the {\bf intersection
    semilattice} of $\ca$. The set if elements of $\cl_\ca$ is
  $\{K\subseteq\dc^n\,|\,\exists\, I\subseteq[t],\text{ such that }
  \bigcap_{i\in I} \ck_i=K\}\cup\{\cv\}$ with the order given by
  reversing inclusions: $x\leq_{\cl_\ca}y$ iff $x\supseteq y$. That
  is, the minimal element of $\cl_\ca$ is~$\cv$, also customarily
  denoted $\hat 0$, and the maximal element is $\bigcap_{K\in\ca}K$.
\end{df} 

The following theorem describes the cohomology groups of the
complement of a~subspace arrangement in terms of the homology groups
of the order complexes of the intervals in the corresponding
intersection lattices.

\begin{thm} (Goresky \& MacPherson, \cite{GoM88}).
  Let $\ca$ be a~central subspace arrangement in $\bc^n$, or in
  $\br^n$, and let $\cl_\ca$ denote its intersection lattice, then
$$\wti H^i(\cm_\ca)\simeq\bigoplus_{x\in\cl_\ca^{\geq \hat 0}}
\wti H_{\text{codim}_{\br}(x)-i-2}(\da(\hat 0,x)).$$
\end{thm}

This provided another strong motivation to study the order complex
construction and to develop technical tools such as lexicographic
shellability, see e.g.,\ ~\cite{Bj94,Ko97}.


\clearemptydoublepage

\chapter{Topological Tools}


\section{Operations on topological spaces} 

\nin For a~topological space $X$, $\wti H_i(X)$, resp.\ $\wti H^i(X)$,
denotes the $i$th reduced homology, resp.\ cohomology, group of $X$;
while $\wti\beta_i(X)$ denotes the $i$th reduced Betti number of~$X$.

Throughout this thesis we use the operations on topological spaces
described in this subsection.

\begin{df}
Let $(X,x)$ and $(Y,y)$ be two pointed topological spaces.

\nin (1) The {\bf wedge} of $X$ and $Y$, denoted $X\vee Y$, is the
pointed topological space $$((X\cup Y)/(x\sim y),z),$$
where the base point $z$ is given by the equivalence class of $x$ 
(and hence of $y$).

\nin (2) The {\bf smash product} of $X$ and $Y$, denoted $X\wedge Y$, is
the pointed topological space obtained as the quotient space
$X\times Y/\sim$, where the equivalence relation is given by:
$(\ti x,y)\sim(x,\ti y)$, for any $\ti x\in X$, $\ti y\in Y$;
with the base point being the equivalence class of $(x,y)$.
\end{df}

The wedge and the smash products enjoy a~variety of properties:
\begin{itemize}
\item they are commutative and associative;
\item $X\vee\text{pt}\,\cong\text{pt}\,\vee X\cong X$;
\item $X\wedge S^0\cong S^0\wedge X\cong X$;
\item $S^n\wedge S^m\cong S^{n+m}$;
\item $X\wedge(Y\vee Z)\cong(X\wedge Y)\vee(X\wedge Z)$.
\end{itemize}

Furthermore, there are important special cases.

\begin{df} 
  Let $X$ be a~topological space, the {\bf suspension} of $X$, denoted
  $\susp X$, is the quotient topological space $X\times[0,1]/\sim$,
  where the equivalence relation is given by $(x_1,0)\sim(x_2,0)$, and
  $(x_1,1)\sim(x_2,1)$, for any $x_1,x_2\in X$.
\end{df}

Clearly, $\susp X\cong X\wedge S^1$. 

\begin{df}
  Let $X$ and $Y$ be two topological spaces. The {\bf join} of $X$ and
  $Y$, denoted $X*Y$, is the quotient topological space $X\times
  [0,1]\times Y/\sim$, where the equivalence relation $\sim$ is given
  by $(x,0,y_1)\sim(x,0,y_2)$, and $(x_1,1,y)\sim(x_2,1,y)$, for any
  $x,x_1,x_2\in X$, and $y,y_1,y_2\in Y$.
\end{df}

For simplicial complexes one can use the alternative, more explicit
definition of the join.

\begin{df}
  Let $X$ and $Y$ be two simplicial complexes, then $Z=X*Y$ is the
  simplicial complex defined by:
\begin{itemize}
\item the set of vertices of $Z$ is equal to the disjoint union of the
  sets of vertices of $X$ and $Y$;
\item the subset $\Sigma$ of the set of vertices of $Z$ is a~simplex
  iff $\Sigma=A\cup B$, where $A$ is a~simplex in $X$ and $B$ is
  a~simplex in $Y$.
\end{itemize}
\end{df}

\section{Spectral sequences} 

\nin A~spectral sequence associated with a~chain complex~$C$ and
a~filtration $F$ on $C$ is a~sequence of 2-dimensional tableaux
$(E_{*,*}^r)_{r=0}^\infty$, where every component $E_{k,i}^r$ is
a~vector space (for simplicity we first consider only field
coefficients), $E_{k,i}^r=0$ unless $k\geq -1$ and $i\geq 0$, and
a~sequence of differential maps $(d^r)_{r=0}^\infty$ such that \be
\item[(0)] $E_{k,i}^0=F_i C_k/F_{i-1}C_k$;
\item[(1)] $d^r\,:\,E_{k,i}^r\longrightarrow E_{k-1,i-r}^r$, 
$\forall\,k,i\in\bbz$;
\item[(2)] $E_{*,*}^{r+1}=H_*(E_{*,*}^r,d^r)$, in other words
\begin{equation}\label{drdef}
E_{k,i}^{r+1}=\text{Ker}\,(E_{k,i}^r\stackrel{d^r}{\longrightarrow}
E_{k-1,i-r}^r)\Big/\text{Im}\,(E_{k+1,i+r}^r\stackrel{d^r}
{\longrightarrow}E_{k,i}^r);
\end{equation}
\item[(3)] for all $k\in\bbz$,
\begin{equation}\label{final}
 H_k(C)=\bigoplus_{i\in\bbz}E_{k,i}^\infty.
\end{equation}
\ee 

\nin{\bf Comments.} 

\vskip4pt
 
\nin
0.\ ~It follows from (0) and (2) that $E_{k,i}^1=H_k(F_i,F_{i-1})$.
 
\vskip4pt
 
\nin
1.\ In the general case $E_{k,i}^\infty$ is defined using the notion
 of convergence of the spectral sequence. We will not explain this
 notion in general, since for the spectral sequence that we consider
 only a~finite number of components in every tableau $E_{*,*}^r$ are
 different from zero, so there exists $N\in\Bbb N$, such that $d^r=0$
 for $r\geq N$. Then, one sets $E_{*,*}^\infty=E_{*,*}^N$, and so
 $H_k(C)=\bigoplus_{i\in\bbz}E_{k,i}^N$.
 
\vskip4pt
 
\nin 2.\ ~For the case of integer coefficients, \eqref{final} becomes
more involved: rather than just summing the entries of
$E^{\infty}_{*,*}$ one needs to solve extension problems to get
$H_*(C)$. This difficulty will not arise in our applications, so we
refer the interested reader to~\cite{McC85} for the detailed
explanation of this phenomena. When considering integer coefficients,
$E^r_{*,*}$ are not vector spaces, but just abelian groups.
 
\vskip4pt
 
\nin
3.\ ~We would like to warn the reader that our indexing is different
 from the standard (but more convenient for our purposes). The
 standard indexing is more convenient for the spectral sequences
 associated to fibrations, an~instance we do not discuss in this
 thesis.

\vskip4pt

Spectral sequences constitute a~convenient tool for computing the
homology groups of a~simplicial complex. A few good sources for 
a~more comprehensive further reading are~\cite{McC85,Sp66,Mas52}.


\clearemptydoublepage

\end{document}